\documentclass[10pt]{article}

\usepackage[bookmarks,bookmarksnumbered,pdfstartview=FitBH,backref=page]{hyperref}
\hypersetup{
colorlinks=true, 
citecolor= blue,
linktoc= all,  
linkcolor= blue, 
urlcolor= green,}

\usepackage[utf8]{inputenc}
\usepackage[T1]{fontenc}

\usepackage{paralist}
\usepackage{mathrsfs} 
\usepackage{amsmath,latexsym,amsthm}
\usepackage{amssymb}
\usepackage{enumitem}
\usepackage{color}

\newtheorem{theo}{Theorem}[section]
\newtheorem{prop}[theo]{Proposition}

\newtheorem{lemm}[theo]{Lemma}

\newtheorem{coro}[theo]{Corollary}
\newtheorem{ques}[theo]{Question}

\theoremstyle{definition}
\newtheorem{rema}[theo]{Remark}
\newtheorem{defi}[theo]{Definition}
\newtheorem{exam}[theo]{Example}

\def\defin#1{\textbf{\emph{#1}}}

\newcommand\Nmath{\mathbb{N}}
\newcommand\Rmath{\mathbb{R}}

\newcommand\Zmath{\mathbb{Z}}

\newcommand\Qmath{\mathbb{Q}}
\newcommand\Cmath{\mathbb{C}}

\newcommand\FF{{\bf F}}

\newcommand{\acting}{\curvearrowright}
\newcommand{\Ima}{\operatorname{Im}}
\newcommand{\Sym}{\operatorname{Sym}}
\newcommand{\cM}{\mathcal{M}}
\newcommand{\cH}{\mathcal{H}}
\newcommand{\coarse}{\mathcal{E}}
\newcommand\RR{{\mathcal{R}}}

\newcommand\LL{{\mathcal{L}}}
\newcommand\LLh{\mathrm{L}}

\newcommand\vraicost{\mathscr{C}}
\newcommand\cost{\mathcal{C}}

\newcommand\ccost{c\mathcal{C}}

\newcommand\cn{\mathcal{N}}

\newcommand{\ch}{\mathcal H}
\newcommand{\Fix}{\mathrm{Fix}}
\newcommand\Stab{\mathrm{Stab}}
\newcommand\id{\mathrm{id}}
\newcommand\SL{\mathrm{SL}}

\newcommand{\lvN}{\mathcal{L}(\cg)}

\newcommand\dom{\operatorname{dom}}
\newcommand\tar{\operatorname{tar}}

\newcommand\uu{\mathfrak{u}}
\newcommand\ul{\mathfrak{u}}
\newcommand\ulae{\mathfrak{u}\textrm{-a.e. }n}
\newcommand{\eps}{\varepsilon}
\newcommand{\ph}{\varphi}
\newcommand{\pmp}{p.m.p.\ }
\newcommand{\rank}{\mathrm{rank}}
\newcommand{\ra}{\mathscr A}
\newcommand{\rb}{\mathscr B}
\newcommand{\cg}{\mathcal G}
\newcommand{\cf}{\mathcal F}
\newcommand{\cd}{\mathcal D}
\newcommand\Lipeq{\overset{\mathcal{L}ip}{\sim}}
\newcommand\ceq{\overset{c}{\sim}}
\makeatletter
\newcommand{\mylabel}[2]{#2\def\@currentlabel{#2}\label{#1}}
\makeatother

\newcommand{\indic}{\mathbf{1}}

\begin{document}

\date{\today}
\title{Non-standard limits of graphs and some orbit equivalence invariants}
\author{Alessandro Carderi, Damien Gaboriau\thanks{C.N.R.S.}, Mikael de la Salle\thanks{C.N.R.S.}}

\maketitle

\begin{abstract}{
We consider probability measure preserving discrete groupoids, group actions and equivalence relations in the context of general probability spaces.
We study for these objects the notions of cost,  $\beta$-invariant and some higher-dimensional variants. We also propose various convergence results about $\ell^2$-Betti numbers and rank gradient for sequences of actions, groupoids or equivalence relations under weak finiteness assumptions.
In particular we connect the combinatorial cost with the cost of the ultralimit equivalence relations. 
Finally a relative version of Stuck-Zimmer property is also considered.}
\end{abstract}



\tableofcontents

\newpage
\section{Introduction}

Limits of discrete structures have attracted a lot of interest in the last decade.
We are interested in sequences of probability measure preserving (p.m.p.) group actions, of graphs or of graphings and in the asymptotic behavior of such invariants as the cost or the $\ell^2$-Betti numbers.
In the course of our investigations we were led to consider graphed \pmp discrete groupoids (in short \pmp groupoids, see Definition~\ref{def: pmp groupoid}) and to introduce the ultraproducts $(\cg_\ul, \Phi_\ul)$ of sequences of such objects. The use of a non principal ultrafilter $\ul$ has the advantage of ensuring convergence
and the existence of limit objects, 
strongly related to more standard objects when some extra assumption provides forms of convergence.
Some other authors have considered a study of sequences of group actions via ultraproducts (see for instance
\cite{Elek-2010-param-test, CKTD,Card-ultra-prod-we-sofic-arxiv,Alekseev-Finn-Sell-sofic-boundaries-2016, Aaserud-Popa}). 
The counterpart in our treatment is that we have been led to develop some orbit equivalence theory in the framework of groupoids on non-standard probability spaces.

\bigskip

The cost of \pmp countable equivalence relations with countable classes is defined in terms of graphings (see \cite{Lev95, Gab-cost}) and is an analogue of the rank of a group.
The extension to \pmp $s$-countable groupoids 
 is quite straightforward and has been considered by several authors (on standard Borel spaces)  for instance by \cite{Ueda-cost-2006, Carderi-Master-thesis-2011, Abert-Nikolov-12, Takimoto-cost-L2-2005}.
The extension to the non-standard setting is also straightforward. 
A graphing $\Phi=(A_i)_{i\in I}$ (Definition~\ref{dfn:graphing}) is a collection of bisections indexed by $I=\Nmath$ or some $\{1, 2, \cdots, N\}$ and the \defin{(groupoid) cost} of the \pmp groupoid $\cg$ is the infimum of the costs of its generating graphings:
\begin{equation}
\cost(\cg)=\inf\left\{\cost(\Phi) \colon \Phi  \text{ generating graphings of } \cg \right\}.
\end{equation} 
 It is the infimum of the ``Haar'' measure of the measurable generating subsets of $\cg$.
When $\cg$ is produced by some \pmp action $\Gamma\acting^{a} (X,\mu)$ of a countable group, one denotes the cost of $\cg$ by $\cost(a)$, not to be confused with the cost of the underlying equivalence relation when the action is not free! 
Another possible source of confusion is the following. A countable group $\Gamma$ is also a \pmp groupoid. As such, its cost is the infimum of the cardinal of a set of generators. In \cite{Gab-cost} and several further references, the cost of a group was defined as the infimum of the costs of all its (standard) free \pmp actions. This latter invariant will subsequently be denoted by $\vraicost_*(\Gamma)$ with a subscript $*$ for the infimum. We observe in Corollary~\ref{cor:monotonicitygroup}(3) that $\vraicost_*(\Gamma)$ is indeed an infimum over all \pmp actions, free or not, standard or not. The  supremum of the costs of all free (standard) \pmp $\Gamma$-actions will be denoted by $\vraicost^*(\Gamma)$.
The group $\Gamma$ has fixed price when $\vraicost_*(\Gamma)=\vraicost^*(\Gamma)$. We introduce the analogous ``fixed price problem'' (Question~\ref{quest: fixed price for groupoids}) for \pmp groupoids.

A key technical result in \cite[\'Etape 1]{Gab-CRAS-cost} and \cite[Proposition IV.35]{Gab-cost} claims that the cost of a \pmp equivalence relation can be computed inside a given Lipschitz class of (finite cost) graphings. The same holds for \pmp groupoids:  $\cost(\cg)=\inf_L\cost_L(\cg,\Phi)$ where  the $L$-\defin{Lipschitz cost} $\cost_L(\cg,\Phi)$
of the graphed groupoid $(\cg,\Phi)$ is the infimum of the costs of all graphings of $\cg$ which are $L$-Lipschitz equivalent to $\Phi$
(see Lemma~\ref{lem:goodgraph} and Lemma~\ref{lem:goodgraphcoarse}).

\medskip
Of course, the speed of convergence of $\cost_L(\cg,\Phi)$ to $\cost(\cg)$ as a function of $L$ depends on the particular graphed groupoid we are considering. This dependance is what makes the spice of the notion of combinatorial cost  for a sequence of graphings (initially introduced by G. Elek \cite{Elek-combinatorial-cost-2007} in the combinatorial context of sequences of finite graphs): one works with a uniform Lipschitz scale $L$ along the sequence (thus obtaining limits of $\cost_L(\cg_n,\Phi)$) and then one let the Lipschitz constant $L$ vary.
If $(\cg_n,\Phi_n)_n$ is a sequence of graphed \pmp groupoids of bounded size, we define the $\ul$-\defin{combinatorial cost} (Definition~\ref{def: ul-comb-cost}) of the sequence as
$$
\ccost_\ul((\cg_n,\Phi_n)_n):=\inf_L\lim_{n\in\ul}\cost_L(\cg_n,\Phi_n).
$$

The Lipschitz framework is appropriate for uniformly bounded sequences of graphings (Section~\ref{subsection:finite_graphings}). 
For the general unbounded setting, we introduce the more appropriate notions of coarse structure (Section~\ref{subsect : coarse structures and graphings}) and coarse equivalence (Section~\ref{subsection:infinite_graphings}) for \pmp groupoids, inspired by Roe's topological analogues \cite{Roe-2003} and we define the $\ul$-combinatorial cost in that context (Definition~\ref{def: comb cost, infinite size}).

For group actions, our first main result is the following:
\begin{theo}[{Theorem~\ref{th: equality of lim cost for fixed price}}]
\label{th-intro:equality of lim cost for fixed price}
  Let $\Gamma$ be a finitely generated group. For every sequence $(a_n)_n$ of \pmp $\Gamma$-actions we have:
  \begin{equation}
\cost(a_\ul)=\ccost_\ul((a_n)_n) \geq \lim_{\ul}\cost(a_{n})\geq \liminf_{n\to \infty}\cost(a_{n}) \geq \vraicost_*(\Gamma).
\end{equation}
Moreover if the action $a_\ul$ is essentially free we have $\vraicost^*(\Gamma)\geq \cost(a_\ul)$ and hence if $\Gamma$ has fixed price we obtain 
 $\lim_\ul\cost(a_n)=\ccost_\ul((a_n)_n)=\cost(a_{\ul})=\vraicost_*(\Gamma)$.
\end{theo}

We observe that this continuity statement fails when one removes the finite generation assumption (see Remark~\ref{rem: counterex revers ineq non bounded}).

This Theorem is the specialization to group actions of  the following more general statement for graphed groupoids:
\begin{theo}[{Theorem~\ref{thm:costccost}}] If a  sequence of \pmp graphed  groupoids $(\cg_n,\Phi_n)$ and its ultraproduct groupoid $(\cg_\ul,\Phi_\ul)$ satisfy $\lim_\ul \cost(\Phi_n) = \cost(\Phi_\ul) <\infty$ (i.e. finite cost and ``no loss of mass at infinity''), then \[\ccost_\ul((\cg_n,\Phi_n)_n)=\cost(\cg_\ul)\geq\lim_{\ul}\cost(\cg_{n}).\]
\end{theo}

This allows us to interpret the original Elek's combinatorial cost \cite{Elek-combinatorial-cost-2007} (defined as $\ccost((\cg_n,\Phi_n)_n):=\inf_L\liminf\limits_{n\to \infty}\cost_L(\cg_n,\Phi_n)$ see Remark~\ref{Rem: Elek comb cost}, Formula (\ref{def: init Elek comb cost with liminf})) in terms of cost for ultraproducts of graphed groupoids.
\begin{coro}[Combinatorial cost vs cost (Corollary~\ref{cor: combcost vs cost})]
If $(G_n)_n$ is  a sequence of finite graphs with uniformly bounded degree, then the combinatorial cost is the infimum over all non principal ultrafilters of the costs of the ultraproduct groupoids $\cg_\ul$ of the associated $(\cg_n,\Phi_n)$:
 \begin{equation*}
\ccost((G_n)_n) = \inf \{\cost(\cg_{\ul}) : \ul \textrm{ non principal ultrafilter}\}.
\end{equation*}
\end{coro}

We now concentrate on sofic groups, a class of groups 
defined in terms of ultraproducts of finite structures (see Section~\ref{sect: sofic groupoids}).

\begin{coro}[Corollary~\ref{cor: sofic gp - sofic approx - comb cost}]
If $\Gamma$ is a sofic group finitely generated by $S$ and $(\cg_n, \Phi_n)_n$ 
is a sofic approximation to $\Gamma$ (with $S$-labelled graphings), then 
  \begin{equation}
\vraicost^*(\Gamma)\geq \cost(\cg_\ul)=\ccost_\ul((\cg_n,\Phi_n)_n) \geq \lim_{\ul}\cost(\cg_{n})\geq \liminf_{n\to \infty}\cost(\cg_{n}) \geq 1.
\end{equation}
If $\Gamma$ has moreover fixed price $\vraicost_*(\Gamma)$, then
$\vraicost_*(\Gamma)=\ccost((\cg_n,\Phi_n)_n) $ and hence if $\Gamma$ has fixed price $\vraicost_*(\Gamma)=1$, we obtain $\ccost((\cg_n,\Phi_n)_n) = \lim_{n\to \infty}\cost(\cg_{n}) = 1$.
\end{coro}

\defin{Farber sequences} are sequences of finite index subgroups $(\Gamma_n)_n$ of $\Gamma$ giving particular instances of sofic approximations via the natural actions $\Gamma\acting \Gamma/\Gamma_n$ (see Section~\ref{sect: Farber sequences}). 
If $\Gamma$ is finitely generated by $S$, we denote by $\mathrm{Sch}(\Gamma/\Gamma_n, S)$ the associated Schreier graph.

\begin{coro}[Corollary~\ref{cor: nested Farber seq}]
\label{cor-intro nested Farber sequences}
If $\Gamma$ is finitely generated by $S$ and $(\Gamma_n)_n$ is a nested Farber sequence, then 
$\cost(a_{\ul})=\ccost_{\ul}(a_n)=\ccost(\mathrm{Sch}(\Gamma/\Gamma_n, S)_n)=\lim\limits_{n\to \infty} \cost(a_n)$.
\end{coro}

\begin{coro}[Corollary~\ref{cor: rk grad vs cost fixed price}]
\label{cor-intro:cor: rk grad vs cost fixed price}
If $\Gamma$ is finitely generated by $S$, has fixed price $\vraicost_*(\Gamma)$, and $(\Gamma_n)_n$ is any
(non necessarily nested) Farber sequence, then 
we have: 
\begin{equation*}
\lim\limits_{n\to \infty} \frac{\rank(\Gamma_n)-1}{[\Gamma:\Gamma_n]}+1=\ccost(\mathrm{Sch}(\Gamma/\Gamma_n, S))=\vraicost_*(\Gamma).
\end{equation*} 
\end{coro}

This corollary is a far reaching extension of a theorem of Ab\'ert and Nikolov \cite[Theorem 1]{Abert-Nikolov-12} (nested and Farber sequences) for which we provide an alternative proof for (nested and) non-necessarily Farber sequences of finite index subgroups -- see Corollary~\ref{cor: Abert-Nikolov}.
\\
The convergence of the sequence in Corollary~\ref{cor: rk grad vs cost fixed price} is not at all clear for non fixed price groups, contrarily to \cite[Theorem 1]{Abert-Nikolov-12} where it is non-increasing.
This also gives a converse to \cite[Theorem 8]{Abert-Nikolov-Gelander} where an inequality 
$\lim\limits_{n\to \infty} \frac{\rank(\Gamma_n)-1}{[\Gamma:\Gamma_n]}+1\leq \ccost(\mathrm{Sch}(\Gamma/\Gamma_n, S))$ was obtained.
The non necessarily fixed price case is considered in Theorem~\ref{cor: rk grad vs cost - non fixed price}.

These corollaries apply potentially to plenty of groups. 
Examples of groups with fixed price $1$ include infinite amenable groups, infinite-conjugacy-class (icc) inner amenable groups, direct products $\Gamma\times \Lambda$ where $\Lambda$ is infinite and $\Gamma$ contains a fixed price $1$ subgroup, groups with a normal (or even just commensurated) fixed price $1$ subgroup $N$, more generally groups with a commensurated fixed price $1$ subgroup $C$ (Corollary~\ref{cor: cost 1 commensurated}) (this applies for instance to $\SL(d, \Zmath[1/p])$, $d\geq 3$), amalgamated free products of fixed price $1$ groups over infinite amenable subgroups, Thompson’s group F, $\mathrm{SL}(n,\Zmath)$ for $n\geq 3$, non-cocompact arithmetic lattices in connected semi-simple algebraic Lie groups of $\Qmath$-rank at least $2$, groups generated by chain-commuting infinite order elements (i.e., infinite order elements whose graph of commutation is connected) sometimes called \defin{Right Angled groups} (this includes Mapping Class Groups $\mathrm{MCG}(\Sigma_g)$ of surfaces with genus $g\geq 3$, $\mathrm{Out}(\FF_n)$ for $n\geq 3$, Right Angled Artin groups with connected defining graphs), more generally groups generated by a sequence of subgroups $(\Gamma_n)$ of fixed price $1$ such that the intersections $\Gamma_{n}\cap \Gamma_{n+1}$ are infinite...

We introduce some terminology. A countable group $\Gamma$ is called a \defin{Stuck-Zimmer group} 
if any \pmp aperiodic $\Gamma$-action on a standard Borel space has a.s. finite stabilizers.
The key point about Stuck-Zimmer groups $\Gamma$ in our context, as observed in \cite{7s-2017}, is
 that under mild divergence assumptions, all sequences $(\Gamma_n)_n$ of finite index subgroups give indeed Farber sequences (see Theorem~\ref{th: Stuck-Zimmer gps and seq f i subgroups}).
If $\Gamma$ is a finitely generated Stuck-Zimmer group and  
$[\Gamma:\Gamma_n]\to \infty$, then 
for every torsion free subgroup $\Lambda<\Gamma$, the induced action $\Lambda\acting  \Gamma/\Gamma_n$ defines a sofic approximation of $\Lambda$.
Thus when  $\Gamma$ is torsion free and has fixed price $\vraicost_*(\Gamma)$, Corollary~\ref{cor: rk grad vs cost fixed price} applies: $\lim_{n}\frac{ \rank (\Gamma_n)-1}{ [\Gamma:\Gamma_n]}=\vraicost_*(\Gamma)-1$.

We introduce a relative version of Stuck-Zimmer groups for which the mild divergence assumption concerns some subgroup $N<\Gamma$. We obtain a general result (Theorem~\ref{th: normal subgr f subg --> rel Stuck-Zimmer}) which specializes to the following:
\begin{theo}[See Theorem~\ref{th: SLd rel Stuck-Zimmer}]
\label{th: intro SLd rel Stuck-Zimmer}
Let $\Lambda<\SL(d,\Zmath)$ be a  finite index subgroup and  $(\Gamma_n)_n$ be any sequence of finite index subgroups of $\Gamma=\Lambda\ltimes \Zmath^d$ such that $[\Zmath^d:\Gamma_n\cap \Zmath^d]\underset{n\to \infty}{\longrightarrow}\infty$.
Then $\Gamma\acting \Gamma/\Gamma_n$ defines a sofic approximation of $\Gamma$ and 
\begin{equation}
\lim_{n}
\frac{\left( \rank
(\Gamma_n)-1\right)}{ [\Gamma:\Gamma_n]}=0.
\end{equation}
\end{theo}
This applies more generally for $\Lambda<\SL(d,\Zmath)$ that admits no $\Lambda$-invariant finite union of infinite index subgroups of $\Zmath^d$ (see Theorem~\ref{th: SLd rel Stuck-Zimmer}). In particular, 
\begin{theo}
\label{th: intro SL2 rel Stuck-Zimmer}
Let $\Lambda$ be any subgroup of $\SL(2,\Zmath)$ containing a hyperbolic element.
Any sequence $(\Gamma_n)_n$ of finite index subgroups of $\Gamma:=\Lambda\ltimes \Zmath^2$ 
such that $[\Zmath^2:\Gamma_n\cap \Zmath^2]\underset{n\to \infty}{\longrightarrow}\infty$ defines a sofic approximation of $\Gamma$ and 
\begin{equation}
\lim_{n} \frac{\left( \rank (\Gamma_n)-1\right)}{ [\Gamma:\Gamma_n]}=0.
\end{equation}
\end{theo}
This holds for instance for the amenable group
$\Gamma=\Lambda\ltimes \Zmath^2=\left\langle \begin{pmatrix} 1 & 1 \\ 1 & 2 \end{pmatrix} \right \rangle \ltimes \Zmath^2$.

\bigskip
L\"uck approximation theorem \cite{Luc94b} considers a descending sequence $(\Gamma_n)_n$ of finite index normal subgroups   with trivial intersection of a countable group $\Gamma$ which acts freely cocompactly on a simplicial complex $L$ and it relates the Betti numbers of $\Gamma_n\backslash L$ to the $\ell^2$-Betti numbers of the action:
\begin{equation}
\lim_{n\to \infty} \frac{b_i(\Gamma_n\backslash L)}{[\Gamma:\Gamma_n]}=\beta^{(2)}_{i}(\Gamma\acting L).
\end{equation}

Bergeron-Gaboriau \cite{berggab} have interpreted this convergence as a continuity result for $\ell^2$-Betti numbers of projective systems of fields of simplicial complexes for equivalence relations.
This allowed them to extend the framework  so as to encompass non-normal subgroups and non trivial intersection descending sequences of subgroups. They established the convergence of the sequence and identified the limit as the $\ell^2$-Betti number of a laminated space \cite[Th\'eor\`eme 0.2, Th\'eor\`eme 3.1]{berggab}.

Our third main result generalizes this continuity statement in several ways: we extend it from equivalence relations to groupoids, from finite quotients to sofic and we consider arbitrary limits instead of projective ones. The use of ultrafilters forces a kind of convergence of the objects as fields of simplicial complexes with groupoid actions and we identify the limits.
We state our result under some boundedness assumption (see Section~\ref{sect: L2 Betti}, see Theorem~\ref{thm:limitgeneral} for the version without the boundedness assumption): 

\begin{theo}[Theorem~\ref{thm:limitbounded}]
  For every $n\in\Nmath$ let $(\cg_n,\Phi_n)$ be a {\bf sofic} graphed \pmp groupoid over $X_n$, let $\Sigma_n$ be a $(\Phi_n,L)$-uniformly locally bounded $\cg_n$-simplicial complex for some fixed $L\in\mathbb N$. Then
\[\lim_\ul\beta^{(2)}_i(\Sigma_n, \cg_n)=\beta^{(2)}_i(\Sigma_\ul,\cg_\ul).\]

In particular if the sequence of $\cg_n$-simplicial complexes $\Sigma_n$ is asymptotically uniformly $k$-connected (see Proposition~\ref{prop:contract}), then for every $i\leq k$
\[\lim_\ul\beta^{(2)}_i(\Sigma_n, \cg_n)=\beta^{(2)}_i(\cg_\ul).\]
\end{theo}

Our approach with ultraproducts allows us  to relate the Laplace operators for the finite objects to that of our limit object (compare \cite[Proposition~6.1]{2005-Elek-Szabo-hyperlinearity} or \cite{thom08} where the finite induced operators do not clearly appear as Laplace operators). Observe however that we make a crucial use of a Lemma of Thom  \cite{thom08} (see Lemma~\ref{lem:lueck_approximation}).
The above theorem is also related to a (recent and independent) work of Schrödl \cite{Schrodel-L2-Betti-simpl-cplx-2018} on $\ell^2$-Betti numbers and Benjamini-Schramm convergence (see also \cite{Elek-10-Betti-test}). Indeed, 
similarly to the case of local-global convergence of graphs (Remark~\ref{rem: local-global cv}), Benjamini-Schramm limits of random simplicial complexes can be expressed in terms of the ultraproduct of the sequence of the associated field of simplicial complexes (see Remark~\ref{rem: BS-CV of complexes} and Remark~\ref{rem: random unimodular simplicial complexes}).

As a corollary of the above theorem, we obtain the following geometric generalization of L\"uck approximation Theorem for (non-nested, non normal) sequences of subgroups: 
\begin{coro}[Corollary~\ref{cor:approxgroups}]
  Let $\Gamma$ be a countable group acting freely cocompactly on a simplicial complex $L$. Let $(\Gamma_n)_n$ be a (non necessarily nested) Farber sequence of finite index subgroups. Then for every $i$ we have:
\[\lim\limits_{n\to \infty}\frac{b_i(\Gamma_n\backslash L)}{[\Gamma:\Gamma_n]}=\beta^{(2)}_i(\Gamma\acting L).\] 
In particular, if $L$ is $k$-connected, then for every $i\leq k$ we have:
\[\lim\limits_{n\to \infty}\frac{b_i(\Gamma_n)}{[\Gamma:\Gamma_n]}=\beta^{(2)}_i(\Gamma).\] 
\end{coro}

 More generally, we obtain the following generalization of \cite[Theorem 3.1]{berggab} for non nested, non Farber sequences:
\begin{coro}[Corollary~\ref{Cor: non Farber sq}]
  Let $(\Gamma_n)_n$ be any sequence of finite index subgroups of a finitely generated group $\Gamma$.
  Assume $\Gamma$ acts freely cocompactly on the simplicial complex $\Sigma$.
Let $a_{\ul}:\Gamma\acting X_{\ul}$ be the ultraproduct of the actions $a_n:\Gamma\acting \Gamma/\Gamma_n$ and $\cn_\ul:=\{(x,\gamma)\in X_{\ul}\times \Gamma: \gamma x=x\}  $ its maximal totally isotropic subgroupoid.
  Then for every $i\leq k$ we have
\[\lim\limits_{n\in \ul}\frac{b_i(\Gamma_n\backslash \Sigma)}{[\Gamma:\Gamma_n]}=\beta^{(2)}_i(\cn_\ul \backslash \left(X_{\ul}\times \Sigma\right),\RR_\ul),\] 
where $\RR_\ul=\cg_\ul/\cn_\ul$ is the \pmp equivalence relation of the ultraproduct action $a_\ul$ and $\cn_\ul \backslash \left(X_{\ul}\times \Sigma\right)$ is the $\RR_{\ul}$-simplicial complex defined in Example~\ref{ex: quotient fibred space}.
\end{coro}

\bigskip
We indicate now how these approximation results must be modified when the boundedness assumptions are not satisfied.

 Let $\Gamma$ be a finitely generated residually finite group and let $(\Gamma_n)_{n\in \Nmath}$ be a nested sequence of finite index normal subgroups with trivial intersection. 
L\"uck approximation theorem $\frac{b_1(\Gamma_n)}{[\Gamma:\Gamma_n]} \to \beta_1^{(2)}(\Gamma)$ does not hold in general when $\Gamma$ is not finitely presented, as examplified by torsion groups (thus $b_1(\Gamma_n)=0$) with positive first $L^2$-Betti number \cite{Lueck-Osin-2011}. 
This defect can be fixed by taking a double limit as exposed below and this idea will then be pushed to the framework of ultraproduct (Theorem~\ref{thm:limitgeneral}).
\begin{prop}
 Let $\Gamma=\langle \gamma_1, \gamma_2, \cdots, \gamma_d\vert r_1, r_2, \cdots, r_k, \cdots\rangle$ be a presentation of $\Gamma$ with finitely many generators. Let $L$ be the associated ($2$-dimensional) Cayley $\Gamma$-complex and let $L_j$ be the $\Gamma$-invariant co-compact sub-complex associated with the first $j$ relations $ r_1, r_2, \cdots, r_j$.
Then 
 \[\beta_1^{(2)}(\Gamma)=
 \lim\limits_{j\to \infty}
 \lim\limits_{n\to \infty}  
 \frac{b_1(\Gamma_n\backslash L_j)}{[\Gamma:\Gamma_n]}.
 \]
\end{prop}
The substance of the counter-examples of \cite{Lueck-Osin-2011} relies on the fact that these limits cannot be exchanged. Indeed on the one hand,  $ \lim\limits_{n\to \infty}  
 \frac{b_1(\Gamma_n\backslash L_j)}{[\Gamma:\Gamma_n]}=\beta_1^{(2)}(\Gamma\acting L_j)$ by L\"uck approximation theorem \cite{Luc94b}, and $\beta_1^{(2)}(\Gamma\acting L_j)$ decreases  with $j$ to $\beta_1^{(2)}(\Gamma)$. On the other hand, for every $n$, we have $b_1(\Gamma_n\backslash L_j)=0$ for large enough $j$ since the canonical surjective maps $\pi_1(\Gamma_n\backslash L_j)\to \pi_1(\Gamma_n\backslash L_k)$ (for $j\leq k$) produce as direct limit 
 the torsion group $\Gamma_n$, so that $\pi_1(\Gamma_n\backslash L_j)$ is generated by torsion elements for large enough $j$ and its abelianization is thus finite.

More generally, in higher dimension, we obtain the following.
\begin{prop}
\label{prop: approx th. non cocompact}
 Let $\Gamma$ be a residually finite group and let $(\Gamma_n)_{n\in \Nmath}$ be a nested sequence of finite index normal subgroups with trivial intersection. 
If $\Gamma$ acts freely on the $d$-connected simplicial complex $L$, we write $L$ as an increasing union $L=\cup\nearrow L_j$ of $\Gamma$-invariant co-compact sub-complexes $L_j$.
 The inclusions $ L_j\subset  L_k$ for $j\leq k$ induce linear maps 
 $ H_d(\Gamma_n\backslash L_j, \Cmath)\to H_d( \Gamma_n\backslash L_k, \Cmath)$ between the $d$-dimensional homology of the compact quotients. 
  Then 
 \[\beta_d^{(2)}(\Gamma)=\lim\limits_{j\to \infty}\lim\limits_{j\leq k,\,  k\to \infty}\lim\limits_{n\to \infty}  
 \dim_{\Cmath}\Ima\Bigl(H_d(\Gamma_n\backslash L_j, \Cmath)\to H_d( \Gamma_n\backslash L_k, \Cmath)\Bigr).\]
 \end{prop}
This proposition follows from Theorem~\ref{thm:limitgeneral} and the fact that, for free \pmp actions, the $L^2$-Betti numbers of the equivalence relation generated by the action agree with the ones of the group itself. The above result holds more generally for Farber sequence of subgroups.

In  \cite{Elek-combinatorial-cost-2007}, Elek introduced another quantity for a graph sequence $(G_n)_n$: the 
 $\beta$-invariant over the field $\mathbb{K}$:
 $$\beta_\mathbb{K}( (G_n)_n):=\inf_q \liminf_{n\to \infty}\frac{\vert E(G_n)\vert - \dim_\mathbb{K} V_q(G_n)}{\vert V(G_n)\vert} -1,$$
where $\vert E(G_n)\vert$ and $\vert V(G_n)\vert$ are the number of edges and vertices and where $ V_q(G_n)$ is the $\mathbb{K}$-vector space spanned by the cycles of length $\leq q$ in the space of $1$-chains.
Equivalently (when $\vert V(G_n)\vert \to \infty$), this is also
 $$\beta_\mathbb{K}( (G_n)_n):=\inf_q \liminf_{n\to \infty}\frac{b_1(G_n^q)}{\vert V(G_n)\vert},$$
where $G_n^q$ is obtained from the graph $G_n$ by gluing $2$-cells along every cycle of length $\leq q$.

Inspired by that, we introduce higher dimensional invariants for uniformly bounded vertex degrees sequences of graphs $(G_n)_n$.
We consider $R^q(G_n)$, the $q$-Rips complex of $G_n$ for all integers $q$
and we introduce the \defin{$d$-dimensional $\beta$-invariant}:
$$\beta_{d,\mathbb{K}}( (G_n)_n):=\lim_{q\to \infty}\lim_{p\to \infty} \lim_\ul \frac{\dim_\mathbb{K} \mathrm{Im} \left(H_d(R^{q}(G_n), \mathbb{K})\to H_d(R^{q+p}(G_n),\mathbb{K})  \right)}{\vert V_n\vert}.$$

If $\mathbb{K}=\Qmath$ and $(G_n)_n$ has uniformly bounded degree, 
Theorem~\ref{thm:limitgeneral} implies that
$\beta_{d,\mathbb{\Qmath}}( (G_n)_n)=\beta_d^{(2)}(\RR_\ul)$
where 
$\RR_\ul$ is the ultraproduct \pmp equivalence relation of the sequence of graphs $(G_n)_n$ (with well defined coarse structure, see Remark~\ref{ex:ultraproduct_of_graphs}).
In particular, if $(G_n)_n$ corresponds to a sofic approximation of a countable group $\Gamma$, one gets 
$\beta_{d,\mathbb{\Qmath}}( (G_n)_n)=\beta_d^{(2)}(\Gamma)$.

\bigskip

We briefly mention an extension we do not consider in this version in order to limit the length of the present article.
There exists also a notion of coarse equivalence between pairs of sequences $(\cg_n,\Phi_n)_n$ and $(\ch_n, \Psi_n)_n$ of graphed groupoids (compare \cite{Alekseev-Finn-Sell-sofic-boundaries-2016} and \cite{Kaiser-2017} which considers finite structures). 
The groupoids $\cg_n$ and $\ch_n$ are in particular stably orbit equivalent, so that they are related by a compression constant $\iota(\cg_n,\ch_n)$.
When the sequence of compression constants is bounded, the ultraproduct groupoids $\cg_\ul$ and $\ch_\ul$ are orbit equivalent (with compression constant $\iota(\cg_\ul,\ch_\ul)=\lim\limits_{\ul} \iota(\cg_n,\ch_n)$) and their costs and $\ell^2$-Betti numbers are related as the standard ones: 
\[\cost(\ch_\ul)-1=\iota(\cg_\ul,\ch_\ul)(\cost(\cg_\ul)-1)\ \text{ and }  \beta_j^{(2)}(\ch_\ul)=\iota(\cg_\ul,\ch_\ul) \beta_j^{(2)}(\cg_\ul).\]

\bigskip

Warning: We have obtained some of our main results  (especially Theorem~\ref{thm:costccost})
 a couple of years ago. 
  We discovered during the “Measured Group Theory” ESI workshop in Vienna (February 2016) that some of our results had a non trivial intersection with some ongoing work of Abért-Tóth (especially the first equality in our Theorem~\ref{th-intro:equality of lim cost for fixed price} in the case of local-global convergent sequences and our Corollary~\ref{cor-intro:cor: rk grad vs cost fixed price}). 
Some accident prevented us from releasing our papers simultaneously, but the overlaps with \cite{Abert-Toth-unif-rk-grad-cost-loc-glob-cv} must be considered parallel and independent results.

\bigskip
\bigskip


\paragraph{Acknowledgments:}
We are grateful to Luiz Cordeiro and László Márton Tóth for their comments and their careful reading of a
preliminary version (November 2018). We also thank Gabor Elek pointing out some references.
The three authors are supported by the ANR project GAMME (ANR-14-CE25-0004).
The first named author is supported by the ERC Consolidator Grant No. 681207 and by the Deutsche Forschungsgemeinschaft (DFG, German Research Foundation) – 281869850 (RTG 2229). 
The second and third named authors are supported by the CNRS and by the LABEX MILYON (ANR-10-LABX-0070) of Université de Lyon, within the program ``Investissements d'Avenir" (ANR-11-IDEX-0007) operated by the French National Research Agency (ANR). The third author is also supported by the ANR project AGIRA (ANR-16-CE40-0022). We are grateful to the Erwin Schr{\"o}dinger Institute for Mathematics and Physics for hospitality during “Measured group theory” program in 2016.

\newpage

\section{Preliminaries}

In this text we will work with probability spaces, often denoted by $(X,\rb,\mu)$ or simply by $(X,\mu)$ whenever the $\sigma$-algebra can be safely omitted from the notation. We will make no further assumption on the probability spaces and in particular we will not assume the $\sigma$-algebra to be countably generated. While working in this context one has to be more careful with some standard definitions since the classical Lusin-Novikov's uniformization theorem and the von Neumann-Jankow selection theorem do not hold. Therefore we will carefully define all the objects we will need.

A \defin{partial isomorphism} from the probability space $(X,\mu)$ to $(Y,\nu)$ is a bi-measurable bijection $\ph: \dom(\ph) \subset X\to \tar(\ph) \subset Y$ between measurable subsets of $X$ and $Y$ which preserves the measure.

\begin{defi}\label{def:fibred_space}
  A countable \defin{fibred space} over $X$ consists in a measure space $(\cf,\ra,\widetilde \mu)$ and a measurable map $\pi\colon\cf\rightarrow X$ such that there is a measurable subset $\widetilde \cf \subset X \times \Nmath$ (for the product $\sigma$-algebra and product measure) and a bi-measurable measure preserving bijection $u \colon \cf \to \widetilde \cf$ such that $u(f) \in \{\pi(f)\} \times \Nmath$ for every $f$.
\end{defi}

\begin{rema}\label{rem: count fibred sp. and measurab}
By definition, countable fibred spaces are the measurable subspaces of the trivial fibred space $X \times \Nmath$. 
Partitioning each $A \in \ra$  into the measurable subsets $A_i:=A\cap X\times\{i\}$ shows that:
 \begin{itemize}
 \item the map $x \in X\mapsto |\pi^{-1}(x) \cap A|$ is measurable and 
\begin{equation}\label{eq: tilde mu on fibred}
\tilde\mu(A)=\int_X\left| \pi^{-1}(x)\cap A\right| d\mu(x);
\end{equation}
\item the subset $\pi(A)=\cup_i \pi(A_i)$ is measurable;
\item if $|\pi^{-1}(x)\cap A|\leq 1$ for almost every $x\in X$, then $\pi\bigr|_A:A\to \pi(A)$ is a measure preserving isomorphism;
\item  any measurable section $\ph:B\subset X\to \cf$ of $\pi$ is a measure preserving isomorphism between $B$ and $\ph(B)$
(consider for $B_i:=\ph^{-1} (X\times \{i\})$).
 \end{itemize}
\end{rema}

When $X$ and $\cf$ are standard Borel spaces and $\pi \colon \cf \to X$ is a measurable map with countable fibers, the standard selection theorem of Lusin-Novikov \cite{Lusin-Novikoff-1935}
produce measurable sections allowing to realize $\cf$ as a Borel subset of $X \times \Nmath$. 
The \defin{height} $H(\cf)$ of a fibred space $\cf$ is the essential supremum of the function $x\mapsto |\pi^{-1}(x)|$ for $x\in X$. We say that $\cf$ has \defin{finite height} if $H(\cf)$ is finite.
The \defin{fibred product} $(\cf_1*\cf_2,\ra,\pi)$ of two countable fibred spaces
$(\cf_1,\ra_1,\pi_1)$ and $(\cf_2,\ra_2,\pi_2)$ over $X$ is defined as usual by 
\begin{equation*}
(\cf_1,\pi_1)*(\cf_2, \pi_2):=\left\{ (x,y)\in \cf_1\times \cf_2:\ \pi_1(x)=\pi_2(y)\right\},
\end{equation*}
(in short $\cf_1*\cf_2$) equipped with the projection $\pi(x,y)= \pi_1(x) = \pi_2(y)$ and the $\sigma$-algebra coming from the injection into $X \times\Nmath \times \Nmath$ if $\cf_1$ and $\cf_2$ are realized inside $X \times \Nmath$. This is easily seen to be independent of the realizations.

\subsection{Measured groupoids}

The main object of study of this text will be \textit{measured groupoids}. 
\begin{defi}
A \defin{groupoid} $\cg$ over the set $X$ is the data of
 \begin{itemize}
 \item a set $\cg$ containing  $X$; the elements of $X$ are called the \textit{units} of $\cg$,
 \item a \textit{source map} $s:\cg\rightarrow X$ and a \textit{target map} $t:\cg\rightarrow X$ which are the identity when restricted to $X$,
 \item for every $g$ and $h\in\cg$ such that $s(g)=t(h)$ there is a \textit{product} $gh\in\cg$ satisfying $s(gh) = s(h)$ and $t(gh)  =t(g)$.
 \item for every $g\in\cg$ there exists a unique \textit{inverse} $g^{-1}\in\cg$ such that $s(g^{-1}) = t(g)$, $t(g^{-1})=s(g)$, $g g^{-1} = t(g)$ and $g^{-1}g=s(g)$.
 \end{itemize}
for which we require that the following identities hold
\begin{itemize}
\item the product is associative, that is $g(hk)=(gh)k$ whenever this is defined,
\item for every $g\in \cg$ we have that $t(g)g=g s(g) = g$.
\end{itemize}
\end{defi}
In this text we will always assume that every groupoid is $s$-\textit{countable}, i.e.,
\begin{itemize}
\item for every $x\in X$ the set $s^{-1}(x)$ is countable. 
\end{itemize}
Observe that the multiplication induces a map $(\cg,s)*(\cg,t)\rightarrow \cg$.

\begin{defi}\label{def: pmp groupoid}
  A \defin{\pmp groupoid} over $(X,\mu)$ is a groupoid $\cg$ over $X$ equipped with a $\sigma$-algebra $\ra$ and a measure $\widetilde \mu$ such that
  \begin{itemize}
  \item $(\cg,\ra,\widetilde \mu,s)$ and $(\cg,\ra,\widetilde \mu,t)$ are both countable fibred spaces over $(X,\mu)$,
  \item the maps $g \in \cg \mapsto g^{-1} \in \cg$ and $(g,h) \in (\cg,s) \ast (\cg,t) \mapsto gh$ are measurable.
  \end{itemize}
\end{defi}
In particular, the measure $\widetilde{\mu}$ coincides with both fibred measures for $s$ and $t$ (see (\ref{eq: tilde mu on fibred})).
Observe that we restrict ourselves to measured groupoids with countable fibers, a.k.a. discrete, or countable, or $r$-discrete (see \cite{2000-Anantharaman-Renault-Amen-gpoids}); but we do not assume that the base space $X$ is standard.

Let $\cg$ be a groupoid over $X$. The \defin{isotropy group} at $x\in X$ is the group 
$\cg_x:=\{g\in\cg: s(g)=t(g)=x\}$. When the isotropy groups are all trivial, the groupoid is called a \defin{principal groupoid} or an \defin{equivalence relation}. For \pmp groupoids, the triviality is required only almost surely.

\begin{exam}[Group action groupoid]
\label{ex: groupoid from pmp action}
  Let $\Gamma$ be a countable group acting on a probability space $(X,\mu)$ preserving the measure. Then $\cg_{\Gamma\acting X}:=\Gamma \times X$ (with the product measured structure) is a \pmp groupoid over the probability space $X$ (identified with the subset $\{1\}\times X$ of $\Gamma \times X$) where $s(\gamma,x)=x$, $t(\gamma,x)=\gamma x$ and $(\gamma,x)(\gamma',x')=(\gamma\gamma',x')$ if $\gamma'x'=x$.
\end{exam}

\begin{exam}
  Lusin-Novikov's uniformization theorem \cite{Lusin-Novikoff-1935} (see also \cite[Theorem 18.10]{Kechris-class-descr-th-1995}), implies that if $(X,\mu)$ is a standard Borel space and $\RR\subset X \times X$ is an Borel \pmp equivalence relation with countable classes (see \cite{FeldmanMoore-1977}), then $\RR$ is a \pmp principal groupoid as in Definition~\ref{def: pmp groupoid} where $s$ and $t$ are respectively the second and first coordinate projections $X \times X \to X$.
\end{exam}

\begin{exam}[Equivalence relation associated to \pmp actions]
\label{ex: example equiv-rel pmp act non meas.} 
We warn the reader that, given a \pmp groupoid $\cg$ over $X$, the underlying equivalence relation $\RR_\cg := \{(s(x),t(x)): x \in \cg\}$ is in general \emph{not a \pmp groupoid} in the sense of Definition~\ref{def: pmp groupoid}, even for the equivalence relation $\RR_{\Gamma \acting X} = \{(x,\gamma x): x\in X,\gamma \in \Gamma \}$
coming from a \pmp action of a countable group $\Gamma$.

For example, consider the space $(X=[-1,1],\rb,\mu)$ when $\rb$ is the $\sigma$-algebra of all Borel-measurable subsets $B$ of $[-1,1]$ satisfying $B=-B$, and $\mu$ the normalized Lebesgue measure. Let $\Omega \subset(0,1]$ be a subset which is not Lebesgue-measurable. The map $\ph:x \mapsto -x$ on $\Omega \cup -\Omega$ and $x \mapsto x$ on its complement is a measurable \pmp involution of $([-1,1],\rb,\mu)$. It defines a \pmp action $\alpha$ of $\Gamma= \Zmath/2\Zmath$ on $[-1,1]$. Then $\cg_{\alpha}$ is a \pmp groupoid (Example~\ref{ex: groupoid from pmp action}), but
$\RR_{\alpha} = [-1,1] \cup \{(x,-x): x\in \Omega \cup -\Omega\}$ is not a \pmp groupoid, since $x\mapsto \vert \RR(x)\vert$ is not measurable (see Remark~\ref{rem: count fibred sp. and measurab} (1)). 
Indeed the field of isotropy groups $x\mapsto \cg_x$ is not measurable in any reasonable sense. 

Anticipating Definition~\ref{def: tot. vs meas isotrop, free}, the map $\ph$ is measurably isotropic but not totally isotropic. If instead, $\Omega=X$, then the corresponding action $\beta:\Gamma \acting X$ is totally free but not measurably free.

However, if the fixed-point set $\Fix (\gamma) := \{x \in X:\gamma x=x\}$ is measurable for every $\gamma \in \Gamma$, then $\RR_{\Gamma \acting X}$ is a \pmp principal groupoid (define the $\sigma$-algebra on $\RR_{\Gamma \acting X}$ as the set of all countable unions of sets of the form $\{(x,\gamma x): x \in A\}$ for $A \subset X$ measurable and $\gamma \in \Gamma$; the fixed-point set condition takes care of the intersections). We warn the reader that the condition that $\Fix(\gamma)$ is measurable is not necessary for $\RR_{\Gamma \acting X}$ to carry a \pmp principal groupoid structure. An example is given by $\Gamma =\Zmath/2\Zmath \times \Zmath/2\Zmath \acting X=[-1,1]$ is defined by $(\varepsilon_1,\varepsilon_2) \cdot x=(-1)^{\varepsilon_2} \alpha(\varepsilon_1) x$ for the action $\alpha$ of $\Zmath/2\Zmath$ considered above. Indeed, $\Fix(0,1)=\Omega \cup -\Omega$ is not measurable, but $\RR_{\Gamma \acting X}$ is just $\{(x,\pm x):x \in X\}$ and clearly carries a \pmp groupoid structure.

More generally given an at most countable collection of partial isomorphisms $\Phi=(\ph_i)_i$ of $(X,\mu)$, we can define an equivalence relation $\RR$ by declaring that two points $(x,y)$ in $X$ are equivalent if there exists $(i_1,\ldots,i_l) \in \Nmath^l$ and $(\varepsilon_1,\ldots,\varepsilon_l) \in \{-1,1\}^l$ such that $\ph_{i_1}^{\varepsilon_{1}}\ldots\ph_{i_l}^{\varepsilon_{l}}(x)=y$. If moreover for every $\Phi$-word $\ph = \ph_{i_1}^{\varepsilon_{1}}\ldots\ph_{i_l}^{\varepsilon_{l}}$, the fixed-point set $\{x : \ph(x)=x\}$ is measurable, that is if the associated groupoid is \textit{realizable}, then $\RR$ is a \pmp groupoid for the $\sigma$-algebra generated by $\{(x,\ph(x)): x\in A\}$ where $\ph$ is a $\Phi$-word and $A\subset \dom(\ph)$ measurable (see Section~\ref{subsect: presentation of groupoids} and \ref{subsect: realizability}).  
\end{exam}

\subsection{Full pseudogroup}

From now on, let us fix a \pmp groupoid $\cg$ over $(X,\mu)$. 

\begin{defi}\label{def: meas. bisection}
A \defin{measurable bisection} of a \pmp groupoid $\cg$ is a measurable subset $\ph \subset \cg$ such that the restriction of $s$ and $t$ to $\ph$ are injective. 
These restrictions are partial isomorphisms (Remark~\ref{rem: count fibred sp. and measurab}). Thus every bisection $\ph$ delivers an associated partial isomorphism $\hat \ph:= t \circ (s\bigr|_\ph)^{-1}$ of $(X,\mu)$.
\end{defi}

\begin{lemm}
  If $\ph$ and $\psi$ are bisections of the \pmp groupoid $\cg$, then the product $\ph\psi:=\{gh:g\in \ph,h \in \psi, s(g)=t(h)\}$ is a bisection.
\end{lemm}
\begin{proof}
  We have to prove that $\ph \psi$ is a measurable subset of $\cg$. By Remark~\ref{rem: count fibred sp. and measurab} it is enough to show that $(s\bigr|_{\ph\psi})^{-1}$ is measurable. First observe that a point $x \in X$ belongs to $s(\ph \psi)$ if and only if there is $g \in \ph$, $h \in \psi$ such that $s(h)=x$ and $s(g) = t(h)=\hat\psi(x)$. So $s(\ph\psi) = s(\psi) \cap \hat{\psi}^{-1}(s(\ph))$ is measurable. Moreover $(s\bigr|_{\ph\psi})^{-1}(x) = (s\bigr|_{\ph})^{-1}(\hat{\psi}(x))(s\bigr|_{\psi})^{-1}(x)$ is measurable, as the composition of the map $x\mapsto ((s\bigr|_{\ph})^{-1}(\hat{\psi}(x)),(s\bigr|_{\psi})^{-1}(x)) \in (\cg,s)\ast(\cg,t)$, which is measurable by definition of fibred products, and of the product map $(\cg,s)\ast(\cg,t) \to \cg$ which is measurable by assumption. 
\end{proof}

We observe that there is a countable family $\{\ph_i\}_i$ of measurable bisections which covers $\cg$ in the sense $\cup_i \ph_i=\cg$. Indeed, the assumption that $(\cg,s)$ (respectively $(\cg,t)$) is a fibred space over $X$ provides countable families $\{\ph_i^s\}_i$ and$\{\ph_j^t\}_j$ of measurable subsets which both cover $\cg$ and on which $s$ (resp. $t$) is injective. The family $\{\ph_i^s\cap \ph_j^t\}_{i,j}$ does the job. The next remark thus follows:
\begin{rema}
\label{rem: X is measurable in G}
The space of units of $X\subset \cg$ is a measurable subset of $\cg$ (since $X = \cup_i \ph_i^{-1} \ph_i$), and the measures $\mu$ and $\widetilde \mu$ coincide on $X$.  
For every bisection $\ph$, its unit set $\ph\cap X$ is thus measurable (contrarily to the isotropy set $\Fix(\ph)=\{g\in \ph: s(g)=t(g)\}$, see Example~\ref{ex: example equiv-rel pmp act non meas.}).
Similarly the product $(g,h)\in (\cg,s)\ast(\cg,t)\mapsto gh\in\cg$ maps measurable subsets to measurable subsets.
\end{rema}

The \defin{full pseudogroup} of $\cg$, denoted $[[\cg]]$, is the set of measurable bisections of $\cg$ modulo null-sets. The full pseudogroup  admits a natural metric, the $L_1$ distance, defined by $d(\ph,\psi) := \widetilde \mu(\ph \Delta \psi)$ and a \defin{trace}, the measure of the fixed points, which can be defined by the formula $\tau(\ph):=\widetilde \mu(X\cap \ph)$. It is easy to observe that this trace is the restriction of the trace of the associated $II_1$ von Neumann algebra, see Section~\ref{sec:vn}.

It can be endowed with the following operations:
\begin{itemize}
\item for $\ph,\psi \in [[\cg]]$, their \textit{product} is $\ph\psi=\{gh:g\in \psi,h \in \psi, s(g)=t(h)\}$;
\item for $\ph\in [[\cg]]$, its inverse is $\ph^{-1}= \{g^{-1}: g \in \ph\}$;
\item for $\ph,\psi \in [[\cg]]$ their \textit{intersection} is $\ph\cap \psi$;
\item if $\ph,\psi\in [[\cg]]$ are such that $s(\ph) \cap s(\psi) = \emptyset$ and $t(\ph) \cap t(\psi)=\emptyset$, then the \textit{union} or \textit{join} $\ph\vee \psi$ is the bisection $\ph \cup \psi$.
\end{itemize}
Moreover the following relations are satisfied
\begin{itemize}
\item the product is associative;
\item all the idempotents ($\ph \ph=\ph$) commute;
\item for every $\ph\in [[\cg]]$ we have $\ph \ph^{-1}\ph=\ph$.
\end{itemize}

\begin{defi}
\label{def: tot. vs meas isotrop, free}
A measurable bisection $\ph$ of a \pmp groupoid $\cg$ is
\begin{enumerate}
\item\defin{measurably isotropic} if for every measurable subset $U\subset s(\ph)$ we have $\mu (\hat\ph(U)\Delta U)=0$;
\item\defin{totally isotropic}: if $s(g)=t(g)$ for $\tilde\mu$-almost every $g\in {\ph}$;
\item \defin{measurably free} if for every measurable subset $A\subset s(\ph)$ of positive measure there exists a measurable subset $B\subset A$ such that $\mu(\hat\ph(B)\Delta B)\not= 0$;
  \item\defin{totally free} if $s(g)\neq t(g)$ for $\tilde\mu$-almost every $g\in {\ph}$.
\end{enumerate}
\end{defi}
In a standard Borel probability space, there is no difference between totally and measurably isotropic (resp. free). 
Example~\ref{ex: example equiv-rel pmp act non meas.} illustrates these notions.
 
\begin{prop}\label{prop:manyinvo-gen}
Every bisection $\ph$ of a \pmp groupoid $\cg$ admits a measurable decomposition 
  $\ph = \Fix_{\rb}(\ph)\sqcup M_F(\ph)$ (unique up to a null set) such that 
\begin{itemize}
  \item $\Fix_{\rb}(\ph)$ is measurably isotropic;
  \item $M_F(\ph)$ is measurably free.
\end{itemize}

  Moreover, the measurably free part $M_F(\ph)$ admits a countable measurable partition $M_F(\ph)=\sqcup_i A_i$ (up to a null set) such that $\mu(s(A_i)\cap t(A_i))=0$ and $\widetilde \mu(\cup_{i=1}^n A_i) \geq (1-(2/3)^n) \widetilde \mu(M_F(\ph))$. 
\end{prop}
Observe that the isotropy set (or fixed-point set) $\Fix(\ph):=\{g\in \ph: s(g)=t(g)\}$ is a possibly non-measurable subset of $ \Fix_{\rb}(\ph)$ (see Example~\ref{ex: example equiv-rel pmp act non meas.}).

\begin{proof}
Define $a \leq 1$ as the supremum of $\widetilde \mu(\cup_i B_i)$ over all countable union of measurable subsets of $\ph$ satisfying $\widetilde \mu(s(B_i) \cap t(B_i))=0$. Observe that $a$ is a maximum: if for every $n$ we have a collection $\{C_j^n\}_j$ satisfying $\widetilde \mu(\cup_j C_j^n) \geq a-1/n$, then the collection $\{B_i\}_i = \{C_j^n\}_{j,n}$ satisfies $\widetilde \mu(\cup_i B_i)=a$. Pick such a family $\{B_i\}_i$ attaining the maximum and set $F:=\ph \setminus \cup_i B_i$. 

$F$ is measurably {isotropic} for otherwise $F$ would contain a measurable subset $V$ such that $\widetilde \mu(s(V)\Delta t(V))\not=0$, and the subset $W=V\cap s^{-1}(s(V) \setminus t(V))$ is non null and satisfies $\mu(s(W)\cap t(W))=0$, thus contradicting the maximality of $\widetilde \mu(\cup_i B_i)$. Set $\Fix_{\rb}(\ph):=F$.

By taking restrictions of each $B_i$ one can assume that the family $\{B_i\}_i$ is a partition of  $M_F(\ph):=\ph\setminus \Fix_{\rb}(\ph)$. 

The next step is to modify the family so that $\widetilde \mu(B_1) \geq \frac 1 3 \widetilde \mu(M_F(\ph))$. This is done, inductively on $i$. Set $B_1^1:=B_1$ and inductively set $B_i'=\{g \in B_i : s(g) \in t(B_{1}^{i-1}) \textrm{ or }t(g) \in s(B_1^{i-1})\}$ and $B_1^i:=B_1^{i-1}\cup(B_i\setminus B_i')$ for $i>1$. Set $B_1'=\cup_i B_1^i$ and observe that $\{B_i'\}_i$ is again a partition of $M_F(\ph)$ satisfying $\mu(s(B_i') \cap t(B_i'))=0$. Moreover, $M_F(\ph)\setminus B_1'$ is contained in the union $(s\bigr|_\ph)^{-1}\circ t(B'_i) \cup (t\bigr|_\ph)^{-1}\circ s(B'_i)$ so that $\widetilde \mu(B_1')\geq \frac 1 3 \widetilde \mu(M_F(\ph))$.  

By repeating the preceding argument inductively for all $B_i$, we obtain in the end that $\widetilde \mu(B_i) \geq \frac 1 3 \widetilde \mu(M_F(\ph) \setminus \cup_{j<i} B_j)$, which finishes the proof.

The uniqueness up to a null set is clear, as a bisection which is
simultaneously measurably isotropic and measurably free has measure
$0$.
\end{proof}

\subsection{Regular representations and von Neumann algebras of a groupoid}
\label{sec:vn}

Let $\cg$ be a \pmp groupoid over $(X,\mu)$.

If $\ph\in [[\cg]]$, define a partial isometry $\lambda(\ph)$ (the left-regular representation of $[[\cg]]$) on $L^2(\cg,\widetilde \mu)$ as follows. If $f \in L^2(\cg,\widetilde \mu)$ one defines $\lambda(\ph)f(g) := 0$ if $t(g) \notin t(\ph)$ and $\lambda(\ph)f(g) := f(h^{-1}g)$ otherwise, where $h$ is the unique element of $\ph$ such that $t(h)=t(g)$. For example, for $\psi\in [[\cg]]$, $\lambda(\ph)\chi_\psi = \chi_{\ph\psi}$. The left von Neumann algebra of $\cg$ is $\LL(\cg)=\lambda([[\cg]])''$. It carries a normal faithful trace $\tau = \langle \cdot \chi_X,\chi_X\rangle$. The von Neumann subalgebra generated by the idempotents of $[[\cg]]$ (the indicator functions of measurable subsets of $X$) is isomorphic to $L^\infty(X,\mu)$. Its commutant is the right von Neumann algebra $\LL'(\cg)$, generated by the right-regular representation $\rho$ which is characterized by $\rho(\ph)\chi_\psi = \chi_{\psi \ph^{-1}}$. 

\subsection{Coarse structures and graphings}
\label{subsect : coarse structures and graphings}
We now define a notion of coarse structure on a \pmp groupoid, inspired by \cite{Roe-2003}.
\begin{defi}
  A \defin{coarse structure} on a \pmp groupoid $\cg$ is a collection $\coarse$ of measurable subsets of $\cg$, called the \defin{controlled sets}, which contains the units and is closed under the formation of measurable subsets, inverses, products and finite unions. A coarse structure is \defin{generating} if $\cg$ is the union of all controlled sets.
A coarse structure on $\cg$ has \defin{bounded geometry} if for every controlled set $E$, the measurable function $x \in X \mapsto |E \cap s^{-1}(x)|$ is essentially bounded.

A \defin{coarse structure} on a sequence $(\cg_n)_n$ of \pmp groupoids over $X_n$ is a collection of sequences $(E_n)_n$ of measurable subsets $E_n \subset \cg_n$, called the \defin{controlled sequences}, which is closed under the same operations. This means that if the sequences $(E_n)_n$ and $(F_n)_n$ are controlled, then $(E_n^{-1})_n$, $(E_n F_n)_n$ and $(E_n \cup F_n)_n$ are also controlled.
A coarse structure $\coarse$ on $(\cg_n)_n$ has \defin{bounded geometry} if for every $(E_n) \in \coarse$, the measurable function $n \in \Nmath,x \in X_n \mapsto |E_n \cap s^{-1}(x)|$ is essentially bounded.
\end{defi}

Note that, if $(\lambda_n)_n$ is a probability measure  on $\Nmath$ with full support, a coarse structure on a sequence $(\cg_n)$ is the same notion as a coarse structure on the \pmp groupoid $\cup_n \cg_n$ over the disjoint union $\cup_n X_n$ equipped with the probability measure $\sum_n \lambda_n \mu_n$. Moreover the notions of bounded geometry agree.

We will almost always assume that the coarse structures are \defin{countably generated}, i.e., that there is a sequence  of controlled sets $E_k$ such that every controlled set is contained in one of the $E_k$ (respectively there is a sequence $(E_{k,n})_n \in \coarse$ such that every controlled sequence $(E_n)$ is contained in one of the $(E_{k,n})_n$).

\begin{exam}\label{ex:coarse_structure_on_graphs} A sequence $(G_n)_n$ of finite graphs with vertex set $X_n$ defines a coarse structure on the sequence of groupoids $X_n \times X_n$ (the transitive equivalence relation on the set $X_n$ equipped with the uniform probability measure), by saying that $(E_n)_n$ is controlled if $\sup_n \sup_{(x,y) \in E_n} d(x,y)< \infty$. This coarse structure has bounded geometry if and only if the sequence $(G_n)_n$ has (uniformly) bounded degree.
\end{exam}

The main source of coarse structures on a \pmp groupoid are graphings, that we now define.

\begin{defi}[Graphings]\label{dfn:graphing}
  Let $\cg$ be a \pmp groupoid over the probability space $(X,\mu)$. A \defin{graphing} of $\cg$ is an at most countable ordered collection of bisections $\Phi=(\ph_i)_i$ indexed by $i \in \Nmath$ or $i \in \{1,\dots,N\}$ for an integer $N$. A graphing is \defin{generating} if the smallest groupoid which contains all the bisections $\ph_i$ is $\cg$ itself.
By a \defin{\pmp graphed groupoid} $(\cg,\Phi)$ we mean a \pmp groupoid $\cg$ for which a generating graphing $\Phi$ has been prescribed. 
\end{defi}
The \defin{size} of $\Phi$, denoted by $|\Phi|$, is the number of bisections which compose $\Phi$.
A \defin{$\Phi$-word of length $k$} is a product of $k$ elements of $\Phi$ and $\Phi^{-1}$.
The set of $\Phi$-words of length at most $k$ is denoted by $\Phi^k$.
We will also denote by $\overline \Phi$ the subset of $\cg$ composed of all elements of elements in $\Phi$, that is $\overline \Phi=\{g\in \cg:\ g\in\ph\in\Phi\}$. Observe that $\overline{\Phi^k}=\overline{\Phi}^k$.
A generating graphing $\Phi=(\ph_i)_i$ of a \pmp groupoid $\cg$ induces also a \defin{length} as follows: 
\begin{equation}\label{eq:length}
  \ell_\Phi:\cg\rightarrow \Nmath,\ g\mapsto \min\left\lbrace k\in\Nmath: g\in \overline{\{\ph_1,\ldots,\ph_k\}}^k\right\rbrace.
\end{equation}
 The empty word (the $\Phi$-word of length $0$) is the bisection $X$ (the space of units).
 The length of every unit is $0$. We will denote by $B_{\ell_\Phi}(L)$ the ball of length $L$: $B_{\ell_\Phi}(L) := \{g \in \cg: \ell_\Phi(g) \leq L\}$.

\begin{rema}[Coarse structures vs graphings]
\label{remark:coarse_structure}
  If $(\cg,\Phi)$ is a graphed \pmp groupoid, then the collection of measurable subsets on which $\ell_\Phi$ is bounded is a countably generated coarse structure with bounded geometry: the \defin{coarse structure associated to $\Phi$}. Similarly, if $(\cg_n,\Phi_n)$ is a sequence of graphed groupoids, then the collection of all the sequences $(E_n \subset \cg_n)_n$ satisfying $\sup_{n \in \Nmath, g_n \in E_n} \ell_{\Phi_n}(g_n)<\infty$ is a countably generated coarse structure with bounded geometry on $(\cg_n)$.

  The precise choice of the length function \eqref{eq:length} is not important, we could as well have chosen a sequence $(a_i)_i$  of positive real numbers tending to $+\infty$, and define instead the shortest-path length obtained by assigning to every $h \in \ph_i$ the length $a_i$:
\begin{equation*}
 g\mapsto \min\left\lbrace \sum_{j=1}^k a_{i_j} : k \in \Nmath, i \in \Nmath^k, \varepsilon \in \{-1,1\}^k, g \in \ph_{i_1}^{\varepsilon_1}\ph_{i_2}^{\varepsilon_2}\dots \ph_{i_k}^{\varepsilon_k}  \right\rbrace.
\end{equation*}

Different choices for $a_i$ give Lipschitz equivalent lengths and they all give the same coarse structure as \eqref{eq:length} on (sequences of) graphed groupoids. 

Conversely, the reader can check that any countably generated coarse structure with bounded geometry on a \pmp groupoid is the coarse structure associated to a graphing. Such a graphing is generating if and only if the coarse structure is.
\end{rema}
\begin{exam}
  A group is naturally a groupoid, a (generating) graphing is a (generating) subset, and when $G$ is finitely generated by $S$, $\Phi=S$ and $a_i=1$ for the elements of $S$, then the above length is the usual word length.
  On the contrary, when the generating set $S$ is infinite, then $\ell_{\Phi}$ is not bounded on $S$ and thus $S$
  it is not controlled for the coarse structure considered in Remark~\ref{remark:coarse_structure}.
\end{exam} 

\subsection{Actions}\label{sec:actions}
\begin{defi}
A (measurable) 
\defin{action} $\theta$ of a \pmp groupoid $\cg$ on the fibred space $(\cf,\pi)$ (both over the probability space $(X,\mu)$) is a measurable map $\theta: (\cg,s) \ast (\cf,\pi) \to \cf, (g,f)\mapsto \theta_g f$ such that:
\begin{itemize}
\item $\theta_g$ is bijection from $\pi^{-1}(s(g))$ to $\pi^{-1}(t(g))$,
\item $\theta_{gh}=\theta_g\theta_h$ whenever $t(h)=s(g)$.
\end{itemize}
\end{defi}
Besides its standard fibred structure $\pi_s:(g,f)\mapsto s(g)=\pi(f)$, 
the fibred product $(\cg,s) \ast (\cf,\pi)$ admits also the target fibred structure $\pi_t:(g,f)\mapsto t(g)$,
and $\theta$ is a morphism of fibred spaces $\theta:((\cg,s) \ast (\cf,\pi), \pi_t) \to (\cf,\pi)$
(i.e. $\pi_t(g,f)=\pi(\theta_g f)$).
As for actions of groups, we will often forget $\theta$ and we will write $gf$ for $\theta_g f$. 
If $\mathcal E \subset \cg$ and 
$\mathcal D \subset \cf$ are measurable subsets, we denote by $\mathcal E \mathcal D$ the image of the restriction of the action to $\mathcal E$ and $\mathcal D$, i.e., the image of $\theta:(\mathcal E,s)*(\mathcal D, \pi)\rightarrow \cf, (g,f)\mapsto \theta_g f$. Observe that $\mathcal E \mathcal D$ is measurable (same proof as Remark \ref{rem: X is measurable in G}).
The \defin{$\cg$-orbit} of $\mathcal D \subset \cf$ is $\cg \mathcal D$.
A \defin{fundamental domain} for the action of $\cg$ on $\cf$ is a measurable subspace $\mathcal D \subset \cf$ such that the action $\theta$
induces an isomorphism of fibred spaces 
$\left(\left(\cg,s\right)*\left(\mathcal D, \pi\right), \pi_t\right)\simeq (\cf, \pi)$
(for the target fibration $\pi_t$).

By a $\cg$-\defin{fibred space}, we mean a fibred space $\cf$ equipped with a prescribed $\cg$-action.

\begin{rema}
\label{ex:quotient fibred}
  Let $\cg$ be a \pmp groupoid and let $\cf$ be a $\cg$-fibred space. Then we can define a \defin{quotient} $\cf/\cg$ to be the quotient of $\cf$ by the equivalence relation induced by the action of $\cg$. If there exists a fundamental domain $\cd\subset \cf$, then this quotient is naturally isomorphic to $\cd$. In particular it has a measurable structure and it is itself a fibred space. Observe also that the induced quotient map from $\cf$ to $\cf/\cg$ is a measurable map of fibred spaces.
\end{rema}

\begin{exam}
\label{ex:subgroupoid}
  A \pmp subgroupoid $\ch$ of a \pmp groupoid $\cg$ over $X$ is a measurable subset $\ch\subset \cg$ which is also a groupoid over $X$ with respect to the restriction of the product, inverse, source and target. Clearly $\cg$ is a subgroupoid of itself. 

If $\ch<\cg$ is a \pmp subgroupoid, then $\ch$ carries two actions on $\cg$:  the \defin{action by left multiplication} $\theta^l$ on $(\cg, t)$ given by $\theta^l_g h=gh$, and the  \defin{action by right multiplication} $\theta^r$ on $(\cg, s)$ given by $\theta^r_g h=hg^{-1}$. Both these actions admit a fundamental domain. Let us show it for the action on the left. If $\{\ph_i\}_i$ are bisections such that $\cup_i\ph_i=\cg$, then we can define $\cd\subset \cg$ inductively as follows. Set $\cd_0:=X$, the set of units of $\cg$. Assume that $\cd_k$ is constructed for some $k\in\mathbb N$ such that for every bisection $\psi$ of $\ch$ we have that $\psi\cd_k\cap \cd_k=\emptyset$. Then by hypothesis $\ch\cd_k\subset \cg$ is measurable and take the first $i$ such that $\tilde \mu(\ph_i\setminus\ch\cd_k)\neq 0$. Put $\cd_{k+1}:=\cd_k\cup (\ph_i\setminus \ch\cd_k)$. It is easy to observe now that $\cd:=\cup_k \cd_k$ is a fundamental domain for the action on the left of $\ch$ on $\cg$.
\end{exam}

\begin{exam}
\label{ex:Cayley}
 If $(\cg,\Phi)$ is a \pmp graphed groupoid, then we can define the \defin{Cayley graph} $Cay(\cg,\Phi)$ of $\cg$ as follows. The set of vertices is $V:=\cg$ which is a fibred space over $X$ with respect to the source map and
 the set of edges is
 $E=\sqcup_{\varphi\in \Phi} \{(h,g):h\in\varphi, g\in \cg, t(g)=s(h)\}$ where the endpoints of the edge $(h,x)$ are $x$ and $h x$.  Then the set of edges $E$ is also a fibred space via $(h,x)\mapsto s(x)$ and $\cg$ acts on the graph $Cay(\cg,\Phi)$ by right multiplication. A fundamental domain of $E$ is given by  $\sqcup_{\varphi\in \Phi} \{(h,s(h)):h\in\varphi\}$.
\end{exam}

\begin{exam}\label{example:fibredProductAction}
  If $\cf_1,\cf_2$ are $\cg$-fibred spaces over $X$, then so is $\cf_1\ast\cf_2$, for the diagonal action $g(f_1,f_2)=(gf_1,gf_2)$.
  \end{exam}

Let us assume that the \pmp graphed groupoid $\cg$ acts on the fibred space $\cf$ and let $\cd\subset\cf$ be a fundamental domain. Similarly to (\ref{eq:length}), we can define a length function $\ell_{\Phi,\cd}$ on $\cf$ by the formula
\begin{equation}\label{eq:lengthact}
\ell_{\Phi,\cd}:\cf\rightarrow \Nmath,\ f\mapsto \min\left\lbrace k\in\Nmath: f\in \overline{\{\ph_1,\dots,\ph_k\}}^k\cd\right\rbrace.  
\end{equation}

By definition we have $\ell_\Phi=\ell_{\Phi,X}$ both for the left and right actions of $\cg$ on itself.

More generally, if $\coarse$ is a coarse structure on a \pmp groupoid $\cg$ acting on a fibred space $\cf$ with a fundamental domain $\cd$, the collection $\coarse_* \cd := \{ E \cd: E \in \coarse\}$ of measurable subsets of $\cf$ defines a notion of bounded sets in $\cf$. Similarly, if $\cg_n$ acts on $\cf_n$ with a fundamental domain $\cd_n$, and $\coarse$ is a coarse structure on $(\cg_n)_n$, then $\coarse_* (\cd_n) := \{ (E_n \cd_n)_n: (E_n)_n \in \coarse\}$ defines a notion of bounded measurable sets on the disjoint unions $\cup_n \cf_n$. When $\coarse$ is given by a (sequence of) graphing(s), then these notions of boundedness coincide with the boundedness according to the length \eqref{eq:lengthact}.

\subsection{Presentation of \pmp groupoids}
\label{subsect: presentation of groupoids}

The notion of presentation of groupoids has been introduced in \cite{Alv08}.

\begin{defi}\label{def:groupoid_homomorphism}
  A \defin{homomorphism} of \pmp groupoids $T:\cg\rightarrow \ch$ is a measurable map commuting with source, target, products and inversion. In particular, it sends the units of $\cg$ to the units of $\ch$, since $X=\{g\in \cg: s(g) = t(g),g^2=g\}$. We say that a homomorphism is \defin{\pmp}if the map between the unit spaces is measure preserving.

 A homomorphism of groupoids is \defin{units-preserving} if the induced map at the level of units is an isomorphism of probability spaces. In particular a units-preserving homomorphism is always a \pmp homomorphism.
\end{defi}
A groupoid $\cn$ is \defin{totally isotropic} if every element $n\in \cn$ is isotropic: $s(n)=t(n)$. As an example every groupoid $\cg$ admits the totally isotropic subgroupoid consisting of all its isotropic elements. 
Observe however that such a subgroupoid is not \pmp in general, see Example~\ref{ex: example equiv-rel pmp act non meas.}. However it is the case whenever $\cg$ is \textit{realizable}, see Definition~\ref{dfn:realizable}.

A totally isotropic subgroupoid $\cn<\cg$ is \defin{normal} if for every $g\in \cg$ and $n\in \cn_{s(g)}$, we have $g n g^{-1}\in \cn_{t(g)}$. As for groups, normal subgroupoids are designed for forming quotients.

\begin{prop}
\label{prop: quot groupoid}
  Let $\cg$ be a \pmp groupoid and let $\cn<\cg$ be a totally isotropic normal \pmp subgroupoid. Then the quotient $\cg/\cn$ is a \pmp grou\-poid and the quotient map $Q:\cg\rightarrow \cg/\cn$ is a \pmp and units-preserving homomorphism of \pmp groupoids. 
\end{prop}
\begin{proof}
  We have already observed in Example~\ref{ex:subgroupoid} that if $\cn<\cg$ is a \pmp subgroupoid, then the left action of $\cn$ on $\cg$ admits a fundamental domain $\cd\subset \cg$. If $\cn$ is a normal totally isotropic groupoid, then this fundamental domain is also a fundamental domain for the action on the right. Indeed, since $\cg = \cn \cd$, every element $g \in \cg$ can be written as $g = nd = d (d^{-1} n d)$ so $\cg = \cd \cn$. Moreover
  if $g,g'\in \cd$ and $n,n'\in \cn$ are such that $gn=g'n'$, then $gng^{-1}g=g' n'g'^{-1}g'$ and hence $g=g'$ and $gng^{-1}=g'n'g'^{-1}$ which implies that $n=n'$. 
This implies that the quotient $\cg/\cn$ is naturally a fibred space with respect to $s(g\cn)=s(g)$ and $t(g\cn)=t(g)$. We can define the multiplication as usual, if $g\cn,h\cn\in\cg/\cn$ are such that $s(g)=t(h)$, then $g\cn h\cn=gh\cn$ and this is well defined. The inverse can be defined analogously and it is a routine check that  these operations turn $\cg/\cn$ into a \pmp groupoid and $Q$ is a measurable map between \pmp groupoids.
\end{proof}

The  argument of the proof above gives also the following.

\begin{exam}[Quotient fibred space action]
\label{ex: quotient fibred space}
If $\cf$ is a $\cg$-fibred space and $\cn<\cg$ is a totally isotropic normal measurable subgroupoid, then the quotient \pmp groupoid $\cg/\cn$ acts measurably on the fibred quotient space $\cn\backslash \cf$.
\end{exam}

As an example we can consider a \textit{realizable} groupoid $\cg$ (see Definition~\ref{dfn:realizable}) and its totally isotropic subgroupoid $\cn$ consisting of all the isotropy elements of $\cg$. It is normal and the quotient $\cg/\cn$ is the associated equivalence relation (principal groupoid) on the unit space $X$ of $\cg$.

If $\cn<\cg$ is a totally isotropic subgroupoid the smallest normal subgroupoid of $\cg$ containing $\cn$ is denoted by $\langle\cn\rangle^\cg <\cg$. It remains totally isotropic.
Observe that there are countably many bisections $\{\ph_i\}_i$ such that $\cg=\cup_i \ph_i$ and countably many bisections $\{\psi_j\}_j$ such that $\cn=\cup_j \psi_j$. Then $\langle\cn\rangle^\cg$ is the subgroupoid generated by $\{\psi_j\ph_i\psi_j^{-1}\}_{i,j}$ and therefore we have that $\langle\cn\rangle^\cg$ is a \pmp subgroupoid of $\cg$. 

Assume now that $\cg$ and $\ch$ are \pmp groupoids and let $T:\cg\rightarrow \ch$ be a \pmp and units-preserving homomorphism of \pmp groupoids. Then the set $\cn:=\{g\in \cg:\ T(g)\text{ is a unit}\}$ is a normal subgroupoid of $\cg$ and moreover $\cg/\cn$ is canonically isomorphic to $\ch$.

A \defin{generating system} $\Psi=(\psi_i)_{i\in I}$ is an at most countable family of partial isomorphisms $\psi_i:\dom(\psi_i)\to \tar(\psi_i)$ of $(X,\rb, \mu)$. The \defin{free \pmp groupoid} $\cf(\Psi)$ over $X$ on $\Psi$ is defined as follows:
\begin{equation*}
\cf(\Psi):=\left\{(w,x): w \text{ is a reduced $\Psi$-word  and } x\in \dom(w) \right \}
\end{equation*}
 with the space of units associated with the empty word:
$X \simeq  \{\emptyset\} \times  X$; with source and target maps $s(w,x):=x$ and $t(w,x):=w(x)$; with product $(w',x')(w,x):=(w'w,x)$ (where $w'w$ stands for its reduced form) whenever the condition $x'=s(w',x')=t(w,x)=w(x)$ holds and with inverse 
$(w,x)^{-1}:=(w^{-1}, w(x))$.
Pick an enumeration $(w_i)_{i\in \Nmath}$ of the reduced $\Psi$-words; consider the two injections 
$J_s:\cf(\Psi)\to X\times \Nmath$, $(w_i, x)\mapsto (x,i)$ and  
$J_t:\cf(\Psi)\to X\times \Nmath$, $(w_i, x)\mapsto (w_i(x),i)$
and observe that they induce the same measurable structure on $\cf(\Psi)$ (the map $(x,i)\mapsto (w_i(x),i)$ is a bimeasurable bijection between the images).

\begin{rema}
  If $\Phi$ is a generating graphing of a \pmp groupoid $\cg$, then the associated partial isomorphisms $\hat \Phi:=(\hat \ph_i)_i$ (see Definition~\ref{def: meas. bisection}) is a generating system. 
    Moreover mapping each bisection $\hat \ph_i$ of $\cf(\hat\Phi)$ to the associated bisection $\ph_i$ of $\cg$ extends to a \pmp and units-preserving homomorphism $T_\Phi:\cf(\hat\Phi)\rightarrow \cg$.
This map is clearly surjective and we denote its kernel by $\cn_\Phi\subset \cf(\hat\Phi)$: it is a normal totally isotropic groupoid. 
\end{rema}

A \defin{relator system} $R=(R_w)_w$ of a generating system $\Psi$ is a collection, indexed by the reduced $\Psi$-words $w$, of measurable subsets of their fixed points sets $R_w\subset \{x\in \dom(w): w(x)=x\}$. Observe that a relator system gives rise to a totally isotropic subgroupoid of $\cf(\Psi)$, 
namely the subgroupoid generated by the restrictions $w\bigr|_{R_w}$.
We denote by $\cn_R^{\Psi}$ the (totally isotropic) normal subgroupoid in $\cf(\Psi)$ it generates.

\begin{defi} 
  A \defin{presentation} of a \pmp groupoid $\cg$ is given by a 
  \\
  -- generating system $\Psi=(\psi_i)_{i\in I}$ on $(X,\rb, \mu)$
\\
-- a relator system $R$ of $\Psi$ 
\\
--  a \pmp and units-preserving isomorphism $T_\Phi: \cf(\Psi)/\cn_R^{\Psi}\rightarrow \cg$
\end{defi}

\subsection{Ultraproducts}

From now on we will denote by $\ul$ a non-principal ultrafilter on $\Nmath$. We say that a property $P(n)$ holds $\ul$-almost surely (denoted $P(n)$ for $\ulae$) if $\{n : P(n)\}$ belongs to $\ul$.
\begin{defi}
  Let $(X_n)_{n\in\Nmath}$ be a sequence of sets. The \defin{ultraproduct} of the sequence $(X_n)_n$ is the set \[X_\ul:=\left(\prod_{n\in\Nmath} X_n\right)/\sim_\ul\quad\text{ where }\quad (x_n)_n \sim_\ul (y_n)_n\ \text{ if }\ x_{n}= y_{n}\textrm{ for }\ulae.\]
\end{defi}

We denote by $x_\ul$ and $A_\ul$ respectively elements and subsets of $X_\ul$. For a sequence $(x_n)_n\in \prod_{n\in\Nmath} X_n$, we denote by $[x_n]_\ul$ its class in $X_\ul$ and similarly for a sequence of subsets $A_n\subset X_n$, we denote its class by $[A_n]_\ul$ (i.e., the image of $\prod_n A_n$ in the quotient $X_\ul$). The construction of ultraproducts of probability spaces is due to Loeb, see \cite{CKTD} or \cite{Card-ultra-prod-we-sofic-arxiv} for a more modern treatment.

\begin{theo}\label{thm:ultraproduct_measure_space} For every $n\in\Nmath$ let $(X_n,\mu_n)$ be a probability space and let $X_\ul$ be the ultraproduct of the sequence $(X_n)_n$. Then there exists a probability measure $\mu_\ul$ on $X_\ul$ such that
  \begin{enumerate}
  \item for every sequence of measurable subsets $(A_n\subset X_n)_n$ the set $[A_n]_\ul$ is $\mu_\ul$-measurable and $\mu_\ul([A_n]_\ul)=\lim_{n\in\ul}\mu_n(A_n)$,
  \item\label{item:measurable_subsets_in_ultraproducts} for every $\mu_\ul$-measurable subset $A_\ul\subset X_\ul$ there exists a sequence of measurable sets $(A_n\subset X_n)_n$ such that $\mu_\ul(A_\ul\Delta [A_n]_\ul)=0$.
  \end{enumerate}
\end{theo}

\begin{defi}\label{dfn:part_iso}
If $\ph_n$ is a partial isomorphisms of $X_n$ for every $n$, then the map $[\ph_n]_\ul:[x_n]_\ul\mapsto [\ph_n(x_n)]_\ul$ is also a partial isomorphism from $[\dom(\ph_n)_\ul]$ to $[\tar(\ph_n)]_\ul$, called the \defin{ultraproduct partial isomorphism}.
\end{defi}

We now define the ultraproducts of \pmp groupoids with a coarse structure.

\begin{prop}[Ultraproducts of coarse \pmp groupoids] 
\label{prop:ultraproduct_of_coarse_structures}
Let  $(\cg_n)_n$ be a sequence of \pmp groupoids on $X_n$.

1--   Let $\coarse$ be a countably generated coarse structure with bounded geometry on $(\cg_n)_n$.
 There is a unique structure of \pmp groupoid over $X_{\ul}$ on
  \[ \cg_\ul := \{(g_n)_n : \exists (E_n)_n \in \coarse, g_n \in E_n\}/\sim_\ul\]
  for the groupoid operations $s([g_n]_\ul) := [s(g_n)]_\ul$, $t([g_n]_\ul) := [t(g_n)]_\ul$, $[g_n]_\ul[h_n]_\ul := [g_n h_n]_\ul$ and $[g_n]_\ul^{-1} := [g_n^{-1}]_\ul$ and for which $\coarse_\ul := \{[E_n]_\ul: (E_n)_n \in \coarse\}$ is a countably generated and generating coarse structure with bounded geometry.

2--
  If $\coarse$ is the coarse structure associated to a sequence of graphings $\Phi_n=(\ph_n^i)_i$, then
  \[ \cg_\ul = \{(g_n)_n : \sup_n \ell_{\Phi_n}(g_n)<\infty\}/\sim_\ul\]
  and $\coarse_\ul$ is the coarse structure associated to the graphing $\Phi_\ul=([\ph_n^i]_\ul)_i$ and this graphing is generating.
\end{prop}

\begin{defi}
The \pmp groupoid $\cg_\ul$ defined in Proposition~\ref{prop:ultraproduct_of_coarse_structures}
is called the \defin{ultraproduct groupoid} of the coarse  sequence $((\cg_n)_n, \coarse)$.
\end{defi}

In the simplest case of bounded degree graphs, we recover the ultraproduct graphs considered for example in \cite[Section 9]{HatLovSze}:
\begin{exam}[Ultraproduct of bounded degree graphs]\label{ex:ultraproduct_of_graphs}
Consider a sequence  $(G_n=(X_n,\mathrm{Edg}_n))_n$ of bounded degree finite graphs. Then the ultraproduct of the groupoids $X_n \times X_n$ for the coarse structure considered in Example~\ref{ex:coarse_structure_on_graphs} is the \pmp equivalence relation $\RR_\ul$ on $X_\ul$ given by $[x_n]_\ul \sim [y_n]_\ul$ if and only if $\lim_\ul d(x_n,y_n)<\infty$. Moreover, the coarse structure on the ultraproduct is generated by
  \[ \mathrm{Edg}_\ul := [\mathrm{Edg}_n]_\ul = \{( [x_n]_\ul,[y_n]_\ul) : \lim_\ul d(x_n,y_n)=1\}.\]  
  \end{exam}

\begin{proof}[Proof of Proposition~\ref{prop:ultraproduct_of_coarse_structures}] The fact that the operations are well defined and turn $\cg_\ul$ into a groupoid over $X_\ul$ is clear. The issue is to define the measure space structure on $\cg_\ul$. For every $(E_n)_n\in \coarse$, we can consider the ultraproduct measure space $[E_n]_\ul \subset \cg_\ul$. By our assumption of bounded geometry, it is a finite height fibred space over $X_\ul$, for both maps $s$ and $t$. Moreover if $(F_n)_n \in \coarse$ contains $(E_n)_n$, then this measure space structure on $[E_n]_\ul$ coincides with the restriction of the measure space structure on $[F_n]_\ul$. Taking a sequence $((E_{k,n})_n)_k \subset \coarse$ generating $\coarse$, we obtain a $\sigma$-algebra and a measure on $\cg_{\ul}$ for which $s$ and $t$ define fibred space structures, and for which every $(E_n)_n \in \coarse$ defines a measurable subset $[E_n]_\ul$ on $\cg_\ul$. The fact that $\{[E_n]_\ul:(E_n) \in \coarse\}$ is a coarse structure is immediate from (\ref{item:measurable_subsets_in_ultraproducts}) in Theorem~\ref{thm:ultraproduct_measure_space}. The measurability of the inverse and product follow from the measurability of these operations on $\cg_n$ and the fact that $\coarse$ is stable under these operations. This coarse structure is generating by definition of $\cg_\ul$, and is countably generated because $\coarse$ was assumed to be countably generated.

The verification that $\coarse_\ul$ coincides with the coarse structure defined by $\Phi_\ul$ is routine.
\end{proof}

\begin{defi} 
Similarly to Definition~\ref{dfn:part_iso} we can define the notion of ultraproduct of bisections. Assume that $(E_n)_n \in \coarse$ and for every $n$ take a bisection $\ph_n$ contained in $E_n$, then $[\ph_n]_\ul\subset \cg_\ul$ is a bisection. 
\end{defi}

\begin{coro}\label{cor:dependance_en_graphage} If two sequences of generating graphings $\Phi_n$ and $\Psi_n$ of  \pmp groupoids $\cg_n$ are \defin{$\ul$-coarsely equivalent} (i.e., for every sequence $g_n \in \cg_n$, $\lim_\ul \ell_{\Phi_n}(g_n) <\infty$ if and only if $\lim_\ul \ell_{\Psi_n}(g_n) <\infty$), then the ultraproduct graphings $\Phi_\ul$ and $\Psi_\ul$ define the same groupoid with coarse structure.
\end{coro}
\begin{proof}
  The assumption of coarse equivalence means that for every $M$, there is $M'$ such that, for $\ulae$, $\ell_{\Phi_n}(g) \leq M \implies \ell_{\Psi_n}(g) \leq M'$, and conversely $\ell_{\Psi_n}(g) \leq M \implies \ell_{\Phi_n}(g) \leq M'$. Indeed, if for every $n$, we define $M_n = \sup\{ \ell_{\Psi_n}(g): g \in \cg_n, \ell_{\Phi_n}(g) \leq M\}$, and pick $g_n \in \cg_n$ with $\ell_{\Phi_n}(g_n) \leq M$ and $\ell_{\Psi_n}(g) = M_n$ if $M_n<\infty$ and $\ell_{\Psi_n}(g)\geq n$ otherwise, we obtain that $M':= \lim_\ul M_n <\infty$. This allows us to identify the groupoids
  \[\{ (g_n)_n : \sup_n \ell_{\Phi_n}(g_n) <\infty\}/\sim_\ul\]
  and
  \[\{ (g_n)_n : \sup_n \ell_{\Psi_n}(g_n) <\infty\}/\sim_\ul,\]
and Proposition~\ref{prop:ultraproduct_of_coarse_structures} shows that this identification is measurable and measure preserving.
\end{proof}

\begin{rema}
In Corollary~\ref{cor:dependance_en_graphage}, it is crucial that the two graphings are $\ul$-coarsely equivalent: sequences of different graphings can give rise to rather unrelated limits. For instance, the transitive equivalence relation on a finite set has many different generating graphings (all connected graphs with this finite set as vertex set). However their ultraproducts can be completely different according, for instance, to whether one takes a sequence of expanders, in which case the ultraproduct equivalence relation is strongly ergodic, or a sequence of Hamiltonian cycles in which case the ultraproduct is not ergodic.
\end{rema}

\begin{rema}\label{rem: order of the graphings and ultra-prod}
The structure of a graphed groupoid does not fundamentally depend on the order of the bisections of the graphing. However, the order becomes crucial when we consider a sequence of $\Phi_n=(\ph_i)_{i\in I_n}$ and its ultraproduct. 
But when the size $(|\Phi_n|)_n$ is bounded, i.e., the $I_n$ are contained in a finite subset $J\subset \Nmath$, then a sequence $(\sigma_n)_n$ of permutations of $J$ is indeed essentially constant under the ultrafilter $\ul$ (indeed the underlying coarse structure is canonical).
It follows  in this case that the ultraproduct groupoid does not depend on the order and the ultraproduct graphing is unique up to a single permutation of the indices.
The same holds when instead of assuming $J$ finite, we assume that the permutations $\sigma_n$ converges to a \emph{permutation} of $\Nmath$, i.e., if for every $j$, $\lim_\ul \sigma_n(j)<\infty$ and $\lim_\ul \sigma_n^{-1}(j)<\infty$ (for example if $\sup_n \sigma_n(j)<\infty$ and $\sup_n \sigma_n^{-1}(j)<\infty$). This is a particular case of Corollary~\ref{cor:dependance_en_graphage}.
\end{rema}

We can now give a similar construction of ultraproduct for $\cg$-fibred spaces. 
\begin{defi}\label{def:ultraaction}
  Let us assume that for every $n\in\Nmath$ the \pmp groupoid $\cg_n$ acts on the fibred spaces $(\cf_n,\pi_n)$ and let $\cd_n\subset\cf_n$ be a fundamental domain such that the sequence $n\mapsto H(\cd_n)$ is bounded\footnote{Recall that $H(\cd_n)$ is the height of $\cd_n$: the essential supremum of the function $x\mapsto |\pi_n^{-1}(x)\cap \cd_n|$.}. Let $\coarse$ be a countably generated coarse structure with bounded geometry on $(\cg_n)_n$, for example the coarse structure coming from a sequence of graphings. Set
\begin{gather*}
  \cf_\ul:=\left\lbrace (f_n)_n\in \prod\cf_n:\ \exists (E_n)_n \in \coarse \textrm{ such that }f_n \in E_n D_n\textrm{ for every }n\right\rbrace/\sim_\ul\\
  (f_n)_n\sim_\ul (f'_n)_n\ \text{ if }\ f_{n}= f'_{n}\textrm{ for }\ulae.
\end{gather*}

We also set $\pi_\ul([f_n]_\ul):=[\pi_n(f_n)]_\ul$. Observe that for every $(E_n)_n \in \coarse$, the fibred spaces $E_n \cd_n$ over $\pi_n(E_n \cd_n)$ have bounded height, and therefore their ultraproduct $[E_n \cd_n]_\ul$ is a fibred space over $[\pi_n(E_n \cd_n)]_\ul$. 
So as for the groupoids we equip $\cf_\ul$ with the inductive limit measure space structure which turns $(\cf_\ul,\pi_\ul)$ into a fibred space, and on which the groupoid $\cg_\ul$ acts by the formula $[g_n]_\ul[f_n]_\ul:=[g_n f_n]_\ul$. Observe that a fundamental domain is given by $\cd_\ul$, the ultraproduct of the sequence $(\cd_n)_n$, whose height is $\lim_\ul H(\cd_n)$. 
\end{defi}

Let us observe that the definition of the ultraproduct does depend on the sequence of fundamental domains. One can also check that if $\cf_n$ is the Cayley graph of $\cg_n$ with respect to $\Phi_n$ as explained in Example~\ref{ex:Cayley}, then the ultraproduct of these fibred spaces is the Cayley graph of $\cg_\ul$ with respect to $\Phi_\ul$. 

Consider the setting of Proposition
\ref{prop:ultraproduct_of_coarse_structures}. If $\ph_n \in [[\cg_n]]$
is a controlled sequence of measurable bisections, on the one hand by
Proposition~\ref{prop:ultraproduct_of_coarse_structures} we can regard
$[\ph_n]_\ul$ as an element of $[[\cg_\ul]]$. On the other hand since
the full pseudogroups of the groupoids are metric spaces we can also
consider the image of $(\ph_n)_n$ in the metric ultraproduct
$\prod_\uu [[ \cg_n]]$ of the sequence $([[\cg_n]])_n$. In order to
avoid confusions, we will for a moment denote by $[\ph_n]_\ul^d$ the
element of $\prod_\uu [[ \cg_n]]$ that represents the sequence of
elements $\ph_n\in [[\cg_n]]$. The following lemma follows easily from
our construction.

\begin{lemm}\label{lem:ultraproduct_of_pseudo_full_groups} 
  There is a unique isometric embedding of $[[\cg_\ul]]$ into
  $\prod_\uu [[ \cg_n]]$ sending $[\ph_n]_\ul$ to $[\ph_n]_\ul^d$ for
  every controlled sequence of bisections. An element
  $[\ph_n]_\ul^d\in \prod_\uu [[ \cg_n]]$ lies in the image if and
  only if for every $\eps>0$ there exists a sequence $\psi_n \subset
  \ph_n$ with $d(\ph_n,\psi_n) \leq \eps$ for $\ulae$ and such that
  $(\psi_n)_n$ is controlled. Moreover the embedding extends to an
  embedding of the von Neumann algebra $\LL(\cg_\uu)$ into the von
  Neumann algebra ultraproduct $\prod_\uu \LL(\cg_n)$.
\end{lemm}
For a sequence $(\varphi_n)_n$ satisfying the condition in Lemma~\ref{lem:ultraproduct_of_pseudo_full_groups}, we will drop the notation $[\ph_n]_\ul^d$ and use only $[\ph_n]_\ul$ to denote both the element of $\prod_\uu [[ \cg_n]]$ and the corresponding element of $[[\cg_\ul]]$.

We have a similar situation for $\cg_n$-fibred spaces.

\begin{lemm}\label{lem:ultraproductactions}
   Let us assume that for every $n\in\Nmath$ the \pmp graphed groupoid $(\cg_n,\Phi_n)$ acts on the fibred space $(\cf_n,\pi_n)$ and let $\cd_n\subset\cf_n$ be a fundamental domain such that the sequence $n\mapsto H(\cd_n)$ is bounded. For every $n$ the full pseudogroup $[[\cg_n]]$ acts on the Hilbert space $\LLh^2(\cf_n)$. Then the Hilbert space $\LLh^2(\cf_\ul)$ embeds into the metric ultraproducts $\prod_\ul \LLh^2(\cf_n)$ and if $\ph_\ul=[\ph_n]_\ul$ is a measurable section of $\cf_\ul$ (where $\ph_n$ are measurable sections of $\cf_n$), then the image of $\chi_{\ph_\ul}$ is $[\chi_{\ph_n}]_\ul^d$. Moreover the above mentioned embedding of $[[\cg_\ul]]$ inside the metric ultraproduct $\prod_\uu [[ \cg_n]]$ is compatible with the action and the passage to von Neumann algebras. 
\end{lemm}

\subsection{Realizability and ultraproducts}
\label{subsect: realizability}

 \begin{defi}\label{dfn:realizable}
A \pmp groupoid  $\cg$ over the probability space $(X,\rb,\mu)$ is said to be \defin{realizable} if every measurably isotropic bisection $\ph$ is totally isotropic (Definition~\ref{def: tot. vs meas isotrop, free}).
\end{defi}

In other words $\cg$ is realizable if and only if every totally free bisection is measurably free or equivalently if and only if $\Fix_{\rb}({\ph})=\Fix ({\ph}):= \{g\in \ph: s(g)=t(g)\}$.  The following is classic.

\begin{lemm}
  Every \pmp groupoid on a standard Borel probability space is realizable.
\end{lemm}

Realizable groupoids are the groupoids that can be studied point-wise and those groupoids that have a \textit{free factor} onto a standard Borel space, see Theorem~\ref{thm: standard factors}. Therefore, up to some extent, realizable groupoids behave like standard one. 

\begin{lemm}[Realizability]
\label{lem:realizab ultralim and fix} 

The ultraproduct $(\cg_\ul, \Phi_\ul)$ of a sequence of realizable graphed \pmp groupoids $(\cg_n, \Phi_n)$ is itself realizable.
Moreover every bisection $\ph_\ul=[\ph_n]_\ul$ satisfies up to null sets
\[[\Fix_{\rb_n}({\ph_n})]_{\ul}= [\Fix({\ph_n})]_{\ul}= \Fix (\ph_{\ul})=\Fix_{\rb_\ul}(\ph_\ul)\] 
and $[M_F(\ph_n)]_\ul=M_F[(\ph_n)]_\ul$.
\end{lemm}
\begin{proof}
  First observe that since the ultraproduct groupoid $\cg_\ul$ is generated by the bisections of the form $[\ph_n]_\ul$ it is enough to prove that for every $\ph_\ul=[\ph_n]_\ul$ we have that $[\Fix({\ph_n})]_{\ul}=\Fix_{\rb_\ul}(\ph_\ul)$. This follows easily from Proposition~\ref{prop:manyinvo-gen}. Indeed consider a decomposition given by Proposition~\ref{prop:manyinvo-gen}, that is $M_F(\ph_n)=\sqcup_i A_n^i$ such that $\mu(\cup_{i=1}^k A_n^i) \geq (1-(2/3)^k) \mu(M_F(\ph_n))$ and  $\mu_n(A_n^i\cap \ph(A_n^i))=0$. Then clearly we get $[M_F(\ph_n)]_\ul=\cup_i [A_n^i]_\ul$ and $\mu_\ul(\ph_\ul([A_n^i]_\ul)\cap [A_n^i])=0$. In particular $M_F[(\ph_n)]_\ul\supset [M_F(\ph_n)]_\ul$. Therefore $\Fix_{\rb_\ul}(\ph_\ul)\subset \dom(\ph_\ul)\setminus [M_F(\ph_n)]_\ul$ which by hypothesis corresponds to $[\Fix({\ph_n})]_{\ul}$. On the other hand we clearly have that $[\Fix({\ph_n})]_{\ul}\subset \Fix(\ph_\ul)$ and hence $[\Fix({\ph_n})]_{\ul}=\Fix_{\rb_\ul}(\ph_\ul)$ as claimed.
\end{proof}

For a \pmp groupoid $\cg$ we will now denote by $\cg^P$ the associated equivalence relation. As we have already remarked, if $\cg$ is realizable, then $\cg^P$ is a \pmp equivalence relation. The following is easy to prove.

\begin{lemm}
  Given a sequence of realizable graphed \pmp groupoids $(\cg_n, \Phi_n)$ we have that the ultraproduct of $\cg_n^P$ is $\cg_\ul^P$. 
\end{lemm}

In general it is not clear what should be the size of a ``class'' or ``orbit'' in a \pmp groupoid, see Example~\ref{ex: example equiv-rel pmp act non meas.}. However for realizable groupoids the situation is much more clear and we have the following. For a groupoid $\cg$ over the space $X$ we will set $|x|_\cg:=|\{t(g):\ s(g)=x\}|$, the cardinality of the class of the associated equivalence relation. 

\begin{lemm}
\label{lem: ultralim has infinite orbits}
Let $(\cg_n, \Phi_n)$ be a sequence of realizable graphed groupoids such that $\lim_\ul|\Phi_n|$ is finite. Then for $\mu_{\ul}$-almost every $[x_n]_\ul\in X_\ul$ we have $|[x_n]_\ul|_{\cg_\ul}=\lim_\ul |x_n|_{\cg_n}$. In particular if for every $k\in \mathbb N$ we have $\lim_\ul\mu_{n}(\{x\in X_n: \vert x \vert_{\cg_n} \leq k\})= 0$, then for $\mu_{\ul}$-almost every $[x_n]_\ul\in X_\ul$ we have that $|[x_n]_\ul|_{\cg_\ul}$ is infinite.
\end{lemm}
\begin{proof}
  By the above lemma we can assume that all the $\cg_n$ are equivalence relations. For $x_n\in X_n$ set $|x_n|_{n,k}:=\vert \{t(g):\ s(g)=x,\ g\in \overline \Phi_n^k\}\vert$ and define similarly $|x_\ul|_{\ul,k}$. Remark that if for some $k$ we have that $|x_n|_{n,k}=|x_n|_{n,k+1}$, then for every $i$ we have that $|x_n|_{n,k}=|x_n|_{n,k+i}$. Observe also that Lemma~\ref{lem:realizab ultralim and fix} implies that  for $\mu_{\ul}$-almost every $[x_n]_\ul\in X_\ul$ we have $|[x_n]_\ul|_{\ul,k}=\lim_\ul |x_n|_{n,k}$ and therefore we can combine the two facts to obtain the desired result.
\end{proof}

 \section{Graphings, cost and combinatorial cost}
\label{sect: cost and comb cost}

The notion of cost of \pmp groupoids (on standard Borel spaces) has been considered by several authors (see for instance \cite{Ueda-cost-2006, Carderi-Master-thesis-2011, Abert-Nikolov-12, Takimoto-cost-L2-2005}).

Let $\cg$ be a \pmp groupoid. Following \cite{Lev95}, we define the \defin{cost} of the graphing $\Phi=(\ph_i)_i$ of $\cg$ as 
\begin{equation}
\cost(\Phi):=\sum_i \widetilde \mu(\ph_i)
\end{equation}
and the \defin{(groupoid) cost} of the \pmp groupoid $\cg$ as the infimum of the costs of its generating graphings:
\begin{equation}
\cost(\cg)=\inf_\Psi\cost(\Psi),
\end{equation} 
where $\Psi$ ranges among all the generating graphings of $\cg$. It is the infimum of the $\widetilde{\mu}$-measures of measurable generating subsets of $\cg$.

\begin{exam}
Let $\Gamma=\FF_p\times \FF_q$ the direct product of two free groups ($p, q\geq 2$) and $\alpha $ be a \pmp $\Gamma$-actions on standard Borel space where $\FF_p$ acts freely and $\FF_q$ trivially.
Then the groupoid $\cg_{\alpha}$ of the $\Gamma$-action has $\cost(\cg_{\alpha})=p$:
The underlying equivalence relation $\RR_{\alpha}=\RR_{\FF_p}$ has cost $p$ \cite{Gab-cost} and every generating graphing of $\cg_{\alpha}$ induces a generating graphing of $\RR_{\alpha}$. This gives $\cost(\cg_{\alpha})\geq p$. Taking the graphing 
$\Phi=(\ph_h)_{h}\cup (\ph_g)_g$ where $h$ (resp. $g$) runs over a set of standard generators of $\FF_p$ (resp. $\FF_q$) 
with $\ph_h:=h\times X$ and 
$\ph_g:=g\times A_g$ for $A_g$ an arbitrarily small $\FF_p$-complete section, one gets a generating graphing of cost close to $p$.

When the space is not standard, replace ``free action'' by ``measurably free'', ``trivial action'' by ``measurably isotropic'' and make use of  "the" standard factor from Theorem~\ref{thm: standard factors}.
\end{exam}

Recall that 
\\
-- $N<\Gamma$ is \defin{q-normal} if $\Gamma$ admits a generating set $S$ such that $s ^{-1}N s\cap N$ is infinite for each $s\in S$
\\
-- $N<\Gamma$ is \defin{weakly-q-normal} 
if  there exists a well ordered family of intermediate subgroups $N=N_0<N_1<\cdots <N_{\eta}=\Gamma$ such that $N'_j:=\cup_{i<j} N_i$ is q-normal in $N_j$.

\begin{prop}\label{prop:cost q-normal}
In the standard or realizable situation:
\\
-- If $N<\Lambda$ is \defin{q-normal} subgroup
and $\alpha$ is a \pmp $\Gamma$-action such that $N$ acts aperiodically,
then $\cost(\cg_{\alpha})\leq \cost(\cg_{\alpha\upharpoonright N})$.
\\
-- If $N<\Gamma$ is an infinite \defin{weakly-q-normal} subgroup
and $\alpha$ is a free \pmp $\Gamma$-action,
then $\cost(\RR_\Gamma)= \cost(\cg_{\alpha})\leq \cost(\cg_{\alpha\upharpoonright N})=\cost(\RR_N)$.
\end{prop}

\begin{proof}
Let $\Phi_0$ be a generating graphing of $\cg_{\alpha\upharpoonright N}$.
If $N<\Lambda$ is q-normal,  $\Lambda$ is generated by a family of $s\in S_\Lambda$ such that $\vert s^{-1} N s \cap N\vert=\infty$. Then $\Lambda$ admits a graphing $\Phi=\Phi_0\cup \Phi_{\Lambda}^*$ where  $\Phi_{\Lambda}^*$ has arbitrarily small cost: a family of bisections 
$\Phi_{\Lambda}^{*}=( s \times A_s)_{s\in S_\rho}$ where the $A_s$ are complete $\left(s^{-1} N  s \cap 
N\right)$-sections with $\sum_s \tilde{\mu}(A_s)\leq \epsilon$ for any prescribed $\epsilon>0$.
This generates since for every $x\in X$, there are $n_1, n_2\in N$ such that $n_1 x\in A_s$ and $s=n_2^{-1}s n_1$.
Then $(s,x)= (n_2^{-1}, sn_1 x). (s, n_1x). (n_1,x)$ and $(n_2^{-1}, sn_1 x), (n_1,x)$ belong to the groupoid generated by $\Phi_0$ while $(s, n_1x)\in s \times A_s$.

The proof for the weakly-q-normal case goes of course by a transfinite induction. 
The induction assumption is that $N_\rho$ admits a graphing $\Phi_\rho=\Phi_0\cup \Phi^*_\rho$ where  $\Phi^*_\rho$ can be chosen of arbitrarily small cost.

If $\rho=j+1$ is a successor, then $N'_\rho=N_j$ and $N_j$ is q-normal in $N_\rho$. Thus the groupoid of the $N_\rho$-action is generated by adding a small cost graphing $\Phi^{**}_\rho$ to $\Phi_j=\Phi_0\cup \Phi^*_j$. 

If $\rho$ is a limit ordinal, then up to enumerating $N'_\rho=\{\gamma_p\}_{p\in \Nmath}$, one selects a increasing sequence 
$(N_{j_p})_p$ (where $j_p$ is the smallest $j$ such that $\gamma_1, \cdots, \gamma_p\in N_j$) and the subgroupoid generated by the restriction of $\alpha$ to $N'_{\rho}$ is generated by $\Phi_\rho=\Phi_0\cup \Phi^*_{j_1}\cup \Phi^*_{j_2}\cup\cdots \cup \Phi^*_{j_p}\cup\cdots$ where $\sum_p \cost(\Phi^*_{j_p})$ can be chosen arbitrarily small by assumption.

Thus $N'_\rho$ being q-normal in $N_\rho$, its groupoid is generated by further adding a small cost graphing $\Phi^{**}_\rho$.

\end{proof}

\begin{coro}
\label{cor: cost 1 commensurated}
Groups with a commensurated fixed price $1$ subgroup $C$ have fixe price $1$.
\end{coro}

\subsection{Graphings with finite size}\label{subsection:finite_graphings}

For the sake of simplicity, we shall first only consider graphings of bounded size and generalize Elek's combinatorial cost in this setting. In that case, the length defined in (\ref{eq:length}) of Section~\ref{subsect : coarse structures and graphings} is bilipschitz-equivalent to the length $\ell'_\Phi(g) = \min\{k \in \Nmath: g \in \overline{\Phi}^k\}$: if $\Phi = (\ph_1, \ph_2, \cdots, \ph_D)$, we have $\ell'_\Phi \leq \ell_\Phi \leq D \ell'_\Phi$ . We will treat the general case in Section~\ref{subsection:infinite_graphings}.

\begin{defi}\label{def:Lipschitz-contain-equiv-graphins}
Let $\Phi$ and $\Psi$ be two graphings of finite size of the \pmp groupoid $\cg$ and $L$ be a natural number. We say that
\begin{enumerate}
\item 
$\Psi$ \defin{$L$-Lipschitz embeds in} $\Phi$ (or  $\Phi$ \defin{$L$-Lipschitz contains} $\Psi$), denoted by $\Psi\subset_L \Phi$, if every bisection $\psi_k\in\Psi$ is piecewise a restriction of words of length at most $L$ in the alphabet $\Phi$ and their inverses ($\widetilde{\mu}$-mod $0$);
\item the graphings $\Phi$ and $\Psi$ are \defin{$L$-Lipschitz
  equivalent}, denoted $\Phi\Lipeq_L\Psi$, if $\Psi\subset_L \Phi$ and
  $\Phi \subset_L \Psi$;
\end{enumerate}
\end{defi}

Two graphings of finite size $\Phi$ and $\Psi$ are said \defin{Lipschitz equivalent} (denoted $\Phi\Lipeq\Psi$) if they are $L$-Lipschitz equivalent for some $L \in \Nmath$. Clearly two graphings $\Phi$ and $\Psi$ are Lipschitz equivalent if and only if their length functions $\ell_\Phi$ and $\ell_\Psi$ on $\cg$ are Lipschitz equivalent. Observe that being Lipschitz equivalent is an equivalence relation (while being $L$-Lipschitz equivalent is not transitive). 

If a graphing $\Phi$ of the \pmp groupoid $\cg$ Lipschitz contains a generating graphing, then $\Phi$ itself is generating.

The $L$-\defin{Lipschitz cost} of the graphed groupoid $(\cg,\Phi)$ is the infimum of the costs of all graphings of $\cg$ which are $L$-Lipschitz equivalent to $\Phi$,

$$\cost_L(\cg,\Phi):=\inf\left\{\cost(\Psi)\colon \Psi\Lipeq_L\Phi \right\}.$$

\begin{rema}\label{rem:for_Lipschitz_cost} In the definition of the $L$-\defin{Lipschitz cost}, we can restrict to graphings $\Psi$ of size bounded by the number of words of length $\leq L$ in the elements of $\Phi$ and its inverse.
\end{rema}

  The following lemma claiming that the cost can be computed inside a given Lipschitz class of (finite cost) graphings was observed for equivalence relations in \cite[\'Etape 1]{Gab-CRAS-cost} and \cite{Gab-cost} (see also Lemma~\ref{lem:goodgraphcoarse}). We give a proof for completeness.   
\begin{lemm}[{Gaboriau \cite[Proposition IV.35]{Gab-cost}}]\label{lem:goodgraph}
Let $(\cg,\Phi)$ be a graphed groupoid where the graphing $\Phi$ has finite size.
The cost of $\cg$ can be computed inside the Lipschitz class of $\Phi$, i.e.,  
\begin{equation}
\cost(\cg)=\inf_L\cost_L(\cg,\Phi). 
\end{equation}
\end{lemm}
\begin{proof}
 Fix $\eps>0$. Let $\Psi$ be a generating graphing of $\cg$ such that $\cost(\Psi)\leq \cost(\cg)+\eps$. Take a positive integer $M$ such that $\tilde\mu(\overline\Phi\setminus\overline\Psi^M)\leq\eps$, where we recall that $\overline\Phi$ is the subset of $\cg$ given by the images of the bisections defining $\Phi$. For $L\in\Nmath$ we define \[\overline{\Psi_L}:=\overline\Psi\cap \overline{\Phi}^L\] and observe that $\cup_L \overline{\Psi_L}=\overline\Psi$ which implies that $\cup_L \overline{\Psi_L}^M=\overline\Psi^M$, thus there exists $L_0$ such that $\tilde\mu(\overline\Phi\setminus \overline{\Psi_{L_0}}^M)\leq 2\eps$. 
  There is a graphing of finite size $\Theta_0$ such that $\overline\Theta_0=\overline {\Psi_{L_0}}\cup(\overline\Phi\setminus \overline{\Psi_{L_0}}^M)$ 
  satisfying $\cost(\Theta_0)\leq \cost(\Psi)+2\eps\leq \cost(\cg)+3\eps$ and $\Theta_0\subset_{L_0}\Phi$. On the other hand, 
  $(\overline\Phi\setminus \overline{\Psi_{L_0}}^M)\subset_1 \Theta_0$ and $(\overline\Phi\cap \overline{\Psi_{L_0}}^M)\subset_{M} \Theta_0$ so that $\Phi\subset_M \Theta_0$.
  \end{proof}

\medskip
We now define the notion of $\ul$-combinatorial cost.
\begin{defi}[Elek \cite{Elek-combinatorial-cost-2007}]
\label{def: ul-comb-cost}
Let $(\cg_n,\Phi_n)_n$ be a sequence of graphed \pmp groupoids of bounded size. The $\ul$-\defin{combinatorial cost} of the sequence is
\begin{equation}\ccost_\ul((\cg_n,\Phi_n)_n):=\inf_L\lim_{n\in\ul}\cost_L(\cg_n,\Phi_n)\end{equation}
\end{defi}

\subsection{Graphings with infinite size}\label{subsection:infinite_graphings}
We consider the general case where we allow graphings with unbounded or infinite size.

\begin{defi}\label{def:Coarsely-contain-equiv-graphins}
Let $\Phi$ and $\Psi$ be two graphings of the \pmp groupoid $\cg$, and $M=(M_k)_{k \in \Nmath} \in \Nmath^\Nmath$ a sequence of natural numbers. We say that
\begin{enumerate}
\item 
$\Psi$ \defin{$M$-coarsely embeds in} $\Phi$ (or  $\Phi$ \defin{$M$-coarsely contains} $\Psi$), denoted by $\Psi\subset_{c,M} \Phi$,
if for every $k\in\Nmath$, the bisection $\psi_k\in\Psi$ is 
piecewise a restriction of words of length at most $M_k$ in the alphabet $\{\ph_1,\ldots,\ph_{M_k}\}$ and their inverses ($\widetilde{\mu}$-mod $0$);
\item the graphings
$\Phi$ and $\Psi$ are \defin{$M$-coarsely equivalent}, denoted 
$\Phi\ceq_M\Psi$, if $\Psi\subset_{c,M} \Phi$ and $\Phi \subset_{c,M} \Psi$;
\end{enumerate}
\end{defi}

The graphings $\Phi$ and $\Psi$ are said \defin{coarsely equivalent} (denoted $\Phi\ceq\Psi$) if they are $M$-coarsely equivalent for some $M \in \Nmath^\Nmath$. Clearly two graphings $\Phi$ and $\Psi$ are coarsely equivalent if and only if their length functions $\ell_\Phi$ and $\ell_\Psi$ on $\cg$ have the same bounded sequences.

When both $\Phi$ and $\Psi$ have finite size, they are coarsely equivalent if and only they are Lipschitz equivalent. So coarse equivalence is the natural generalization of Lipschitz equivalence for graphings with infinite size.

If a graphing $\Phi$ of the \pmp groupoid $\cg$ coarsely contains a generating graphing, then $\Phi$ itself is generating.

The $M$-\defin{coarse cost} of the graphed groupoid $(\cg,\Phi)$ is defined as the infimum of the costs of all graphings of $\cg$ which are $M$-coarsely equivalent to $\Phi$,

$$\cost_M(\cg,\Phi):=\inf\left\{\cost(\Psi)\colon \Psi\ceq_M\Phi \right\}.$$

We remark that $\Phi$ coarsely embeds into $\Psi$ if and only if for every $k$, there is $k'$ such that $\{\ph_1,\dots,\ph_k\}$ Lipschitz embeds into $\{\psi_1,\dots,\psi_{k'}\}$. This allows us to generalize in a straightforward way all the statements proved for graphings with bounded size. For example, Lemma~\ref{lem:goodgraph} still holds in this setting for graphings of finite cost.
\begin{lemm}\label{lem:goodgraphcoarse}
Let $(\cg,\Phi)$ be a graphed groupoid where the graphing $\Phi$ has finite cost.
The cost of $\cg$ can be computed inside the coarse class of $\Phi$, i.e.,  
\begin{equation}
\cost(\cg)=\inf_M\cost_M(\cg,\Phi). 
\end{equation}
\end{lemm}
\begin{proof}
  Let $\Phi=(\phi_i)_i$. Let $\eps>0$. Choose $k$ such that $\cost((\phi_i)_{i>k}) \leq \eps$. Consider the finite graphing $\Phi'=(\phi_i)_{i \leq k}$, and define $\cg' \subset  \cg$ the groupoid generated by $\Phi'$. Clearly $\cost(\cg) \leq \cost(\cg')+\eps$. By Lemma~\ref{lem:goodgraph} there is an integer $L_0$ and a finite size graphing $\Psi'$ which is $L_0$-Lipschitz equivalent to $\Phi'$ and such that $\cost(\Psi) \leq \cost(\cg')+\varepsilon$. Then the graphing $\Psi$ obtained as the union of $\Psi'$ and $(\phi_i)_{i>k}$ is coarsely equivalent to $\Phi$ and satisfies $\cost(\Psi') \leq \cost(\cg)+2\varepsilon$.
\end{proof}
\begin{rema}\label{rem:goodgraphcoarse} The proof shows more: for every $\varepsilon>0$ there is $i_\varepsilon$ and a graphing $\Psi$ of finite size which is coarsely equivalent to $\{\ph_1,\dots,\ph_{i_\varepsilon}\}$ and such that the graphing $\Psi \cup (\ph_j)_{j > i_\varepsilon}$ has cost $\leq \cost(\cg)+\varepsilon$.
\end{rema}

We generalize Definition~\ref{def: ul-comb-cost} of combinatorial cost to possibly unbounded graphings.
\begin{defi}
\label{def: comb cost, infinite size}
Let $(\cg_n,\Phi_n)_n$ be a sequence of graphed \pmp groupoids. The $\ul$-\defin{combinatorial cost} of the sequence is
\begin{equation}\ccost_\ul((\cg_n,\Phi_n)_n):=\inf_{M \in \Nmath^\Nmath} \lim_{n\in\ul}\cost_M(\cg_n,\Phi_n)\end{equation}
\end{defi}
When $\Phi_n$ have size bounded by $D$, we recover Definition~\ref{def: ul-comb-cost}. This deserves a small justification. Consider first some integer $L$. It follows from Remark~\ref{rem:for_Lipschitz_cost} that $\cost_{L}(\cg,\Phi_n)$ can be computed for a graphing $\Psi_n$ of size bounded by $(2D)^{L}$. This graphing is in particular $M$-coarsely equivalent to $\Phi_n$ for $M$ the constant sequence $((2D)^L)_{k \in \Nmath}$. In particular $\cost_M(\cg,\Phi_n) \leq \cost_{L}(\cg,\Phi_n)$ and \[\inf_{M \in \Nmath^\Nmath} \lim_{n\in\ul}\cost_M(\cg_n,\Phi_n) \leq \lim_{n \in \ul} \cost_{L}(\cg,\Phi_n).\] Conversely, if $M \in \Nmath^\Nmath$ and $\Psi_n$ is $M$-coarsely equivalent to $\Phi_n$, put  $L_0 := \max_{i \leq D} M_i$ and $L_1:=\max_{i \leq L_0} M_i$. Then we have $\Phi_n \subset_{L_0} \{\psi_1,\dots,\psi_{L_0}\}$ and $\{\psi_1,\dots,\psi_{L_0}\} \subset_{L_1} \Phi_n$. In particular $\cost_{\max(L_0,L_1)}(\cg_n,\Phi_n) \leq \cost_M(\cg_n,\Phi_n)$ and
\[\inf_{L \in \Nmath} \lim_{n\in\ul}\cost_L(\cg_n,\Phi_n) \leq \lim_{n \in \ul} \cost_{M}(\cg,\Phi_n).\]

\subsection{Limit of graphed groupoids} 

\begin{prop}\label{prop:compar}
 Consider a sequence of graphed \pmp groupoids $(\cg_n,\Phi_n)$ of bounded size. Take $L\in\Nmath$.
 \begin{enumerate}[label=(\roman*)]
 \item \label{item: limit of L-Lip embedd gr are L-Lip embedd} If for $\ul$-almost every $n$ the bisection $\psi_n$ of $\cg_n$ is $L$-Lipschitz embedded in $\Phi_n$, then $[\psi_n]_\ul$ $L$-Lipschitz embeds in $[\Phi_n]_\ul$.

 \item If the bisection $\psi_\ul$ of $\cg_\ul$ is $L$-Lipschitz embedded in $\Phi_\ul$, then there exists a sequence of bisections $\psi_n$ of $\cg_n$ such that $[\psi_n]_\ul=\psi_{\ul}$ and such that $\psi_n$ is $L$-Lipschitz embedded in $\Phi_n$ for $\ul$-almost every $n\in \Nmath$.
 
 \item Let $\Psi_\ul$ be a finite size graphing of the ultraproduct $\cg_\ul$ which $L$-Lipschitz contains $\Phi_\ul$. Then there exists a sequence of bounded size graphings $\Psi_n$ of $\cg_n$ such that $[\Psi_n]_\ul=\Psi_{\ul}$ and such that $\Psi_n$ $L$-Lipschitz contains $\Phi_n$ for $\ul$-almost every $n\in \Nmath$. Moreover if $\Psi_\ul$ is $L$-Lipschitz equivalent to $\Phi_\ul$, then we can choose $\Psi_n$ to be $L$-Lipschitz equivalent to $\Phi_n$.
  \end{enumerate}
\end{prop}
\begin{proof}
  Let $D:=\max_n |\Phi_n|$. We let $W$ be the finite set of formal words of length at most $L$ in the letters $\ph^i$ for $i\in \{1,\ldots,D\}$ and their inverses. For $w\in W$, denote by $w_n$ and $w_{\ul}$ the corresponding \pmp bisection of $\cg_n$ and $\cg_\ul$ respectively. We set \[A_n^w:=\{x\in \dom(\psi_n):\ \psi_n(x)=w_n(x)\}.\]

Let us first show \textit{(i)}. By hypothesis $\cup_{w\in W}A_n^w=\dom(\psi_n)$ for $\ulae$. Since $W$ is finite we have that $\dom([\psi_n]_\ul)=[\dom(\psi_n)]_\ul=\cup_w [A_n^w]_\ul$ and moreover $[\psi_n]_\ul\bigr|_{[A_n^w]_\ul}=w_\ul\bigr|_{[A_n^w]_\ul}$. This is exactly that $\psi_\ul$ $L$-Lipschitz embeds in $\Phi_\ul$.
\vspace{0.3cm}

The proof of \textit{(ii)} is similar. For every $w\in W$ we can find $B_n^w\subset X_n$ such that 
\begin{itemize}
\item $B_n^w\cap B_n^{w'}=\emptyset$ for $w,w'\in W$ with $w\neq w'$;
\item $[B_n^w]_\ul= A^w_\ul$  up to a nullset for every $w\in W$;
\item $\cup_{w\in W} [B_n^w]_\ul=\dom(\psi_\ul)$ up to a nullset. 
\end{itemize}
For every couple of words $w,w'\in W$ with $w\neq w'$ we observe that $\tar\left(w_\ul\bigr|_{[B_n^w]_\ul}\right)$ is disjoint from $\tar\left(w'_\ul\bigr|_{[B_n^{w'}]}\right)$ since they are restriction of the same bisection on two disjoint sets. Therefore, since $W$ is finite, we can assume that the measurable sets $B_n^w$ also satisfy the following property:
\begin{itemize}
\item $\tar\left(w_n\bigr|_{B_n^w}\right)\cap\tar\left(w'_n\bigr|_{B_n^{w'}}\right) =\emptyset$ for $w,w'\in W$ with $w\neq w'$.
\end{itemize}

Set $\psi_n$ to be the join of the partial authomorphism $w_n\bigr|_{B_n^w}$. Clearly $\psi_n$ $L$-Lipschitz embeds in $\Phi_n$ and $[\psi_n]=\psi_\ul$ outside a nullset.

For \textit{(iii)}, let us only prove the moreover part. The first is similar and left to the reader. Denote $\Phi_n = (\ph_n^k)_{1\leq k \leq D}$ and $\Psi_\ul = (\psi_\ul^k)_{1 \leq k \leq D'}$. Suppose we have $\Psi_\ul\Lipeq_L [\Phi_n]_\ul$ for some $L\in\Nmath$. \textit{(ii)} gives us a sequence $(\Psi_n)_n$ of size $D'$ such that $\Psi_n$ $L$-Lipschitz embeds into $\Phi_n$. Let $V$ be the finite set of formal words of length at most $L$ in the letters $\psi^i$ for $i\in \{1,\ldots,D'\}$ and their inverses and for $v\in V$ we denote by $v_n$ and $v_{\ul}$ the corresponding \pmp bisection of $\cg_n$ and $\cg_\ul$ respectively. For every $k \leq D$ we set
\[U^k_n:=\{x\in\dom(\ph^k_n):\ \ph^k_n(x)=v_n(x)\text{ for some }v\in V\}.\]

Observe that $\lim_\ul\mu_n(\dom(\ph^k_n)\setminus U^k_n)=0$. We define $\Psi'_n$ to be the union of $\Psi_n$ and the restriction of $\ph^i_n$ to $\dom(\ph^i_n)\setminus U_n^{i}$ for $i\leq k$. Then clearly 
$\Phi_n$ $L$-Lipschitz embeds in  $\Psi'_n$. Since we have added only restrictions of elements in $\Phi_n$, it follows that $\Psi'_n$  still $L$-Lipschitz embeds in $\Phi_n$. Finally $\Psi'_n$ and $\Psi_n$ differ only for finitely many elements whose support has arbitrarily small measure, therefore $[\Psi_n]_\ul=[\Psi'_n]_\ul=\Psi_\ul$. 
\end{proof}

Let us observe that if $\Phi_n$ have uniformly bounded size, then $\lim_\ul \cost(\Phi_n) = \cost(\Phi_\ul)<\infty$ and in general this condition amounts to a uniform summability of the measures of the domains of the bisections of $\Phi_n$. As a consequence of the above proposition we obtain the following theorem. 

\begin{theo}\label{thm:costccost} 
If a  sequence of \pmp graphed  groupoids $(\cg_n,\Phi_n)$ and its ultraproduct groupoid $(\cg_\ul,\Phi_\ul)$ satisfy $\lim_\ul \cost(\Phi_n) = \cost(\Phi_\ul) <\infty$, then \[\ccost_\ul((\cg_n,\Phi_n)_n)=\cost(\cg_\ul)\geq\lim_{\ul}\cost(\cg_{n}).\]
\end{theo}

\begin{rema}
\label{rem: counterex revers ineq non-group}
The reverse inequality
\begin{equation}
\cost(\cg_\ul)\leq \lim_{\ul}\cost(\cg_{n})
\end{equation}
is generally not true. 
Consider, for instance, the free group on $2$ generators (more generally a finitely generated group $\Gamma$ of cost $\vraicost_{*}(\Gamma)>1$), with finite generating set $S$ and with a decreasing sequence $(\Gamma_n)_n$ of finite index normal subgroups with trivial intersection.
The sequence of finite Schreier graphs for $(\Gamma/\Gamma_n, S)$ defines a sequence $\Phi_n=(\ph_s)_{s\in S}$ of graphings of the transitive equivalence relation $\cg_n$ on the finite spaces $X_n=\Gamma/\Gamma_n$ thus of cost $1-\frac{1}{[\Gamma:\Gamma_n]}$. The ultraproduct groupoid $\cg$ is given by a free action of $\Gamma$ on $X_\ul$ which has cost $>1$.
This example is exactly what justifies the use of the combinatorial cost \cite{Elek-combinatorial-cost-2007} or of the groupoid cost \cite{Abert-Nikolov-12} (see Cor.~\ref{cor: Abert-Nikolov}, where instead of a naked Schreier graph as above, the isotropy subgroups are taken into account).
\end{rema}

\begin{rema}\label{rem: counterex revers ineq non bounded}
The following inequalities hold in general without any boundedness assumption:
\begin{eqnarray}
\label{eq: C Gu geq limu C Gn}
\ccost_\ul((\cg_n,\Phi_n)_n) &\geq& \lim_\ul\cost(\cg_n)\\
\label{eq: C Gu leq cCu Gn general}
\ccost_\ul((\cg_n,\Phi_n)_n) &\geq& \cost(\cg_\ul).
\end{eqnarray}
(\ref{eq: C Gu geq limu C Gn}) is clear from the definition of combinatorial cost, while
(\ref{eq: C Gu leq cCu Gn general}) holds since the first part of the argument below is still valid in this context.

On the contrary, the other inequalities
\begin{eqnarray}
\cost(\cg_\ul) &\geq& \lim\limits_\ul \cost(\cg_n)
\label{eq:cGu geq lim C(Gn)}
\\
\cost(\cg_\ul) &\geq& \ccost_\ul((\cg_n,\Phi_n)_n)
\label{eq:cGu geq cC(Gn)}
\end{eqnarray}
both require some boundedness assumption as examplified by a residually finite fixed price $1$ group such as $\Gamma=\FF_{\infty}\times \Zmath$ with a decreasing sequence of finite index normal subgroups with trivial intersection.
The groupoid cost for the action on any finite quotient is infinite, and similarly for the combinatorial cost. But the ultraproduct free action has cost $1$ \cite[Proposition VI.23]{Gab-cost}.

Even with $\sup_n \cost(\Phi_n) < \infty$, uniform summability is needed:
If instead of $\Phi_n=(\ph_s)_{s\in S}$ in the example of Remark~\ref{rem: counterex revers ineq non-group} one subdivides the bisections $\ph_s$ into singletons, we obtain the generating graphing 
$\Phi'_n=(\ph'_{v,s}=(v,sv))_{v\in \Gamma/\Gamma_n, s\in S}$, $\mu(\dom(\ph'_{v,s})) = \frac {1}{ [\Gamma:\Gamma_n]}$.
Then $\cg_\ul$ is the trivial groupoid $X_\ul$ and in particular has cost $0$, whereas $\ccost_\ul (\cg_n,\Phi_n) \geq \lim_\ul \cost(\cg_n) = 1$, thus contradicting (\ref{eq:cGu geq lim C(Gn)}) and (\ref{eq:cGu geq cC(Gn)}).

\end{rema}

\begin{proof}[Proof of Theorem~\ref{thm:costccost}] Let $M \in \Nmath^\Nmath$ and $\Psi_n$ be a graphing of $\cg_n$ such that $\Phi_n \ceq_M \Psi_n$. By Corollary~\ref{cor:dependance_en_graphage}, $\Psi_\ul$ is a generating graphing of $\cg_\ul$ and in particular
  \[ \cost(\cg_\ul) \leq \cost(\Psi_\ul) \leq \lim_\ul \cost(\Psi_n).\]
  Taking the infimum over $\Psi_n$ we obtain
  \[ \cost(\cg_\ul) \leq \lim_\ul \cost_M(\cg_n,\Phi_n).\]
This shows (without the uniform integrability assumption) that
\[\cost(\cg_\ul) \leq \ccost_\ul((\cg_n,\Phi_n)_n).\]

For the converse inequality, fix $\eps>0$. By Remark~\ref{rem:goodgraphcoarse} there are integers $i_0$ and $L$ and a finite graphing of $\cg_\ul$, $\Psi'_\ul \Lipeq_L (\ph_\ul^1,\dots,\ph_\ul^{i_0})$ such that the graphing $\Psi_\ul :=\Psi'_\ul \cup (\ph_\ul^i)_{i > i_0}$ has cost $\leq \cost(\cg_\ul)+\varepsilon$. By Proposition~\ref{prop:compar} for every $n$ there is a graphing $\Psi'_n\Lipeq_L (\ph_n^1,\dots,\ph_n^{i_0})$ such that $[\Psi'_n]_\ul =\Psi'_\ul$. In particular, there is $M \in \Nmath^\Nmath$ such that $\Psi_n := \Psi'_n \cup (\ph_n^i)_{i >i_0}$ is $M$-coarsely equivalent to $\Phi_n$. Finally, 
\[\lim_\ul \cost(\Psi_n) = \lim_\ul \cost(\Psi'_n) + \lim_\ul \cost((\ph_n^i)_{i>i_0}).\]
The first term is equal to $\cost(\Psi'_\ul)$ by linearity, and the second term  is equal to $\cost((\ph_\ul^i)_{i>i_0})$ by the uniform integrability assumption. In particular we obtain
\[ \lim_\ul \cost(\Psi_n) = \cost(\Psi_\ul) \leq \cost(\cg_\ul)+\eps.\]
Taking the infimum over all $\eps$ we obtain $\ccost_\ul((\cg_n,\Phi_n)_n) \leq \cost(\cg_\ul)$, which proves the theorem.\end{proof}

\begin{rema}[Elek's combinatorial cost]
\label{Rem: Elek comb cost}
Initially, the combinatorial cost \cite{Elek-combinatorial-cost-2007} was introduced as an invariant for sequences of finite graphs with uniformly bounded degree $\leq D$. 
If $(G_n)_n$ is such a sequence of graphs, then we have already discussed a canonical coarse structure on the sequence, see Example~\ref{ex:coarse_structure_on_graphs}, and commented about the ultraproduct of the sequence in Example~\ref{ex:ultraproduct_of_graphs}. 

 The original combinatorial cost $\ccost$ of the sequence $(G_n)_n$ uses $\liminf\limits_{n\to \infty}$ instead of $\lim\limits_{n\in \ul}$ (compare Definition~\ref{def: ul-comb-cost}):
 \begin{equation}\label{def: init Elek comb cost with liminf}
 \ccost((G_n)_n):=\inf_L\liminf\limits_{n\to \infty}\cost_L(\cg_n,\Phi_n).
 \end{equation}
 We claim that:
 \begin{equation}\label{eq:inequ combcost vs u-combcost}
\ccost((G_n)_n) = \inf \{\ccost_{\ul}((\cg_n,\Phi_n)_n) : \ul \textrm{ non principal ultrafilter}\}.
\end{equation}
The inequality $\leq $ holds since, $\liminf\limits_{n\to \infty}\cost_L(\cg_n,\Phi_n)\leq \lim\limits_{n\in \ul}\cost_L(\cg_n,\Phi_n)$
for every $L$.
For the reverse inequality, for every $\varepsilon>0$, there is an $L$ such that $\liminf\limits_{n\to \infty}\cost_L(\cg_n,\Phi_n)\leq \ccost((G_n)_n) +\varepsilon$ and there is a non principal ultrafilter $\ul$ such that  $\liminf\limits_{n\to \infty}\cost_L(\cg_n,\Phi_n)=\lim\limits_{n\in \ul}\cost_L(\cg_n,\Phi_n)$.
Observe that, interleaving the terms of two sequences of graphs with different combinatorial costs, indicates that a strict inequality
 $\ccost((G_n)_n) <\ccost_{\ul}((\cg_n,\Phi_n)_n)$ may happen for some ultrafilters.
 Moreover, it is not clear to us whether the infimum in (\ref{eq:inequ combcost vs u-combcost}) is realized by some ultrafilter.
\end{rema}

\begin{coro}[Combinatorial cost vs cost]\label{cor: combcost vs cost}
If $(G_n)_n$ is  a sequence of finite graphs with uniformly bounded degree, then the combinatorial cost is the infimum over all non principal ultrafilters of the costs of the ultraproduct groupoids $\cg_\ul$ of the associated coarse groupoids $\cg_n$:
 \begin{equation*}
 \ccost((G_n)_n) = \inf \{\cost(\cg_{\ul}) : \ul \textrm{ non principal ultrafilter}\}.
\end{equation*}
\end{coro}

\subsection{Weak containment}
We now mimic the powerful notion of weak containment introduced by Kechris \cite[Section 10 (C)]{Kechris-Global-aspects-2010}.

\begin{defi}[Weak containment for graphed groupoids]
The  graphed groupoid $(\ch,\Psi)$ on $(Y,\ra, \nu)$ is \defin{weakly contained in} the graphed groupoid $(\cg, \Phi)$ on $(X, \rb, \mu)$, in symbol,
\[(\ch,\Psi) \prec (\cg, \Phi)\] if $\Psi=(\psi_j)_{j\in J}$ and
$\Phi=(\ph_j)_{j\in J}$ have the same sets of indices and for any finite set $K \subset [[\cH]]$ and every $\eps>0$, there is a map $\pi \colon K \to [[\cg]]$ satisfying $\pi(\psi_j) = \ph_j$ for every $\psi_j \in K \cap \Psi$, $d(\pi(\psi\psi'),\pi(\psi)\pi(\psi')) \leq \eps$ and $|\tau(\psi) - \tau(\pi(\psi))|<\eps$ for every $\psi,\psi' \in K$ such that $\psi\psi' \in K$.
\end{defi}
Note that if $B \in K$ is a \defin{unit bisection} of $\ch$, (meaning $B \subset Y$ or equivalently $B^2=B$), then the condition $d(\pi(B),\pi(B)^2) \leq \eps$ implies that $\pi(B)$ is at distance $\leq \eps$ from a unit bisection. So in the definition of weak containment, we can always assume that $\pi$ maps unit bisections to unit bisections, and also that disjoint unit bisections are mapped to disjoint unit bisections. Moreover, if $\ph$ is any bisection, then the \defin{conjugation} satisfies $\ph B \ph^{-1}=\overline{\ph}(B)$. One thus recovers Kechris' notion of weak containment  for group actions. 

If $\ra$ is separable\footnote{When $\ra$ is not separable, the same holds if one allows an ultrafilter on a larger index set than $\Nmath$, as in Proposition~\ref{prop:relation_weakequivalence_factor}.}, then $(\ch,\Psi) \prec (\cg, \Phi)$ if and only if there is a trace preserving embedding of $[[\ch]]$ into the metric ultraproduct of $\prod_\ul [[\cg]]$ sending each $\ph_j$ to $[\psi_j]_\ul$, see Proposition~\ref{prop:relation_weakequivalence_factor}.
\begin{prop}\label{prop:properties weak-contain graphings}
Assume $(\ch,\Psi) \prec (\cg, \Phi)$ and let $\psi$ be a $\Psi$-word and let $\ph$ be the corresponding $\Phi$-word, then 
\begin{enumerate}
\item
$\widetilde\nu(\psi)=\widetilde\mu(\varphi)$; thus $\cost(\Psi)=\cost(\Phi)$;

\item $\tau(\psi)=\tau(\ph)$ (equal measure of units);

\item
$\widetilde \nu(M_F(\psi))\leq \widetilde \mu(M_F(\ph))$ (the measurably free part from Proposition~\ref{prop:manyinvo-gen} is monotonic).
\end{enumerate}
\end{prop}
\begin{proof}
The first two items are easy and left to the reader.

For the measurably free part, use finitely many elements from a partition given by Proposition~\ref{prop:manyinvo-gen}.
\end{proof}

\begin{defi}[Weak equivalence for graphed groupoids]
\label{def: wek-equiv graphed groupoid}
The  graphed grou\-poids $(\ch,\Psi)$ and $(\cg, \Phi)$ are weakly equivalent, in symbols, 
\[(\ch,\Psi) \sim_{w} (\cg, \Phi)\]
if $(\ch,\Psi) \prec (\cg, \Phi)$ and $(\cg, \Phi)\prec (\ch,\Psi) $.
\end{defi}

\begin{rema}[Local-global convergence]
\label{rem: local-global cv}
  Hatami, Lov\'asz and Szegedy \cite{HatLovSze} have defined a notion of local-global convergence for sequences of bounded-degree graphs. This enters rather well in our setting, the only difference being that they consider unlabeled and unoriented graphs whereas we have mainly focused on labeled graph(ing)s. In order to interpret the local-global convergence in our setting, let us say that a subset $\Phi \subset \cg$ of a \pmp groupoid is a \defin{bounded degree combinatorial graphing} if $\Phi$ is measurable, symmetric ($\Phi^{-1} = \Phi$) and $x \mapsto |s^{-1}(x) \cap \Phi|$ is essentially bounded on $X$. Note that a symmetric subset of $\cg$ is a combinatorial graphing if and only if it is a finite union of bisections.

If $\cg$ is the equivalence relation on a finite set, we recover the classical notion of combinatorial graph. If $\cg$ is a \pmp equivalence relation on a standard probability space, then we recover the notion of graphing as in \cite{Kechris-Miller}. In that paper Hatami, Lov\'asz and Szegedy first define a notion of convergence, and then they identify the limiting objects as graphings in their sense. To be more precise, we wish to point out that the limiting object is not unique: the well-defined limiting object is a \emph{weak-equivalence class of combinatorial graphings}, as we define now.
  
Let $\cg$ and $\ch$ be \pmp graphings with combinatorial graphings $\Phi$ and $\Psi$. Write $\Phi = \cup_{j=1}^n \ph_j$ as a finite union of bisections. We say that $(\ch,\Psi)$ is weakly contained in $(\cg,\Phi)$ if for any finite set $K \subset [[\cH]]$ and every $\eps>0$, there is a map $\pi \colon K \to [[\cg]]$ satisfying $\cup_{j=1}^n \pi(\psi_j) = \Phi$ if $\psi_j \in K$ for every $j$, $d(\pi(\psi\psi'),\pi(\psi)\pi(\psi')) \leq \eps$ and $|\tau(\psi) - \tau(\pi(\psi))|<\eps$ for every $\psi,\psi' \in K$ such that $\psi\psi' \in K$.

As above if $\ra$ is separable, then $(\ch,\Phi)$ is weakly contained in $(\cg,\Psi)$ if and only if there is a trace preserving embedding of $[[\ch]]$ into the metric ultraproduct of $\prod_\ul [[\cg]]$ sending each $\Phi$ to $\Psi$. The same holds \emph{mutatis mutandis} for non separable $\ra$. In particular this notion does not depend on the decomposition $\Phi = \ph_1\cup\dots\cup \ph_n$.

One says that $(\cg,\Phi)$ and $(\ch,\Psi)$ are weakly equivalent if $(\cg,\Phi)$ is weakly contained in $(\ch,\Psi)$ and conversely  $(\ch,\Psi)$ is weakly contained in $(\cg,\Phi)$. 

Then one can easily check that a sequence $(G_n)$ of bounded-degree graphs converges locally-globally in the sense of \cite{HatLovSze} if and only if the weak equivalence class of the graphed equivalence relation $(\cg_\ul,\mathrm{Edg}_\ul)$ (see Example~\ref{ex:ultraproduct_of_graphs}) does not depend on the free ultrafilter $\ul$. More generally, the sequence $(G_n)$ local-global converges along $\ul$ towards a \pmp equivalence relation $\RR$ with a combinatorial graphing $\Psi$ if and only if $(\RR,\Psi)$ is weakly equivalent to $(\cg_\ul,\mathrm{Edg}_\ul)$.

Let us mention that ultraproducts of bounded degree finite graphs as combinatorial graphings were already considered in \cite{Elek-2010-param-test} (we thank G. Elek for pointing this article out after a preliminary version of our text has circulated). 
\end{rema}

The following extends a result of Kechris \cite[end of Section 10 (C)]{Kechris-Global-aspects-2010} from group actions to graphed groupoids:
\begin{theo}[weak containment and cost]
If $(\ch,\Psi) \prec (\cg, \Phi)$ and $\Phi$ has finite cost,  then $\cost (\ch) \geq  \cost(\cg)$.
In particular, the cost is an invariant of the weak equivalence class of $(\cg, \Phi)$.
\end{theo}
\begin{proof}

  First observe that the graphings $\Phi=(\ph_j)_{j\in J}$ and $\Psi=(\psi_j)_{j\in J}$ have the same cost (see Proposition~\ref{prop:properties weak-contain graphings}).

  Let $\eps_0>0$. By (the proof of) Lemma~\ref{lem:goodgraphcoarse}, there exists a generating graphing $\Psi'=(\psi'_i)_{i\in I}$ of $\ch$
such that:
\\
(1) $\cost (\Psi')\leq \cost(\ch)+\eps_0$;
\\
(2) there are finite subsets $J_0\subset J$ and $I_0\subset I$ such that $\Psi'_{0}:=(\psi'_i)_{i\in I_0}$ generates the same groupoid as $\Psi_{0}:=(\psi_j)_{j\in J_0}$ (even Lipschitz equivalent, but this will be of no use here),
\\
(3) there is a bijection $r:J\setminus J_0\to I\setminus I_0$ such that $\psi_j=\psi'_{r(j)}$ for each $j\in J \setminus J_0$.

Let $K \subset[[\cH]]$ be a finite subset containing $\Psi'_0$ and $\eps>0$, and let $\pi \colon K \to [[\cg]]$ be as in the definition of weak containment. Denote by $\cg_{K,\eps}$ the groupoid generated by $\pi(\Psi'_0)=\{ \pi(\psi'_i)\}_{i \in I_0}\}$. By (2) and the definition of weak containment, $\lim_{K,\eps} d(\ph_j,\ph_j \cap \cg_{K,\eps})=0$ for every $j \in J_0$. In particular there is $K,\eps$ such that $\sum_{j \in J_0} d(\ph_j,\ph_j \cap \cg_{K,\eps}) \leq \eps_0$. Taking $\eps$ small enough we can also assume that $\cost(\pi(\Psi'_0)) \leq \cost(\Psi'_0)+\eps_0$. It follows that the graphing of $\cg$ obtained as the union of $\pi(\Psi'_0)$, of $(\ph_j \setminus (\ph_j \cap \cg_{K,\eps}))_{j \in J_0}$ and of $(\ph_j)_{j \in J\setminus J_0}$ is generating and has cost less than
\[ \cost(\Psi'_0)+\eps_0+\eps_0 + \cost((\ph_j)_{j \in J\setminus J_0}) = \cost(\Psi')+2\eps_0 \leq \cost(\ch)+3\eps_0.\]
One concludes that $\cost(\cg) \leq \cost(\ch)$ as $\eps_0$ was arbitrary.
\end{proof}

\subsection{Factors}

Recall that homomorphisms of \pmp groupoids are defined in Definition~\ref{def:groupoid_homomorphism}.
\begin{defi}\label{dfn:factor}
A homorphism $T$ between the \pmp groupoids $\cg$ and $\cg'$ over the probability spaces $(X,\mu)$ and $(X',\mu')$ is a \defin{factor map} if $T$ maps $\widetilde \mu$ to $\widetilde \mu'$. In this case $\cg'$ is a \defin{factor} of $\cg$.
\end{defi}

Be aware that the condition on the measures $\widetilde \mu$ and $\widetilde \mu'$ entails that the factor map is essentially onto and is ``class-bijective'', i.e., $T(g)$ is a unit if and only if $g$ itself is a unit $\widetilde \mu$-almost surely; in particular $\widetilde \mu(T^{-1}(X')\setminus X)=0$ and $T$ is \pmp in the sense of Definition~\ref{def:groupoid_homomorphism}. 

The second author takes the opportunity to mention that this notion was introduced in the framework of standard equivalence relations in
 \cite[Definition p. 1815 and Lemma 2.3 and 2.4]{Gab05a}
 and to correct an obvious misprint there: replace one-to-one by bijective.

For instance, a group morphism $\Gamma\to \Gamma'$ is a factor in our sense if and only if it is an isomorphism. In fact, we can completely characterize the groupoids which factors over a given group.
\begin{rema}[\pmp actions of \pmp groupoids]
\label{Rem: pmp actions of pmp groupoids}
 If $\cg'$ is a countable group $\Gamma$ (seen as a groupoid over a point), then every \pmp action of $\Gamma$ on $X$ gives rise to a factor $\cg_{\Gamma'\acting X} \to \Gamma$ sending $(\gamma,x)$ to $\gamma$, and conversely every factor $\cg \to \Gamma$ arises this way. In general one can interpret factors $\cg \to \cg'$ as \defin{\pmp actions of the \pmp groupoid $\cg'$ on $X$} (not to be confused with the different notion of actions -- not \pmp -- of a groupoid considered in Section~\ref{sec:actions}). Free actions correspond to factors $\RR \to \cg'$ where $\RR$ is a principal \pmp groupoid (i.e. equivalence relation).
\end{rema}

 This leads to the fixed price problem for general groupoids.
 \begin{ques}[Fixed price problem for \pmp goupoids]
 \label{quest: fixed price for groupoids}
  When  $\cg$ is a \pmp groupoid its infimum cost $\vraicost_*(\cg)$ is the infimum of the costs of all 
\pmp groupoids $\ch$ that factor onto it:
\[\vraicost_*(\cg):=\inf\{\cost (\ch): \ch \text{ \pmp groupoid that  factors onto } \cg  \}.\]
  The \pmp groupoid 
$\cg$ has fixed price if $\cost (\ch)=\vraicost_*(\cg)$ for every principal (i.e. with trivial isotropy) \pmp groupoid $\ch$ which factors onto $\cg$.
Does there exist a \pmp groupoid that does not have fixed price?
(Compare \cite[Question 1.8]{Gab-cost}).
\end{ques}
Observe that a similar question has been raised by
Popa--Shlyakhtenko--Vaes \cite{PSV-2018-classif-sub-hyp}
in the context of standard equivalence relations.

If $T:\cg\rightarrow \cg'$ is a factor map, then for every bisection $\ph'\subset \cg'$, its preimage $\ph:=T^{-1}(\ph')$ is a bisection $\widetilde{\mu}$-mod $0$. Indeed let us check that $s$ is essentially injective on $\ph$ (the same argument works for $t$). Assume that $g,h\in \ph$ have the same source. Then $s(T(g))=s(T(h))$ and hence $T(g)=T(h)$. Therefore $T^{-1}(T(g)T(h)^{-1})\subset X$ (up to some $\widetilde{\mu}$ null set) and in particular $gh^{-1}$ is a unit and therefore $g=h$. As a corollary we obtain that the preimage of a graphing $\Phi'$ is still a graphing, which we call the \defin{pull-back graphing} and denote by $T^{-1} \Phi'$. We also get the following.

\begin{prop}
  If $T:\cg\rightarrow \cg'$ is a factor map, then $\cost(\cg)\geq\cost(\cg')$. 
\end{prop}

 Observe that isotropy elements are always sent to isotropy elements: if $s(g)=t(g)$, then $s(T(g))=t(T(g))$.

A factor $T:\cg\rightarrow\cg'$ is called \defin{free} if it moreover satisfies one of the following equivalent conditions:
\begin{itemize}
\item for almost every $g\in\cg$ if $s(T(g))=t(T(g))$, then $s(g)=t(g)$;
\item for almost every $x\in X$ the factor $T$ realizes an isomorphim between the isotropy groups of $x$ and $T(x)$.
\end{itemize}

A factor $T:\cg\rightarrow\cg'$ is called \defin{measurably free} if 
for every bisection $\ph'\subset \cg'$, the measurably free parts (see Proposition~\ref{prop:manyinvo-gen}) of $\ph'$ and $\ph=T^{-1}(\ph')$ coincide: 
$M_F(\ph)=T^{-1}(M_F(\ph'))$.

Clearly every measurably free factor of a realizable \pmp groupoid is automatically free, see Definition~\ref{dfn:realizable}. 

  Remark that if $\RR$ and $\RR'$ are \pmp equivalence relation over the standard probability spaces $(X,\mu)$ and $(X',\mu')$ and $T$ is a morphism between $\RR$ and $\RR'$ such that 
\begin{inparaenum}[(a)]
 \item $T$ restricted to $X$ is surjective and measure preserving onto $X'$,
  \item $T$ is bijective on every equivalence class,
 \end{inparaenum}
 then $T$ is automatically a factor map.  
 
\begin{prop}\label{prop: projective limit and costs}
Let $(\cg_n,\Phi_n)$ be a sequence of graphed groupoids and assume that  for every $n$, $T_n: (\cg_{n+1},\Phi_{n+1})\to (\cg_n,\Phi_n)$  is a factor map such that $\Phi_{n+1} = T_n^{-1} \Phi_n$. Assume furthermore that $|\Phi_1|=|\Phi_n|$ is finite. Then $(\cost(\cg_n))_n$ is non increasing with $n$ and for every non principal  $\ul$
 \[\ccost_\ul((\cg_n,\Phi_n)_n)=\cost(\cg_\ul)=\lim_{n\to \infty}\cost(\cg_{n}).\]
\end{prop}
\begin{proof}
All the $\Phi_n$ have the same size $|\Phi_1|$.
Let $\Psi_n$ be $L_n$-Lipschitz equivalent with $\Phi_n$ and satisfy $\cost(\Psi_n)\leq \cost(\cg_n)+2^{-n}$. 
The $L_n$-Lipschitz expression of each bisection in $\Psi_n$ can be pulled-back in $\Phi_p$ when $p\geq n$.
Then the canonical pull-back $\Psi_{n,p}$ of $\Psi_n$  in $\cg_p$ is $L_n$-Lipschitz equivalent with $\Phi_p$ and generates $\cg_p$. It follows that $\ccost_\ul ((\cg_n, \Phi_n)_n)\leq \lim\limits_\ul \cost(\Psi_n)=\lim\limits_{\ul} \cost (\cg_n)$. Then  Theorem~\ref{thm:costccost} gives the reverse inequality.
By monotonicity 
$\lim_{\ul} \cost (\cg_n)=\lim_{n\to \infty} \cost (\cg_n)$
and thus 
$\ccost_\ul ((\cg_n, \Phi_n)_n)=\cost(\cg_\ul)=\lim\limits_{n\to \infty}\cost(\cg_{n})$ for every non principal $\ul$.
\end{proof}

A different way of proving the above proposition is to show that the projective limit of the groupoids $\cg_n$ is weakly equivalent to the ultraproduct $\cg_\ul$. 

\begin{theo}[Standard factor]
\label{thm: standard factors}
 Every graphed \pmp groupoid $(\cg, \Phi)$ over any probability space $(X,{\mathscr B},\mu)$ admits some \pmp factor $(\cg', \Phi')$ on a standard Borel probability space such that:
    \begin{enumerate}
      \item $(\cg', \Phi')$ is weakly-equivalent to $(\cg, \Phi)$;
      \item the pull-back graphing of $\Phi'$ is $\Phi$;
      \item $\cost(\cg)=\cost(\cg')$ (even when $\cost(\Phi)=\infty$); \label{it: standard factors, cost}
      \item the factor $\cg\to \cg'$ is measurably free.
 \end{enumerate}

Moreover if $\cg$ is realizable, then the factor is free.
\end{theo}

The existence of a standard factor is an incarnation of the \textit{downward Löwenheim–Skolem theorem} in model theory which roughly states that every model of a countable theory has a separable substructure. In the case of countable groups the existence of the standard factor was also explicitly stated and proved in \cite{Card-ultra-prod-we-sofic-arxiv}.

\begin{exam}[Example~\ref{ex: example equiv-rel pmp act non meas.} continued]
For both actions $\alpha$ and $\beta$ of $\Gamma=\Zmath/2\Zmath$ on $[-1,1]$ from Example~\ref{ex: example equiv-rel pmp act non meas.}, the groupoids $\cg_{\alpha}$ and $\cg_{\beta}$ have two natural standard factors:
$[0,1]$ with constant isotropy  $\equiv \Zmath/2\Zmath$ and the group $ \Zmath/2\Zmath$ itself. Both have cost $=1$.
Only the first one is weakly equivalent to $\cg_{\alpha}$ (resp. $\cg_{\beta}$). 

On the other hand, in any standard factor, every measurably isotropic element must be sent to a totally isotropic element.
This forces every standard factor from $\cg_{\alpha}$ and $\cg_{\beta}$ to transit through the quotient of $[-1,1]$ by "$x\mapsto -x$ on some measurable subset containing $\Omega$".
\end{exam}

\begin{proof}[Proof of Theorem~\ref{thm: standard factors}]
Consider an auxiliary graphing $\Psi$ which contains $\Phi$ as well as a sequence of graphings $\Phi_n$ whose costs tend to $\cost (\cg)$. Let 
$\mathcal{D}\subset \rb$ be the countable collection which is invariant under all the $\hat{\Psi}$-words and which contains:
\begin{itemize}
\item  for each $\Psi$-word $m$, its unit subset $N_m:=m\cap X$ and the elements $s(A_i)$ and $s(\Fix_{\rb}(m))$ for $(A_i)_i$ and $\Fix_{\rb}(m)$ given by Proposition~\ref{prop:manyinvo-gen};
\item  for each $\Phi$-word $w$, a countable collection $\mathcal{C}_{w}\subset \rb$ of Borel subsets such that the set $\{\mu(\hat{w}(C_j)\cap C_k)): C_j, C_k\in \mathcal{C}_{w}\}$ is dense in the set 
$\{\mu(\hat{w}(B)\cap B')): B, B'\in \rb\}$.
\end{itemize}
  Let $X':=\{0,1\}^{\mathcal{D}}$ with its canonical Borel $\sigma$-algebra $\rb'$
  and define a (measurable) map $\pi:X\to X'$ by $x\mapsto \{\chi_{D}(x)\}_{\mathcal{D}}$ (where $\chi_{D}$ is the characteristic function of $D\in \mathcal{D}$). Let $\mu'=\pi_*\mu$. Then $(X',\mu')$ is a standard Borel probability space. Observe that $\pi^{-1}(\rb')$ is the $\sigma$-algebra generated by $\mathcal{D}$.

Each $\hat{\psi}\in \hat{\Psi}$ delivers a partial isomorphism $\hat{\psi}':\pi(\dom(\hat{\psi}))\to \pi(\tar(\hat{\psi}))$ given for $x\in \dom(\hat{\psi})$ by $\hat{\psi}'(\pi(x))=\pi(\hat{\psi} (x))$.
Let's check that this is well defined:
If $x$ and $y$ are not separated by $\mathcal{D}$, then the same holds for $\hat{\psi}(x)$ and $\hat{\psi}(y)$ by invariance of $\mathcal{D}$.
If $U'\in \rb'$ is contained in $\tar(\hat{\psi}')$, then $\hat{\psi}'^{-1}(U')=\pi(\hat{\psi}^{-1}(\pi^{-1} (U')))$ is $\rb'$-measurable since $\pi^{-1} (U')$ and $\hat{\psi}^{-1}(\pi^{-1} (U'))$ are measurable for $\pi^{-1}(\rb')$.

 Let $\cf$ and $\cf'$ be the free groupoids generated by $\hat{\Psi}$ and $\hat{\Psi}':=\{\hat{\psi}': \hat{\psi}\in \hat{\Psi}\}$ respectively. Observe that $\pi$ defines a factor $\tilde{T}$ from $\cf$ to $\cf'$. Consider the natural map $\rho:\cf\to\cg$ mapping each bisection $\psi$ of $\cf$ to the corresponding one in $\cg$. Observe that the kernel of this map $\mathcal N\subset \cf$ is a  totally isotropic normal subgroupoid of $\cf$, defined by the measurable sets $N_m=m\cap X$ (for each $\hat{\Psi}$-word $m$). 
Its image $\mathcal N':=\tilde{T}(\mathcal N)<\cf'$ is the normal saturation in $\cf'$ of the relator system defined, for each $\Psi'$-word $m'$ associated to the $\hat{\Psi}$-word $m$, by $R'_{m'}:=\{x'\in X': \text{ whose } N_m\text{-th coordinate }=1\}$.
It is a totally isotropic normal subgroupoid.

"The" standard \pmp factor is now defined as the quotient $\cg':=\cf'/\mathcal N'$ (see Proposition~\ref{prop: quot groupoid}). 

Checking the rest is routine:
$\tilde{T}$ goes down to a map $T:\cg\to \cg'$. 
The first condition in the choice of $\mathcal{D}$ ensure that we obtain a measurably free factor.
Since $\Phi\subset \Psi$ is a graphing of $\cg$, it follows that $\Phi'=T(\Phi)$ is a graphing of $\cg'$.
The second condition ensures that $(\cg, \Phi)\prec (\cg',\Phi')$.
Even if $\Phi$ does not have finite size, the condition that $\Psi$ contains graphings $\Phi_n$ of cost arbitrarily close to that of $\cg$ allows, using the graphings $T(\Phi_n)$, to ensure that $\cost(\cg')\leq \cost(\cg)$. 

On the other hand, $\cg'$ being a factor of $\cg$, we thus have $(\cg, \Phi)\sim_{w} (\cg',\Phi')$ and the reverse inequality
$\cost(\cg')\geq \cost(\cg)$.
\end{proof}

The next statements illustrate the relationships between the notions of factors and weak containment. We start with two easy observations.

\begin{lemm}\label{lem:ultrapower_weakly_equiv} The ultrapower $(\cg_\ul,\Phi_\ul)$ of a constant sequence $(\cg, \Phi)$ of graphed groupoids is weakly equivalent to $(\cg,\Phi)$.
\end{lemm}

\begin{lemm} 
\label{Lem: factors vs full pseudogroups}
If $T\colon \cg \to \cH$ is a factor, then the map $T^{-1}\colon [[\cH]] \to [[\cg]]$ preserves the trace and all operations.

    Conversely, if $Y$ is a standard measure space, then any such embedding $[[\cH]] \to [[\cg]]$ comes from a factor.
\end{lemm}
\begin{proof} Only the second statement deserves a proof. We will use the following classical fact (see \cite[Theorem 343B]{fremlin2011measure}, or \cite{vonNeumann32} for the case when both spaces are standard): every isometric embedding of the measure algebra of a standard $\sigma$-finite space $\Omega_1$ to the measure algebra of an arbitrary measure space $\Omega_2$ sending $\Omega_1$ to $\Omega_2$ and preserving the boolean algebras operations comes from a measure-preserving map $\Omega_2 \to \Omega_1$. We recall a proof for completeness: one readily reduces to probability spaces and, by the classification of standard probability spaces, one can assume that $\Omega_1=\{0,1\}^\Nmath$ for a Borel probability measure. The required map $T \colon \Omega_2 \to \Omega_1$ is then defined by $T(x)(n)=1$ if $x$ belongs to the image of the cylinder $\{\omega \in \Omega_1: \omega(n)=1\}$ and $T(x)(n)=0$ otherwise.

An isometric embedding of the pseudogroups $[[\cH]] \to [[\cg]]$ preserves disjointness, so extends to an isometric embedding of the measure algebra of $\cg$ into the measure algebra of $\ch$, which as recalled above in turn comes from a measure-preserving map $\cg \to \ch$ as $\ch$ is standard. The fact that the original embedding $[[\cH]] \to [[\cg]]$ preserves the pseudogroup operations implies that the associated map $\cg \to \ch$ preserves the groupoid operations almost everywhere, so is a factor.
  \end{proof}

The following proposition generalizes \cite[Theorem 1]{CKTD} in our context.

\begin{prop}\label{prop:relation_weakequivalence_factor} 
Let  $(\cH,\Psi)$ and  $(\cg, \Phi)$ be two graphed groupoids over $(Y,\ra, \nu)$ and $(X, \rb, \mu)$ respectively. Then $(\cH,\Psi) \prec (\cg, \Phi)$ if and only if there is an ultrafilter $\ul$ on some index set $I$ and a trace preserving pseudogroup homomorphism from $[[\ch]]$ into the full pseudogroup $[[\cg_\ul]]$ of the ultrapower $(\cg_\ul, \Phi_\ul)$ of the constant family $(\cg, \Phi)_{i\in I}$, which sends $\Psi$ to $\Phi_\ul$.

Moreover, if (the measure algebra of) $\ra$ is separable, then $I$ can be taken to be equal to $\Nmath$ and $\ul$ can be any nontrivial ultrafilter.
\end{prop}
\begin{proof} The $\Leftarrow$ direction is a direct consequence of Lemma~\ref{lem:ultrapower_weakly_equiv}. Conversely, assume that $(\cH,\Psi)$ is weakly contained in $(\cg,\Phi)$. Denote by $I$ the set of all tuples $(N,\varepsilon,A_1,\dots,A_N)$ with $N$ an integer, $\varepsilon>0$ and $A_1,\dots,A_N \in \ra$. The order $(N,\varepsilon,A_1,\dots,A_N)\leq (N',\varepsilon',A'_1,\dots,A'_{N'})$ if $N \leq N'$, $\varepsilon>\varepsilon'$, $\{A_1,\dots,A_N\} \subset \{A'_1,\dots,A'_{N'}\}$ turns $I$ into a directed set. Choose $\ul$ an ultraproduct on $I$ finer than the order filter. For every $i=(N,\varepsilon,A_1,\dots,A_N) \in I$ let $B_1,\dots,B_N \in \rb$ given as in the definition of weak containement, and define a map $u_i\colon \Psi \cup \{A_1,\dots,A_N\} \to [[\cg]]$ by saying $u_i(\psi_k)=\ph_k$ and $u_i(A_k)=B_k$. Then the map $u\colon \Psi \cup Malg(Y)\to \Phi \cup Malg(X)$ defined by $u(A) = [u_i(A)]_\ul$ (1) is well defined as $A$ belongs to the domain of $\ul$-ae $u_i$ and (2) extends to a trace preserving homomorphism.

Assume now that $\ra$ is separable, and let $(A_k)_k$ be a dense sequence in the measure algebra of $\ra$. For every integer $N$, let $B_1,\dots,B_N \in \rb$ given as in the definition of weak containement for $(N,1/N,A_1,\dots,A_N)$, and define as above $u_N\colon \Psi \cup \{A_1,\dots,A_N\} \to [[\cg]]$ by saying $u_N(\ph_k)=\psi_k$ and $u_N(A_k)=B_k$. Similarly the map $u\colon \Psi \cup \{A_k\}_k\to \Phi \cup Malg(X)$ defined by $u(A) = [u_i(A)]_\ul$ extends fisrt by continuity to $\Psi \cup Malg(Y)$, and then to a trace preserving homomorphism.
\end{proof}

\subsection{Sofic groupoids}
\label{sect: sofic groupoids}
The notion of sofic \pmp equivalence relations has been introduce by Elek-Lippner \cite{Elek-Lippner-2010}.
We recall the definition of a sofic groupoid, which is a variation of Ozawa's definition of sofic equivalence relations \cite{Ozawa-hyp-sofic}, see also \cite{Cordeiro-2017}.

A \pmp groupoid is said to be finite if it is finite as a set.

\begin{defi}\label{def:sofic_groupoid} A \pmp groupoid $\ch$ over a probability space $(X,\ra,\mu)$ is \defin{sofic} if for every finite subset $K \subset [[\ch]]$ and every $\varepsilon>0$, there is a finite groupoid $\cg$ and a map $\pi \colon [[\ch]] \to [[\cg]]$ satisfying $d(\pi(\ph\psi),\pi(\ph)\pi(\psi)) \leq \eps$ and $|\tau(\ph) - \tau(\pi(\ph))|<\eps$ for every $\ph,\psi \in K$.
\end{defi}
By \cite[Proposition 1.4]{Cordeiro-2017}, we can replace \emph{finite groupoid} by \emph{transitive equivalence relations on a finite set} in the definition. Moreover if $\ra$ is separable, then $\ch$ is sofic if and only if $[[\ch]]$ embeds isometrically (equivalently preserving the trace) into the ultraproduct $\prod_\ul[[\cg_n]]$ for a sequence $\cg_n$ of finite groupoids (resp. finite transitive equivalence relations)\footnote{without this separability assumption, the same holds but we have to allow ultrafilters on larger sets of indices than $\Nmath$; see the proof of Proposition~\ref{prop:relation_weakequivalence_factor} for details.}. By Proposition~\ref{lem:ultraproduct_of_pseudo_full_groups}, this holds if and only 
there are some finite generating graphings $\Phi_n$ of $\cg_n$ such that $[[\ch]]$ embeds isometrically in $[[\cg_\ul]]$ 
where $\cg_\ul$ is the ultraproduct groupoid of the sequence $(\cg_n,\Phi_n)_n$.
 Using Lemma~\ref{Lem: factors vs full pseudogroups},
 we can summarize this discussion for standard groupoids as follows.

\begin{prop} 
\label{prop: sofic vs weak-cont and factors}
For a standard \pmp groupoid $\ch$ with generating graphing $\Psi$, the following properties are all equivalent.
\begin{enumerate}
  \item $\ch$ if sofic,
\item $(\ch,\Psi)$ is weakly contained in the ultraproduct $(\cg_\ul,\Phi_\ul)$ of a sequence of finite graphed groupoids $(\cg_n,\Phi_n)$.
\item $(\ch,\Psi)$ is weakly contained in the ultraproduct $(\RR_\ul,\Phi_\ul)$ of a sequence of finite transitive graphs $(\RR_n,\Phi_n)$.
\item \label{it: sofic=factor-fin-gpoid} $\ch$ is a factor of the ultraproduct $(\cg_\ul,\Phi_\ul)$ of a sequence of finite graphed groupoids $(\cg_n,\Phi_n)$.
\item
\label{it: sofic=factor-fin-trans-graph}
$\ch$ is a factor of the ultraproduct $(\RR_\ul,\Phi_\ul)$ of a sequence of finite transitive graphs $(\RR_n,\Phi_n)$.
\end{enumerate}
\end{prop}
Any such sequence $(\cg_n,\Phi_n)_n$ or  $(\RR_n,\Phi_n)_n$ as in Items \ref{it: sofic=factor-fin-gpoid}. or \ref{it: sofic=factor-fin-trans-graph}. of Proposition~\ref{prop: sofic vs weak-cont and factors} is called a \defin{sofic approximation} of $\ch$.

\section{Group actions and cost}
\label{sect: gp act and cost}

 We  now state the main facts proved in the preceding sections in the context of actions of countable groups. 

 If $a: \Gamma\acting X$ is a \pmp  $\Gamma$-action, we denote by $\cost(a)$ and call it the \defin{cost of the action} the cost of the \pmp groupoid $\cg_{\Gamma\acting X}$. Similarly for a sequence of $\Gamma$-actions $a_n:\Gamma\acting X_n$, we denote by $\ccost_\ul(a_n)$ the combinatorial cost of the sequence of graphed groupoids $(\cg_{\Gamma\acting X_n}, S)$ for any generating subset $S<\Gamma$ (this is independent of the choice of $S$ by definition, see Definition~\ref{def: comb cost, infinite size}).

    Denote\footnote{We introduce this notation $\vraicost_{*}(\Gamma)$ in order to avoid confusion with the cost $\cost(\Gamma)$ of $\Gamma$ seen as a groupoid which is the rank of $\Gamma$.} by $\vraicost_*(\Gamma)$ the infimum of the costs over all free \pmp standard actions of $\Gamma$ (it is realized be some free action \cite[Proposition VI.21]{Gab-cost}) and by
    $\vraicost^*(\Gamma)$ the \defin{supremum cost} (i.e., the supremum of the costs of the free \pmp actions of $\Gamma$). It is realized by any Bernoulli shift action of $\Gamma$ (Ab\'ert-Weiss \cite{Abert-Weiss-2013}).
Recall that the group has \defin{fixed price} when $\vraicost^*(\Gamma)= \vraicost_*(\Gamma)$.
We warn the reader that there is a risk of confusion between $\vraicost(\Gamma)$ (as defined in \cite{Gab-cost}) and the cost $\cost(\Gamma)$ of $\Gamma$ seen as a groupoid, which is the rank of $\Gamma$.

 The following is a corollary of Theorem~\ref{thm: standard factors}.
 \begin{coro}\label{cor:monotonicitygroup}
  Consider a countable group $\Gamma$.
  \begin{enumerate}
  \item The cost of \pmp $\Gamma$-actions is non-decreasing under factors.
 
   \item The cost of any \pmp $\Gamma$-action is the minimum of the cost of all its standard factors. 
   
   \item The infimum cost $\vraicost_{*}(\Gamma)$ equals the infimum of the (groupoid) costs over all \pmp actions of $\Gamma$, free or not, standard or not.
  \end{enumerate}
\end{coro}
We also reformulate Theorem~\ref{thm:costccost} for groups.

\begin{theo}\label{th: equality of lim cost for fixed price}
  Let $\Gamma$ be a finitely generated group. For every sequence $(a_n)_n$ of \pmp $\Gamma$-actions we have:
  \begin{equation}\label{eq: ineqs an, bn}
\cost(a_\ul)=\ccost_\ul((a_n)_n) \geq \lim_{\ul}\cost(a_{n})\geq \liminf_{n\to \infty}\cost(a_{n}) \geq \vraicost_*(\Gamma).
\end{equation}
Moreover if the action $a_\ul$ is essentially free we have $\vraicost^*(\Gamma)\geq \cost(a_\ul)$ and hence if $\Gamma$ has fixed price we obtain 
 $\lim_\ul\cost(a_n)=\ccost_\ul((a_n)_n)=\cost(a_{\ul})=\vraicost_*(\Gamma)$.
\end{theo}

\begin{ques}\label{qst:cont} 
Can the inequality
$\cost(a_\ul)\geq\lim_\ul\cost(a_n)$ in Theorem~\ref{th: equality of lim cost for fixed price} be strict for some finitely generated group $\Gamma$?
\end{ques}
Observe that such an example would provide an example to the {\em fixed price problem} (whether there are countable groups for which $\vraicost_*(\Gamma)\not=\vraicost^*(\Gamma)$ \cite[Question I.8]{Gab-cost}).
 Compare Remark~\ref{rem: counterex revers ineq non-group} where the approximating groupoids are not given by an action of the group
 or Remark~\ref{rem: counterex revers ineq non bounded} where the group is not finitely generated.

We also obtain the following corollary of Theorem~\ref{thm:costccost} for sofic approximations:
\begin{coro}
\label{cor: sofic gp - sofic approx - comb cost}
If $\Gamma$ is a sofic group finitely generated by $S$ and $(\cg_n, \Phi_n)_n$ 
is a sofic approximation to $\Gamma$ (with $S$-labelled graphings), then 
  \begin{equation}
\vraicost^*(\Gamma)\geq \cost(\cg_\ul)=\ccost_\ul((\cg_n)_n) \geq \lim_{\ul}\cost(\cg_{n})\geq \liminf_{n\to \infty}\cost(\cg_{n}) \geq 1.
\end{equation}
If $\Gamma$ has moreover fixed price $\vraicost_*(\Gamma)$, then
$\vraicost_*(\Gamma)=\ccost((\cg_n)_n) $ and hence if $\Gamma$ has fixed price $\vraicost_*(\Gamma)=1$, we obtain 
$\ccost((\cg_n)_n) = \lim_{n\to \infty}\cost(\cg_{n}) = 1$.
\end{coro}
This applies for instance to any sofic approximation (if it exists!) of the Higman group $H_4$ (which has fixed price $=1$).

Recall that if $\Lambda$ has finite index in $\Gamma$, then the action of $\Gamma$ on the coset space $\Gamma/\Lambda$ has (groupoid) cost equal to
\begin{equation}
\cost(\Gamma\acting \Gamma/\Lambda)=1+\frac{\rank(\Lambda)-1}{[\Gamma:\Lambda]}.
\end{equation}

The following theorem is due to Ab\'ert and Nikolov \cite{Abert-Nikolov-12}. They stated it for nested Farber sequences (Definition~\ref{def: Farber sequence}), and their proof indeed applies for general nested sequences and groupoid cost.
We provide a new proof.
\begin{coro}[Ab\'ert-Nikolov {\cite[Theorem 1]{Abert-Nikolov-12}}]\label{cor: Abert-Nikolov}
 Let $(\Gamma_n)_n$ be a nested sequence of finite index subgroups of the finitely generated group $\Gamma$. The (groupoid) cost of the profinite action $\Gamma\acting X_{\mathrm{prof}}:=\projlim \Gamma/\Gamma_n$ is equal to the rank gradient of the sequence plus $1$, that is $\cost(\Gamma\acting X_{\mathrm{prof}})=1+\lim\limits_{n\to \infty}\frac{\rank(\Gamma_n)-1}{[\Gamma:\Gamma_n]}.$
\end{coro}
\begin{proof}
The coset actions $a_n:\Gamma\acting \Gamma/\Gamma_n$ are factors of the profinite action which is itself a factor of the ultraproduct action $a_\ul$ (for any non principal $\ul$),  so by monotonicity under factors we have
$1+\inf_n \frac{\rank(\Gamma_n)-1}{[\Gamma:\Gamma_n]}\geq 
\cost(\Gamma\acting X_{\mathrm{prof}})\geq \cost(a_{\ul})$. By Theorem~\ref{th: equality of lim cost for fixed price}, $\cost(a_{\ul})=1+\lim_\ul\frac{\rank(\Gamma_n)-1}{[\Gamma:\Gamma_n]}$, and the result follows.
\end{proof}

\subsection{Farber sequences}
\label{sect: Farber sequences}

A group $\Gamma$ is \defin{sofic} if it is sofic as a groupoid.

\begin{rema}
\label{rem:soficity and freeness of ultralimit action}
A sequence of maps $\sigma_n:\Gamma\to \Sym (X_n)$ where $X_n$ is a finite set equipped with the equiprobability measure $\mu_n$ defines for any generating subset $S$ of $\Gamma$ an equivalence relation $(\RR_n, \Phi_n)$ 
(by its graphing $\Phi_n=(\sigma_n(s))_{s\in S}$) which is a sofic approximation to $\Gamma$ 
if and only if the ultraproduct $(\RR_\ul, \Phi_\ul)$ (with $\Phi_\ul=(\ph_s)_{s\in S}$ and $\ph_s=[\sigma_n(s)]_\ul$)  generates  a free $\Gamma$-action modulo a null-set, for every non-principal ultrafilter $\ul$.
By Corollary~\ref{cor:dependance_en_graphage}, 
$\RR_\ul$ and its coarse structure is independent of the choice of $S$.
\end{rema}

 \begin{defi}[Farber sequence]
 \label{def: Farber sequence}
 A \defin{Farber sequence} is a sequence 
   $(\Gamma_n)_n$ of finite index subgroups of $\Gamma$ such that every $\gamma\in \Gamma\setminus\{\id\}$ satisfies 
\begin{equation*}
\frac{\left|\left\lbrace g\in \Gamma/\Gamma_n:\ \gamma\in g\Gamma_ng^{-1}
\right\rbrace\right|}{[\Gamma:\Gamma_n]}\underset{n\to \infty}{\longrightarrow} 0.
\end{equation*}
\end{defi}
Originally, Farber \cite{Farber-1998-L2-approx} considered \defin{nested sequences} (i.e., $\Gamma_n\geq \Gamma_{n+1}$ for every $n$) of finite index subgroups in order to extend L\"uck approximation theorem.

\begin{rema} 
\label{rem: Farber vs free ultralim act}
Observe the following:
\begin{enumerate}
\item 
A sequence  $(\Gamma_n)_n$ of finite index subgroups is Farber if and only if for every non-principal ultrafilter $\ul$ the ultraproduct $\Gamma$-action $a_{\ul}$  of the sequence of actions $a_n: \Gamma \acting \Gamma/\Gamma_n$ is essentially free.

\item 
If $S$ is a generating subset $S$ of $\Gamma$, the sequence $(\Gamma_n)_n$ of finite index subgroups is Farber if and only if the $\ul$-ultraproduct of the sequence of associated Schreier graphs $\cg_n:=\mathrm{Sch}(\Gamma/\Gamma_n, S)$ seen as $S$-graphed principal groupoids (i.e., equivalence relations) is a free $\Gamma$-action $b_{\ul}$, for every non principal $\ul$.
In this case, the free $\Gamma$-actions $a_{\ul}$ and $b_{\ul}$ coincide.
\end{enumerate}
\end{rema}

\begin{coro}\label{cor: nested Farber seq}
If $\Gamma$ is finitely generated by $S$ and the sequence $(\Gamma_n)_n$ is a nested Farber sequence, then 
$\cost(a_{\ul})=\ccost_{\ul}(a_n)=\ccost(\mathrm{Sch}(\Gamma/\Gamma_n, S)_n)=\lim\limits_{n\to \infty} \cost(a_n)$.
\end{coro}
\begin{proof} Since $a_n$ is a factor of $a_p$ for $p\geq n$, apply  Proposition~\ref{prop: projective limit and costs} and Remark~\ref{rem: Farber vs free ultralim act}.\end{proof}

\subsection{Sequences of actions and non nested Farber sequences}

Theorem~\ref{th: equality of lim cost for fixed price} applies for instance to non-nested Farber sequences.
Recall that $\ccost$ denotes the original combinatorial cost \cite{Elek-combinatorial-cost-2007}, see Remark~\ref{Rem: Elek comb cost}.
\begin{theo}
\label{cor: rk grad vs cost - non fixed price}
Let $\Gamma$ be a group finitely generated by $S$.
For every (non necessarily nested) Farber sequence $(\Gamma_n)_n$ of finite index subgroups, we have: 
\begin{equation*}
\vraicost^*(\Gamma)\geq \ccost((\mathrm{Sch}(\Gamma/\Gamma_n, S))_n) \geq \liminf\limits_{n\to \infty}\frac{\rank(\Gamma_n)-1}{[\Gamma:\Gamma_n]} +1\geq \vraicost_*(\Gamma)\geq \beta^{(2)}_1(\Gamma)+1.
\end{equation*}
And moreover, 
$
\vraicost^*(\Gamma)\geq \limsup\limits_{n\to \infty} \frac{\rank(\Gamma_n)-1}{[\Gamma:\Gamma_n]}+1.
$
\end{theo}

\begin{proof} Consider the sequence of canonical actions $a_n:\Gamma\acting \Gamma/\Gamma_n$. The cost of the associated groupoid is $\cost(a_n)=\frac{\rank(\Gamma_n)-1}{[\Gamma:\Gamma_n]}+1$ and the ultraproduct action $a_\ul$ is free.

By Theorem~\ref{th: equality of lim cost for fixed price}, we have
\[
\vraicost^*(\Gamma)\geq \cost(a_\ul)=\ccost_\ul((a_n)_n) \geq \lim_{\ul}\cost(a_{n})\geq \liminf_{n}\cost(a_{n}) \geq \vraicost_*(\Gamma)
.\]
By Theorem~\ref{thm:costccost}, $\ccost_\ul((\mathrm{Sch}(\Gamma/\Gamma_n, S))_n)=\cost(b_\ul)=\cost(a_\ul)$
(see Remark~\ref{rem: Farber vs free ultralim act}), for every non-principal $\ul$.
The conclusion follows from (\ref{eq:inequ combcost vs u-combcost}) in Remark~\ref{Rem: Elek comb cost}
and the fact that $\limsup$ and $\liminf$ of a bounded sequence are simply the supremum (resp. the infimum) of the limits along all non principal ultrafilters. The inequality $\vraicost_*(\Gamma)\geq \beta^{(2)}_1(\Gamma)+1$  
is an application of \cite[Corollaire 3.23]{gab_ihes}.
\end{proof}

\begin{coro}[Rank gradient, combinatorial cost and fixed price]
\label{cor: rk grad vs cost fixed price}
If $\Gamma$ is finitely generated by $S$, has fixed price $\vraicost_*(\Gamma)$, and $(\Gamma_n)_n$ is any
(non necessarily nested) Farber sequence, then 
we have: 
\begin{equation*}
\lim\limits_{n\to \infty} \frac{\rank(\Gamma_n)-1}{[\Gamma:\Gamma_n]}+1=\ccost(\mathrm{Sch}(\Gamma/\Gamma_n, S))=\vraicost_*(\Gamma).
\end{equation*} 
\end{coro}

\subsection{Stuck-Zimmer property}


A countable group $\Gamma$ is called a \defin{Stuck-Zimmer group} (resp. \defin{strong Stuck-Zimmer group}) if each \pmp $\Gamma$-action on a standard Borel space satisfies the dichotomy:
\\
-- either it has a non-null invariant set of finite orbits; or 
\\
-- it has a.s. finite stabilizers (resp. trivial stabilizers, i.e., $\Stab_{\Gamma}(x)=\{\id\}$). 

\begin{exam}

The integers $\Zmath$ as well as all Torsion Tarski Monsters (more generally all infinite simple groups with only countably many subgroups) are examples of strong Stuck-Zimmer groups.

The main source for Stuck-Zimmer property examples come from algebraic groups with higher rank conditions such as irreducible lattices in semisimple real Lie groups (see \cite{Stuck-Zimmer-1994} and \cite{Creutz-Peterson-2017}) or some  $\Gamma=\mathbf{G}(\mathcal{O}_S)$ 
for the ring $\mathcal{O}_S$ of $S$-integers associated with some global field $K$ (see \cite[Theorem 10.8]{Creutz-Peterson-2017}). 

If  $\mathbf{G}$ is a simple algebraic group defined over the global field $K$ with $v$-rank
at least two for some place $v$ and such that the $v_{\infty}$-rank  is  at  least  two  for  every  infinite  place $v_{\infty}$, then $\mathbf{G}(K)$ is strongly Stuck-Zimmer
\cite[Corollary 10.9]{Creutz-Peterson-2017}.

This last example is particularly interesting since it produces algebraic examples without Kazhdan property (T).
For instance the non Kazhdan property (T) group $\mathbf{G}(\Qmath)$ is strongly Stuck-Zimmer
when $\mathbf{G}=\mathrm{SO}(q)$ for the quadratic form 
$q=x^2 - a y^2 + pz^2 - w^2 + s^2$ ($p$ prime and $a$ not a square modulo $p$).
\end{exam}

\begin{prop}
\label{prop: Stuck-Zimmer gps and u-limits}
Let $\Gamma$ be a finitely generated Stuck-Zimmer group
and let  $(a_n)_n$ be a sequence of realizable $\Gamma$-actions satisfying
$\mu_{n}(\{x\in X_n: \vert a_n(\Gamma) (x) \vert \leq k\})\underset{n\to \infty}{\longrightarrow} 0$, for every $k>0$.
Then for every torsion free subgroup $\Lambda<\Gamma$ and every non principal ultrafilter $\ul$, the 
restriction to $\Lambda$ of the 
ultraproduct action $a_{\ul}$
is free.
\end{prop}
Recall (Lemma~\ref{lem:realizab ultralim and fix}) that standard \pmp groupoids are realizable and that ultraproducts of realizable groupoids are realizable.

\begin{proof}
The orbits of the ultraproduct $\Gamma$-action $a_{\ul}$ are almost all infinite by Lemma~\ref{lem: ultralim has infinite orbits}.
The same holds for any of its (measurably free) weak equivalent standard factors (given by Theorem~\ref{thm: standard factors}) which thus has finite stabilizers.
Since stabilizers only increase under factors, $a_{\ul}$ also has finite stabilizers. The restriction of $a_{\ul}$  to the torsion free subgroup $\Lambda$ is thus essentially free
\end{proof}

\begin{theo}
\label{th: Stuck-Zimmer gps and seq f i subgroups}
Let $\Gamma$ be a finitely generated Stuck-Zimmer group and let $(\Gamma_n)_n$ be a sequence of finite index subgroups of $\Gamma$ such that $[\Gamma:\Gamma_n]\to \infty$.

\begin{enumerate}
\item 
\label{it: torsion free subg and sofic approx}
For every torsion free subgroup $\Lambda<\Gamma$, the induced action $\Lambda\acting  \Gamma/\Gamma_n$ defines a sofic approximation of $\Lambda$.\footnote{In particular, if $\Gamma$ is torsionfree (or strongly Stuck-Zimmer), then $(\Gamma_n)_n$ is a Farber sequence.}
 \item
\label{it: rk grad vs fixed price}
If $\Gamma$ is torsion free and has fixed price $\vraicost_*(\Gamma)$, then 
\begin{equation}
\lim_{n}
\frac{\left( \rank
(\Gamma_n)-1\right)}{ [\Gamma:\Gamma_n]}=\vraicost_*(\Gamma)-1.
\end{equation}

\item 
If $\Gamma$ has fixed price $=1$ and admits a finite index torsion free subgroup, then 
\begin{equation}
\lim_{n} \frac{\left( 
\rank
(\Gamma_n)-1\right)}{ [\Gamma:\Gamma_n]}=0.
\end{equation}
\end{enumerate}

\end{theo}

As an application of the third item, one recovers for instance Theorem 2 of \cite{Abert-Nikolov-Gelander}: $\Gamma$ is a lattice of a higher rank simple Lie group (thus a Stuck-Zimmer group) and is right angled (thus has fixed price $=1$ \cite{Gab-cost}).

\begin{proof}[Proof of Theorem~\ref{th: Stuck-Zimmer gps and seq f i subgroups}]
The  arguments for each item will hold for every non-principal ultrafilter, thus obtaining genuine limits.

Applying Proposition~\ref{prop: Stuck-Zimmer gps and u-limits} to the sequence of actions $a_n:\Gamma\acting \Gamma/\Gamma_n$, item~\ref{it: torsion free subg and sofic approx} follows (see Remark~\ref{rem:soficity and freeness of ultralimit action}).
Concerning  item~\ref{it: rk grad vs fixed price}, the groupoid cost of $\Gamma\acting  \Gamma/\Gamma_n$ is exactly $1+ \frac{\left( 
\rank
(\Gamma_n)-1\right)}{ [\Gamma:\Gamma_n]}$. One concludes by 
Theorem~\ref{th: equality of lim cost for fixed price}.
As for the last claim, the finite index subgroup $\Lambda$ may be assumed to be normal. It has also fixed price $=1$ (\cite{Gab-cost}). In a weak equivalent standard factor action $b$ of $a_{\ul}$ (Theorem~\ref{thm: standard factors}), we pick a cheap graphing for $\Lambda$ and add representatives of the classes $\Gamma/\Lambda$ defined on arbitrarily small Borel sets that meet all $\Lambda$-orbits. This gives a cheap graphing for the $\Gamma$-action $b$. Thus $a_{\ul}$ has cost $=1$.
 One concludes by 
Theorem~\ref{th: equality of lim cost for fixed price} since $\cost(\Gamma\acting \Gamma/\Gamma_n)=1+\frac{\rank
(\Gamma_n)-1 }{ [\Gamma:\Gamma_n]}\geq 1$.
\end{proof}

\subsection{``Relative'' Stuck-Zimmer property}

We recall that a subgroup of a group $\Gamma$ is called \defin{almost normal} if its normalizer has finite index in $\Gamma$, or equivalently if it has finitely many $\Gamma$-conjugates.

\medskip

A pair $N<\Gamma$ of countable groups is said to satisfy the  \defin{relative Stuck-Zimmer property} if each \pmp $\Gamma$-action on a standard Borel space satisfies the dichotomy:
\\
-- either there is an infinite subgroup $N_0 < N$ that is almost normal in $\Gamma$ and 
that acts trivially on a subset of positive measure; or
\\
--
almost every stabilizer $\Stab_{\Gamma}(x)$ is finite.

\begin{theo}
\label{th: normal subgr f subg --> rel Stuck-Zimmer}
Let $\Gamma$ be a finitely generated group and $N\triangleleft \Gamma$ be an infinite normal subgroup such that
\begin{enumerate}[label=\alph*)]
\item \label{it : assumpt count many sbgr}
$N$ has  countably many subgroups (for instance $N$ finitely generated abelian),\item \label{cond: N and commut of gamma}
for every $\gamma\in \Gamma\setminus N$, the intersection $N_{\gamma}:=\{n\in N:  \gamma n=n \gamma \}$ of $N$ with the centralizer of $\gamma$ has infinite index in $N$.
\end{enumerate}

The following holds:
\begin{enumerate}
\item 
\label{it: rel Stuck-Zimmer}
For every \pmp $\Gamma$-action on a standard probability space, 
\begin{enumerate}[label=(\roman*)]
\item
\label{it: inf normal subgr trivial act}
either the set of $x \in X$ with infinite $N$-stabilizer has positive measure; and in this case there is an infinite subgroup $N_0 < N$ that is almost normal in $\Gamma$ and that acts trivially on a subset of $X$ of positive measure; or
\item 
\label{it: or act is free}
almost every stabilizer is finite and contained in $N$.
\end{enumerate}

\item\label{it: index on N to inft -> sofic}
Assume $N$ is finitely generated torsion free and every infinite subgroup $N_0<N$ which is almost normal in $\Gamma$ has finite index in $N$. Then for any sequence $(\Gamma_n)_n$ of finite index subgroups of $\Gamma$ such that $[N:N\cap \Gamma_n]\to \infty$, the actions $a_n:\Gamma\acting \Gamma/\Gamma_n$ gives rise to a sofic approximation of $\Gamma$.

 If moreover $\Gamma$ has fixed price, then
\begin{equation}
\lim_{n} \frac{
\rank
(\Gamma_n)-1}{ [\Gamma:\Gamma_n]}=\vraicost_*(\Gamma)-1.
\end{equation}

 If moreover $N$ has fixed price $1$ (for instance $N$ amenable), then
\begin{equation}
\lim_{n} \frac{
\rank
(\Gamma_n)-1}{ [\Gamma:\Gamma_n]}=0.
\end{equation}

\end{enumerate}
\end{theo}
\begin{rema}
If $N$ is torsion free, the dichotomy in Theorem~\ref{th: normal subgr f subg --> rel Stuck-Zimmer}-\ref{it: rel Stuck-Zimmer}. becomes simply
\begin{enumerate}[label=(\roman*)]
\item
either $N$ does not act essentially freely
(and in this case there is an infinite subgroup $N_0 < N$ which is almost normal in $\Gamma$ and which acts trivially on a subset of $X$ of positive measure);
\item 
or the $\Gamma$-action is free.
\end{enumerate}

\end{rema}\begin{rema}\label{rem:stuckZimmerRelatif_not_realizable} The conclusion of Theorem~\ref{th: normal subgr f subg --> rel Stuck-Zimmer}-\ref{it: rel Stuck-Zimmer}. remains valid, suitably reworded, for \pmp actions on general (not necessarily standard) probability spaces. Let us say that a \pmp action $\Gamma \acting X$ is \defin{measurably free} if for every $\gamma \in \Gamma \setminus\{1\}$, the bisection $\{\gamma\} \times X$ of the groupoid $\cg_{\Gamma \acting X}$ is measurably free in the sense of Definition~\ref{def: tot. vs meas isotrop, free}. On the opposite, we say that the action $\Gamma \acting X$ is \defin{measurably trivial} on a measurable subset $Y \subset X$ if for every $\gamma \in \Gamma$, the bisection $\{\gamma\} \times Y$ is measurably isotropic as in the same definition. With this terminology, we obtain for example the following generalization of the first step of the theorem:

Let $\Gamma$ be a finitely generated group and $N\triangleleft \Gamma$ be an infinite normal subgroup satisfying conditions \ref{it : assumpt count many sbgr} and \ref{cond: N and commut of gamma} above. For every \pmp $\Gamma$-action on any probability space,
\begin{enumerate}[label=(\roman*)]
\item either there is an infinite subgroup $N_0 < N$ which is almost normal in $\Gamma$ and which acts measurably trivially on a subset of $X$ of positive measure,
\item 
or for positive measure subset $A \subset X$ its measurable stabilizer $\{\gamma \in \Gamma: A \subset \Fix_{\rb}(\gamma)\}$  is a finite subgroup of $N$ (and in particular every element of $\Gamma\setminus N$ acts measurably freely).
\end{enumerate}
In particular if the (groupoid of the) action is realizable, then the conclusion \ref{it: rel Stuck-Zimmer}. of Theorem~\ref{th: normal subgr f subg --> rel Stuck-Zimmer} holds verbatim.
\end{rema}

\begin{proof}
Consider a \pmp action of $\Gamma$ on a standard probability space $(X,\mu)$.
By assumption~\ref{it : assumpt count many sbgr} we have a countable decomposition of $X$ as $X= \sqcup_{L} X_L$ where $L<N$ and $X_L:=\{x\in X: \Stab_{\Gamma}(x)\cap N=L\}$.
\\
(i)
If for some infinite $L$ the $X_{L}$ is non-null, then since $\gamma\cdot X_L=X_{\gamma L \gamma^{-1}}$ the $\gamma\cdot X_L$ are either disjoint or equal and they all have the same measure. It follows that $L$ has finitely many conjugates by $\Gamma$, i.e., that $N_0:=L$ is almost normal in $\Gamma$.
\\
(ii)
Otherwise, $X = \sqcup_{L\textrm{ finite}} X_L$. Let $\gamma\in \Gamma\setminus N$ and consider $\Fix(\gamma):=\{x\in X: \gamma x=x\}$.
Assume $\mu(\Fix(\gamma))>0$. There is a finite subgroup $L<N$ such that $\Fix(\gamma) \cap X_L$ has positive measure. Denote by $Y$ the $\Gamma$-invariant subset $\cup_{\gamma' \in \Gamma} \gamma' X_L$. 
As noted above, the group  $L$ has finitely many $\Gamma$-conjugates $L_1,\dots,L_k$ and $Y=\cup_{i=1}^k X_{L_i}$. Consider $V := \{ n \gamma n^{-1} : n \in N\}$ the set of $N$-conjugates of $\gamma$. Consider $y \in Y$, and $i$ such that $y \in X_{L_i}$. We claim that the set of $N$-conjugates $v$ of $\gamma$ such that $y\in \Fix(v)$ is contained in some $L_i$-coset, and in particular its cardinality is bounded above by $|L|$. Indeed if $y\in \Fix(n_1 \gamma n_1^{-1})\cap \Fix(n_2 \gamma n_2^{-1})\cap X_{L_i}$ for $n_1,n_2 \in N$, then $n_1\gamma n_1^{-1}n_2\gamma^{-1}n_2^{-1}$ fixes $y$ and belongs to $N$ because $N$ is normal. Hence $n_1\gamma n_1^{-1}n_2\gamma^{-1}n_2^{-1}$ must be an element of $L_i$, or in other words $n_1\gamma n_1^{-1} \in L_i n_2\gamma n_2^{-1}$. So we deduce
\[ \sum_{v \in V} \mu(\Fix(v) \cap Y) = \int_Y \sum_{v \in V} 1_{y \in \Fix(v)} d\mu(y) \leq |L| \mu(Y).\]
Remark also that $\Fix(n\gamma n^{-1}) \cap Y=n(\Fix(\gamma) \cap Y)$ and in particular $\mu(\Fix(v)\cap Y)$ does not depend on $v$ in $V$. So the previous inequality is $|V| \mu(\Fix(\gamma) \cap Y) \leq |L| \mu(Y)$, which implies $|V|<\infty$. This is excluded by assumption \ref{cond: N and commut of gamma}.

We note that the above proof works with no change if $(X,\mu)$ is not a standard measure space but the action is assumed to be realizable. If the action is not assumed to be realizable, we obtain Remark~\ref{rem:stuckZimmerRelatif_not_realizable} if in the previous proof we replace $X_L$ by
\[ \left(\cap_{\gamma \in L} s(\Fix_{\rb}(\{\gamma\} \times X))\right) \cap \left(\cap_{\gamma \in N \setminus L} s(M_F(\{\gamma\} \times X))\right)\]
and $\Fix(\gamma)$ by $s(\Fix_{\rb}(\{\gamma\}\times X))$, 
where $\Fix_{\rb}(\gamma)$ and $M_F(\gamma)$ are the measurably isotropic and measurably free parts of the bisection $\{\gamma\} \times X$ of the groupoid $\cg_{\Gamma \acting X}$ given by Proposition~\ref{prop:manyinvo-gen}.

Let us prove item~\ref{it: index on N to inft -> sofic}. We have to prove that under the conditions in item~\ref{it: index on N to inft -> sofic}, the ultraproduct action  $a_\ul$ of $a_n:\Gamma\acting \Gamma/\Gamma_n$ (which by Lemma~\ref{lem:realizab ultralim and fix} is realizable for any non principal ultrafilter) is essentially free. To do so, we apply item~\ref{it: rel Stuck-Zimmer}. In the $N$-action on $\Gamma/\Gamma_n$, all the orbits are of the form $N/(N\cap \gamma \Gamma_n \gamma^{-1})$ for $\gamma$ in $\Gamma$, and therefore all have the same cardinality $[N:N\cap \Gamma_n]$ by normality of $N$. The cardinality of the $N$-orbits thus tends to infinity with $n$ and the ultraproduct $N$-orbits are $\mu_{\ul}$-a.a. infinite by Lemma~\ref{lem: ultralim has infinite orbits}. This implies that for any finite index subgroup $N_0 <N$ (and in particular any infinite subgroup $N_0<N$ which is almost normal in $\Gamma$), the $N_0$-orbits are $\mu_{\ul}$-a.a. infinite. Therefore by item~\ref{it: rel Stuck-Zimmer} and the assumption that $N$ is torsion free, the action $a_\ul$ is essentially free. 

If $N$ has fixed price $=1$, then $\Gamma$ itself has fixed price $=1$ (\cite[Crit\`ere VI.24.(2)]{Gab-cost}). When $\Gamma$ has fixed price, we conclude by Theorem~\ref{th: equality of lim cost for fixed price} since $\cost(\Gamma\acting \Gamma/\Gamma_n)=\frac{\rank
(\Gamma_n)-1 }{ [\Gamma:\Gamma_n]}+1$.
\end{proof}

\begin{theo}
\label{th: SLd rel Stuck-Zimmer}
Let $\Lambda<\SL(d,\Zmath)$ a subgroup whose action on $\Rmath^d$ is \defin{$\Zmath$-strongly irreducible} (there is no finite union of infinite index subgroups of $\Zmath^d$ which is $\Lambda$-invariant).
Let $(\Gamma_n)_n$ be any sequence of finite index subgroups of $\Gamma=\Lambda\ltimes \Zmath^d$ such that $[\Zmath^d:\Gamma_n\cap \Zmath^d]\underset{n\to \infty}{\longrightarrow}\infty$.
Then it defines a sofic approximation of $\Gamma$ and 
\begin{equation}
\lim_{n}
\frac{\left( \rank
(\Gamma_n)-1\right)}{ [\Gamma:\Gamma_n]}=0.
\end{equation}
\end{theo}

\begin{proof}[Proof of Theorem~\ref{th: SLd rel Stuck-Zimmer} (and \ref{th: intro SL2 rel Stuck-Zimmer})]
Set $N= \Zmath^d$. Condition \ref{cond: N and commut of gamma} of Theorem~\ref{th: normal subgr f subg --> rel Stuck-Zimmer} is satisfied since an element $\gamma=\lambda n\in \Lambda\ltimes N$ (with $\lambda\in \Lambda$ and $n\in N$) commutes with $h\in N$ 
iff 
$\lambda n h n^{-1} \lambda^{-1}=h$ iff 
$h$ is a $1$-eigenvector of the associated linear transformation.
This happens on rank $=d$ subgroup iff $\lambda$ acts as the identity.
The irreducibility assumption allows us to use item \ref{it: index on N to inft -> sofic} of Theorem~\ref{th: normal subgr f subg --> rel Stuck-Zimmer}. Indeed, if $N_0 < N$ is infinite and almost normal, then if $E$ denotes its $\Rmath$-linear span, then the union $\cup_{\gamma \in \Gamma} \gamma E$ is $\Gamma$-invariant and finite by almost normality. So the strong irreducibility implies that $E=\Rmath^d$, equivalently that $N_0$ has finite index in $N$.

The hyperbolic element ensures the $\Lambda$-action in Theorem~\ref{th: intro SL2 rel Stuck-Zimmer} is $\Zmath$-strongly irreducible.
\end{proof}

\begin{coro}
Theorem~\ref{th: normal subgr f subg --> rel Stuck-Zimmer}-\ref{it: rel Stuck-Zimmer} allows us to claim that the ergodic IRS's of  $\Gamma=\Lambda\ltimes \Zmath^d$ are all induced by those of $\Lambda_0\ltimes (\Zmath^d/N_0)$ for some finite index $\Lambda_0<\Lambda$ and $\Lambda_0$-invariant finite index $N_0<\Zmath^d$.

\end{coro}

\section{\texorpdfstring{$\ell^2$-Betti numbers}{L2-Betti numbers}}
\label{sect: L2 Betti}
Let $(\cM,\tau)$ be a von Neumann algebra with a faithful normal trace $\tau$ normalized by $\tau(1)=1$ and $\cH$ is a Hilbert $\cM$-module. The von Neumann dimension of $\cH$ is defined as $\tau'(1)$, for $\tau'$ the canonical trace on the commutant $\cM'$ of $\cM$ in $B(\cH)$ \cite{takesakiVol1}. For example, if $\cH = L^2(\cM,\tau) \otimes \cH_1$ for a Hilbert space $\cH_1$, then $\cM'$ is equal to $\cM^{op} \bar \otimes B(\cH_1)$ (where $\cM^{op}$ acts on $L^2(\cM,\tau)$ by right multiplication) and $\tau' = \tau \otimes Tr$. In general, $\cH$ is unitarily equivalent to a subspace of $L^2(\cM,\tau) \otimes \cH_1$ for some $\cH_1$, and $\tau'$ is the restriction of $\tau \otimes Tr$ to the corner $p (\cM^{op} \bar \otimes B(\cH_1)) p$, where $p \in \cM^{op} \bar \otimes B(\cH_1)$ is the orthogonal projection on $\cH$.

We shall use that if $f \colon \cH_1 \to \cH_2$ is a bounded $\cM$-equivariant map between Hilbert $\cM$-modules, then $\mathrm{dim}_\cM(\cH_1)= \mathrm{dim}_\cM(\ker f) + \mathrm{dim}_\cM(\overline{\Ima f})$. This implies that if a sequence
\[ 0 \xrightarrow{f_0} \cH_1 \xrightarrow{f_1} \cH_2 \xrightarrow{f_2} \dots \cH_n \xrightarrow{f_n} 0\]
is weakly exact in the sense that $\ker f_{k+1} = \overline{\Ima f_k}$ for all $k$, then we have the \emph{additivity formula}
\[ \sum_{k, k \textrm{ odd }} \mathrm{dim}_\cM(\cH_k) = \sum_{k, k \textrm{ even }} \mathrm{dim}_\cM(\cH_k).\]

We will introduce now the $\cg$-simplicial complexes in the sense of \cite{gab_ihes} that we will use to compute $\ell^2$-Betti numbers.

\begin{defi}\label{def:gsimplicial}
  Let $\cg$ be a \pmp groupoid over the probability space $(X,\mu)$. A $\cg$-\defin{simplicial complex} $\Sigma$ (in the sense of \cite{gab_ihes}) is given by a fibred space $(\Sigma^{(i)}, \pi_i)$ over $X$ for every non negative integer $i$ such that  
\begin{enumerate}
  \item \label{defaction} 
  the groupoid $\cg$ acts on each $\Sigma^{(i)}$ and $\Sigma^{(0)}$ admits a $\cg$-fundamental domain
  
  \item\label{item:subfibred} for every $i\geq 1$ the fibred space $\Sigma^{(i)}$ (the $i$-\textit{cells}) is contained as a $\cg$-fibred space in $\Sigma^{(0)}*\stackrel{i+1}{\ldots}*\Sigma^{(0)}$,
  
  \item\label{defboundary} 
  if $(v_0,\ldots,v_{i+1})\in\Sigma^{(i+1)}$, then $(v_0,\ldots,v_{k-1},v_{k+1},v_{i+1})\in\Sigma^{(i)}$ for every $k$,
  
  \item 
  $\Sigma^{(i)}$ is invariant under the permutation group $\Sym(i+1)$ which acts on the coordinates of  $\Sigma^{(0)}*\ldots*\Sigma^{(0)}$ and this action commutes with the $\cg$-action,
  
  \item if $(v_0,\ldots,v_{i})\in\Sigma^{(i)}$, then for every $j\neq k$ we have that $v_j\neq v_k$.
\end{enumerate}
The $\cg$-action on $\Sigma^{(0)}*\stackrel{i+1}{\ldots}*\Sigma^{(0)}$ in \ref{defaction}. is the action described in Example~\ref{example:fibredProductAction}.

  Let $\cg$ be a \pmp groupoid and let $\Phi$ be a generating graphing of $\cg$. We  say that the $\cg$-simplicial complex $\Sigma$ is $(\Phi,L)$-\defin{ULB} (\textit{uniformly locally bounded}) if the fibered space $\Sigma^{(0)}$ admits a finite fundamental domain $\cd^0$ with $H(\cd^0)\leq L$ and $\Sigma^{(1)}$ has a fundamental domain $\cd^1$ with $\cd^1\subset B_{\ell_\Phi}(L)\cd^0*\cd^0$, that is each edge in $\Sigma^{(1)}$ corresponds to a pair $(v_1,v_2)\in \Sigma^{(0)}*\Sigma^{(0)}$ for which there exist $\Phi$-words $w_1, w_2$ with the properties that $w_iv_i\in \cd^0$ and $w_1^{-1}w_2$ has length $\leq L$.
\end{defi}

We  say that the $\cg$-simplicial complex $\Sigma$ is \defin{ULB} if it is $(\Phi,L)$-ULB for some generating graphing $\Phi$ of $\cg$ and some integer $L$. Note that if $\Sigma$ is ULB, then every skeleton $\Sigma^{(i)}$ has a finite-height fundamental domain.

For each $x\in X$ we denote by $\Sigma[x]$ the (maximal) simplicial complex contained in $\Sigma$ whose $0$-cells are $\pi^{-1}_0(x)$. This simplicial complex has at most countably many cells in each dimension. 

\begin{defi}
    Let $\cg$ be a \pmp groupoid and let $\Sigma$ be a $\cg$-simplicial complex. We say that $\Sigma$ is $k$-\defin{connected} if for almost every $x\in X$ the simplicial complex $\Sigma[x]$ has trivial homotopy groups up to dimension $k$.
\end{defi}

\begin{exam}\label{ex:universal}
  Let $(\cg,\Phi)$ be a graphed \pmp groupoid over the probability space $(X,\mu)$. Fix a natural number $L$. For $i\in\Nmath$ we define 
  \begin{align*}
    \Sigma_{\Phi,L}^{(i)}:=&\left\{(g_0,\ldots,g_i)\in(\cg,t)*\stackrel{i+1}{\ldots}*(\cg,t):\ \forall j\neq k\text{ we have }1\leq \ell_\Phi(g_j^{-1}g_k)\leq L\right\},\\
    \cd^i:=&\left\{(g_0,\ldots,g_i)\in \Sigma_{\Phi,L}^{(i)}
    : g_0\in X\right\}.  
  \end{align*}
Then $\Sigma_{\Phi,L,k}:=(\Sigma_{\Phi,L}^{(i)})_{i\leq k}$ is a $\cg$-simplicial complex where the action is given by the formula $g(g_0,\ldots,g_i)=(gg_0,\ldots,gg_i)$ and $\cd^i\subset \Sigma_{\Phi,L}^{(i)}$ are fundamental domains satisfying the above conditions and hence $\Sigma_{\Phi,L,k}$ is $(\Phi,L)$-ULB. Note also that $\cup_L \Sigma_{\Phi,L,k}$ is $k$-connected.
\end{exam}

Let $\cg$ be a \pmp groupoid over the probability space $(X,\mu)$ and let $\Sigma$ be a uniformly bounded $\cg$-simplicial complex. Then the Hilbert space $\LLh^2(\Sigma^{(i)})$ is an $\lvN$-module and carries a unitary representation of $\Sym(i+1)$ which commutes with the $\lvN$-action. We denote by $C_i^{(2)}(\Sigma)$ its subspace consisting of functions such that $\sigma \cdot f = \varepsilon(\sigma) f$ where $\varepsilon \colon \Sym(i+1)\to \{-1,1\}$ is the signature. We also define the \defin{boundary operators} $\partial_i \colon C_i^{(2)}(\Sigma) \to C_{i-1}^{(2)}(\Sigma)$ by (using the notation \ref{defboundary} of Definition~\ref{def:gsimplicial})
\begin{equation}\label{eq:def_of_boundary}  (\partial_i f)(v_0,\dots,v_{i-1}) := \sum_{v_i: (v_0,\dots,v_i) \in \Sigma^{(i)}} f(v_0,\dots,v_i),\end{equation}
which allow us to define the following chain complex
\[ \ldots C_{i+1}^{(2)}(\Sigma)\stackrel{\partial_{i+1}}{\rightarrow}  C_{i}^{(2)}(\Sigma)\stackrel{\partial_{i}}{\rightarrow} C_{i-1}^{(2)}(\Sigma)\ldots \]

The $i$-th \defin{reduced homology} of $\Sigma$ is the $\lvN$-module
\[ \overline H_i^{(2)}(\Sigma) := \ker \partial_i / \overline{\Ima \partial_{i+1}}.\]

The $i$-th $\ell^2$-Betti number of $\Sigma$ is the von Neumann dimension of the above homology, $\beta^{(2)}_i(\Sigma,\cg):=\mathrm{dim}_{\LL(\cg)}(\overline H_i^{(2)}(\Sigma))$. As explained after Definition 3.6 of \cite{gab_ihes}, if we define the Laplacian $\Delta_i:=\partial_i^*\partial_i+\partial_{i+1}\partial_{i+1}^*$, then $\overline H_i^{(2)}(\Sigma)$ is naturally isomorphic to $\ker\Delta_i$ and in particular $\beta^{(2)}_i(\Sigma,\cg)=\mathrm{dim}_{\LL(\cg)}(\ker(\Delta_i))$.

\begin{theo}[Gaboriau]\label{thm:bettibound}
 Let $\cg$ be a \pmp groupoid over a probability space and $k\in\Nmath$. 
 All the uniformly bounded and $k$-connected $\cg$-simplicial complexes have the same  $k$-th $\ell^2$-Betti number.
\\
This common value is denoted by $\beta^{(2)}_k(\cg)$.

 Furthermore, for a \pmp (free or not) action $\Gamma\acting X$ of a countable group,  
  the group and the groupoid of the action $\cg_{\Gamma\acting X}$ (see  Example~\ref{ex: groupoid from pmp action}) have the same $\ell^2$-Betti numbers
 $\beta^{(2)}_k(\cg_{\Gamma\acting X})=\beta^{(2)}_k(\Gamma)$.
 
\end{theo}

\begin{rema}
  Gaboriau proved the theorem  for equivalence relations over standard Borel spaces, the case of general standard \pmp groupoids was considered by Sauer (adapted by Takimoto in our langage), see \cite{Sauer-Betti-groupoids-2005,Takimoto-cost-L2-2005}.
  The proofs from \cite{gab_ihes, Takimoto-cost-L2-2005} also work in our non standard setting. Alternatively one can use the existence of a measurably free standard factor, Theorem~\ref{thm: standard factors}.
See also \cite{Aaserud-Popa}. 
\end{rema}

\subsection{\texorpdfstring{Ultraproducts of $\cg$-simplicial complexes}{Ultraproducts of G-simplicial complexes}}
\label{sect: ultraprod of cg-simplicial complexes}

We will now define the ultraproduct of a sequence of $\cg$-\defin{simplicial complex}.

\begin{prop}\label{prop:simpliciallimit}
   Let $(\cg_n,\Phi_n)_n$ be a sequence of \pmp graphed groupoids over the probability spaces $(X_n,\mu_n)$, let $L$ be a natural number and for every $n\in\Nmath$ let $\Sigma_n$ be a $(\Phi_n,L)$-ULB $\cg_n$-simplicial complex. For every $i\in\Nmath$ consider the ultraproduct $\Sigma_\ul^{(i)}$ as in Definition~\ref{def:ultraaction} with respect to the fundamental domains $\cd_n^i$. Then $\Sigma_\ul:=(\Sigma_\ul^{(i)})_i$ is a $(\Phi_\ul,L)$-ULB $\cg_\ul$-simplicial complex.
\end{prop}
\begin{proof}
  By construction for every $i\in\Nmath$ the fibred space $\Sigma_\ul^{(i)}$ is a $\cg_\ul$ fibred space. We have to show that $\Sigma_\ul^{(i)}$ is a measurable subset of $\Sigma_\ul^{(0)}*\stackrel{i+1}{\ldots}*\Sigma_\ul^{(0)}$. For this, note that if $v_\ul\in \Sigma_\ul^{(i)}$, then we have $v_\ul=[(v_n^0,\ldots,v_n^i)]_\ul$. There is a sequence $g_n \in \cg_n$ such that $\ell_{\Phi_n}(g_n)$ is bounded and such that $g_nv_n^k\in B_{\ell_\Phi}(L)\cd_n^0$ for every $k\in\{0,\ldots,i\}$ and $\ulae$. So that for every $k$ we have $[g_n]_\ul[v_n^k]_\ul=[g_nv_n^k]_\ul\in [B_{\ell_{\Phi_n}}(L)\cd_n^0]_\ul\subset \Sigma_\ul^{(0)}$ and hence \[[g_n]_\ul v_\ul\in [B_{\ell_{\Phi_n}}(L)\cd_n^0]_\ul*\ldots*[B_{\ell_{\Phi_n}}(L)\cd_n^0]_\ul\subset \Sigma_\ul^{(0)}*\stackrel{i+1}{\ldots}*\Sigma_\ul^{(0)} .\] 

If for every $i,n\in\Nmath$ we let $\cd_n^i$ be the fundamental domain as in the definition of uniform local boundedness, then their ultraproducts $\cd_\ul^i$ are fundamental domains of $\Sigma_\ul^{(i)}$ which satisfy the required conditions and hence $\Sigma_\ul$ is a  $(\Phi_\ul,L)$-ULB $\cg_\ul$-simplicial complex.
\end{proof}

\begin{rema}[Benjamini-Schramm convergence]
\label{rem: BS-CV of complexes}
  Let $(\cg_n,\Phi_n)_n$ be a sequence of \pmp bounded size graphed groupoids over the probability spaces $(X_n,\mu_n)$.
If $(\Sigma_n)_n$ is a sequence of $(\Phi_n,L)$-ULB $\cg_n$-simplicial complexes which is convergent in the sense of Benjamini-Schramm, then all the ultraproducts $\Sigma_\ul^{(i)}$, for any non-principal ultrafilter $\ul$, have the same local statistics. I.e., the measure of each isomorphism class of bounded  simplicial complex rooted at the fundamental domain. 
This can be taken as a definition for Benjamini-Schramm convergent sequence. Similarly to Remark~\ref{rem: local-global cv} this is independent of the marking (labels and orientation on the graphings).
The same holds for all their weakly-equivalent standard factors.
\end{rema}

We want now to understand under which conditions $\Sigma_\ul[x_\ul]$ is $k$-connected in order to apply Gaboriau's theorem (Theorem~\ref{thm:bettibound}). Let $(\cg,\Phi)$ be a graphed \pmp groupoid, let $\Sigma$ be a $(\Phi,L)$-ULB $\cg$-simplicial complex and let $\cd^i\subset B_{\ell_\Phi}(L)\cd^0*\ldots * B_{\ell_\Phi}(L)\cd^0\subset \Sigma_n^{(i)}$ be a fundmental domain. Let $x\in X$ and let $N\in\Nmath$ be a natural number. We define $\Sigma[x]\bigr|_N$ to be the topological subspace of $\Sigma[x]$ consisting of (the geometric realization of) those simplexes such that their $\ell_{\Phi,\cd^i}$-length is less than $N$ (for the appropriate $i$). We say that $x\in X$ is $(N,M;k)$-\defin{connected} if every homotopy sphere of dimension less than $k$ which is lying in $\Sigma[x]\bigr|_N$ is homotopic inside $\Sigma[x]\bigr|_M$ to a point.  

\begin{prop}\label{prop:contract}
  Let $(\cg_n,\Phi_n)_n$ be a sequence of \pmp graphed groupoids over the probability spaces $(X_n,\mu_n)$, let $L\in\Nmath$ and for every $n\in\Nmath$ let $\Sigma_n$ be a $(\Phi_n,L)$-ULB $\cg$-simplicial complex. Let $\Sigma_\ul$ be the $(\Phi_\ul,L)$-ULB $\cg_\ul$-simplicial complex constructed in Proposition~\ref{prop:simpliciallimit}. If for every $\eps>0$ and $N\in\Nmath$ there exists $M$ and a sequence of measurable subsets $A_n\subset X_n$ such that for any $x\in A_n$ we have that $\Sigma_n[x]$ is $(N,M;k)$-connected and $\lim_\ul\mu(A_n)\geq 1-\eps$, then $\Sigma_\ul[x]$ is almost surely $k$-connected.
\end{prop}

We will call such a sequence of simplicial complexes \defin{asymptotically uniformly }$k$-\defin{connected}.

\begin{proof}
 It is enough to observe that for every $M\in\mathbb N$ and $[x_n]_\ul\in X_\ul$ the sequence $\left(\Sigma_n[x]\bigr|_M\right)_{n\in \Nmath}$ is a sequence of uniformly finite simplicial complexes and hence it is constant $\ul$-almost surely. 
\end{proof}

\subsection{\texorpdfstring{Limits of $\ell^2$ Betti numbers: uniformly bounded case}{Limits of L2 Betti numbers: uniformly bounded case}}

\begin{theo}\label{thm:limitbounded}
  For every $n\in\Nmath$ let $(\cg_n,\Phi_n)$ be a {\bf sofic} graphed \pmp groupoid over $X_n$, let $\Sigma_n$ be a $(\Phi_n,L)$-ULB $\cg_n$-simplicial complex for some fixed $L\in\mathbb N$. Then
\[\lim_\ul\beta^{(2)}_i(\Sigma_n,\cg_n)=\beta^{(2)}_i(\Sigma_\ul,\cg_\ul).\]

In particular if the sequence of $\cg_n$-simplicial complexes $\Sigma_n$ is asymptotically uniformly $k$-connected, then for every $i\leq k$
\[\lim_\ul\beta^{(2)}_i(\Sigma_n,\cg_n)=\beta^{(2)}_i(\cg_\ul).\]
\end{theo}

When $X_n$ is a finite set and $\cg_n$ is the transitive equivalence relation on $X_n$, then the quantity $\beta^{(2)}_i(\Sigma_n,\cg_n)$ corresponds to the average of the usual Betti numbers $b_i(\Sigma_n[x])$ of the finite complexes $\Sigma_n[x]$ for $x\in X_n$ (see \cite[Proposition 2.3]{berggab}),
\[\beta^{(2)}_i(\Sigma_n,\cg_n)=\frac{1}{|X_n|}\sum_{x\in X_n} b_i(\Sigma_n[x]).\]

As a corollary of Theorem~\ref{thm:limitbounded} we obtain the following extension of the approximation theorems of L\"uck and Farber \cite{Luc94b,Farber-1998-L2-approx}.

\begin{coro}\label{cor:approxgroups}
  Let $\Gamma$ be a countable group acting freely cocompactly on a $k$-connected simplicial complex. Let $(\Gamma_n)_n$ be a (non necessarily nested)
  Farber sequence of finite index subgroups  (Definition~(\ref{def: Farber sequence})). Then for every $i\leq k$ we have
\[\lim\limits_{n\to \infty}\frac{b_i(\Gamma_n)}{[\Gamma:\Gamma_n]}=\beta^{(2)}_i(\Gamma).\] 
\end{coro}
More generally, we obtain the following generalization of \cite[Th\'eor\`eme 3.1]{berggab} for non nested, non Farber sequences:
\begin{coro}
\label{Cor: non Farber sq}
  Let $(\Gamma_n)_n$ be any sequence of finite index subgroups of a finitely generated group $\Gamma$.
  Assume $\Gamma$ acts freely cocompactly on the simplicial complex $\Sigma$ such that $\Gamma\backslash \Sigma$ is a simplicial complex.
Let $a_{\ul}:\Gamma\acting X_{\ul}$ be the ultraproduct of the actions $a_n:\Gamma\acting \Gamma/\Gamma_n$ and $\cn_\ul:=\{(x,\gamma)\in X_{\ul}\times \Gamma: \gamma x=x\}  $ its maximal totally isotropic subgroupoid.
  
  Then for every $i\leq k$ we have
\[\lim\limits_{n\in \ul}\frac{b_i(\Gamma_n\backslash \Sigma)}{[\Gamma:\Gamma_n]}=\beta^{(2)}_i(\cn_\ul \backslash \left(X_{\ul}\times \Sigma\right), \RR_{\ul}),\] 
where $\RR_\ul=\cg_\ul/\cn_\ul$ is the \pmp equivalence relation of the ultraproduct action $a_\ul$ and $\cn_\ul \backslash \left(X_{\ul}\times \Sigma\right)$ is the $\RR_{\ul}$-simplicial complex defined in Example~\ref{ex: quotient fibred space}.
\end{coro}

\begin{proof}[Proof of Corollary~\ref{Cor: non Farber sq}]
This is a simple application of Theorem~\ref{thm:limitbounded}: Set \begin{inparaenum}[(a)]
\item $X_n=\Gamma/\Gamma_n$, \item  $\cg_n$ the equivalence relation of $a_n$ graphed using a finite generating set $S\subset \Gamma$, 
\item  $\cn_n:=\{(x,\gamma)\in X_{n}\times \Gamma: \gamma x=x\} =\coprod_{x\in X_n} \mathrm{Stab}_{a_n}(x) $ the maximal totally isotropic subgroupoid of the groupoid of $a_n$, 
\item $\Sigma_n:= \cn_n\backslash\left( X_n\times \Sigma\right)=\coprod_{x\in X_n} \mathrm{Stab}_{a_n}(x)\backslash \Sigma$
\end{inparaenum}
and observe as in 
\cite[Proposition 2.3]{berggab} that 
$\beta^{(2)}_i(\Sigma_n,\cg_n)=\frac{1}{|X_n|}\sum_{x\in X_n} b_i(\Sigma_n[x])=\frac{b_i(\Gamma_n\backslash \Sigma)}{[\Gamma:\Gamma_n]}$. Then, choosing fundamental domains $\cd^{0}, \cd^{1}$ for the $\Gamma$-action on $\Sigma^{(0)}$ (resp. $\Sigma^{(1)}$),  observe that the ultraproduct $\Sigma_\ul$ (Proposition~\ref{prop:simpliciallimit}) of the $\Sigma_n$ with fundamental domains the image $\cd^{\epsilon}_n$ of $X_n\times \cd^{\epsilon}$ (for $\epsilon=0,1$)   in $\cn_n\backslash\left( X_n\times \Sigma\right)$ is precisely $\cn_\ul \backslash (X_{\ul}\times \Sigma)$.
\end{proof}

\begin{proof}[Proof of Corollary~\ref{cor:approxgroups}]
Under the assumption the sequence is Farber, the ultraproduct action is essentially free. We thus recover in Corollary~\ref{Cor: non Farber sq} the $\ell^2$-Betti number of the group by Theorem~\ref{thm:bettibound}, for any non principal ultrafilter.
\end{proof}

The proof of Theorem~\ref{thm:limitbounded} relies on the following fundamental lemma, due to Thom \cite{thom08} under the additional assumption that $\sup_n \|x_n\|<\infty$. 

\begin{lemm}\label{lem:lueck_approximation} Let $(\cM,\tau)$ be a finite von Neumann algebra and $x \in \cM$, let $x_n \in M_{k_n}(\Zmath)$ such that $x=x^*$ and $x_n=x_n^*$ for all $n$ and assume that $\tau(P(x)) = \lim_n \frac 1 {k_n} Tr( P(x_n))$ for every polynomial $P$. Then 
\[ \tau(\chi_{\{0\}}(x))  = \lim_n \frac 1 {k_n} Tr\left( \chi_{\{0\}}(x_n)\right).\]
\end{lemm}
\begin{proof}
Denote by $\mu_n$ the spectral measure of $x_n$ and $\mu$ the spectral measure of $x$. These are the compactly supported probability measures on $\Rmath$ such that $\tau(f(x))=\int fd\mu$ and $\frac{1}{k_n} Tr (f(x_n)) = \int f d\mu_n$ for every Borel function $f$ on $\Rmath$. By classical approximation considerations \cite[Theorem 30.2]{billinglsey}, the assumption that $\tau(P(x)) = \lim_n \frac 1 {k_n} Tr ( P(x_n))$ implies in particular that $\int f d\mu = \lim_n \int f d\mu_n$ for every bounded continuous function on $\Rmath$.  We want to replace $f$ by the indicator function of $\{0\}$, and for this it's enough to prove that $\lim_{\varepsilon \to 0} \sup_n \mu_n([-\varepsilon,\varepsilon]\setminus \{0\}) = 0$.

Consider the Borel function $\varphi(t) = \log |t|$ if $t \neq 0$ and $\varphi(t)=0$ if $t=0$. Then \[\int \varphi d\mu_n = \frac 1 {k_n} \sum_{\lambda \neq 0} \log |\lambda|=\frac 1 {k_n} \log\left(\prod_{\lambda\neq 0}|\lambda|\right),\] where the sum and the product range over all the non zero eigenvalues of $x_n$, counted with multiplicity. The value $\prod_{\lambda\neq 0}|\lambda|$ is one of the coefficients of the characteristic polynomial of $x_n$ and hence it is an integer number. Therefore $ \int_\Rmath \varphi(t) d\mu_n(t)\geq 0$.

Using the inequality $\varphi(t) \leq   (\log \varepsilon) \chi_{[-\varepsilon,\varepsilon]\setminus \{0\}}(t)  +  t^2$, we obtain the inequality
\[ 0 \leq (\log \varepsilon) \mu_n([-\varepsilon,\varepsilon]\setminus \{0\}) + \frac{1}{k_n} Tr(x_n^2).\]
This implies that $\lim_{\varepsilon \to 0} \sup_n \mu_n([-\varepsilon,\varepsilon]\setminus \{0\}) = 0$ and concludes the proof.
\end{proof}

Note that, like in \cite{thom08}, the same proof works to show that for every algebraic integer $\theta$,
\[ \tau(\chi_{\{\theta\}}(x))  = \lim_n \frac 1 {k_n} Tr( \chi_{\{\theta\}}(x_n)).\]
One just replaces the function $\varphi(t)$ by $\log|P(t)|$ if $P$ is a monic polynomial with integer coefficients having $\theta$ as a root.

Let us say that a self-adjoint element $x$ of a finite von Neumann algebra is \defin{approximated by integer matrices} if there is a sequence $(x_n)_n$ as in the statement of Lemma~\ref{lem:lueck_approximation}. For example, it follows from Definition~\ref{def:sofic_groupoid} that if $\cg$ is a sofic \pmp groupoid, then every self-adjoint matrix whose entries are integer combinations of elements in the full pseudogroup, seen as an element of $M_k(\LL'(\cg))$, is approximated by integer matrices. By a diagonal argument, we deduce the following more general statement.
\begin{lemm}\label{lem:lueck_approximation_sofic} Let $(\cM,\tau)$, and for every integer $n$, $(\cM_n,\tau_n)$ be finite von Neumann algebras and $x \in \cM$, $x_n \in \cM_n$ such that $x=x^*$ and $x_n=x_n^*$ for all $n$ and assume that $\tau(P(x)) = \lim_n \tau_n(P(x_n))$ for every polynomial $P$. Assume also that each $x_n$ is approximated by integer matrices in the sense above. Then 
\[ \tau(\chi_{\{0\}}(x))  = \lim_n\tau_n\left( \chi_{\{0\}}(x_n)\right).\]
And therefore, this holds with both $\lim_n$ replaced by $\lim_\ul$.
\end{lemm}

\begin{proof}[Proof of Theorem~\ref{thm:limitbounded}]
We first deal with the particular case when $X_n$ is a finite set with the uniform probability measure, and $\cg_n$ is an equivalence relation on $X_n$. Let $i \geq 0$. For every $n\in\Nmath$ consider the $i$-th Laplacian $\Delta_{i,n}: C_{i}^{(2)}(\Sigma_n)\rightarrow  C_{i}^{(2)}(\Sigma_n)$. Observe that $\Delta_{i,n}$ restricted to each simplicial complex $\Sigma_n[x]$ is a classical Laplacian and therefore it can be expressed as a matrix with integer coefficients. The complex dimension of its kernel equals the sum of the $i$-th Betti number of $\Sigma_n[x]$ for $x\in X_n$, therefore $|X_n|^{-1}\dim_{\Cmath}(\ker(\Delta_{i,n}))=\beta^{(2)}_i(\Sigma_n,\cg_n)$. By Lemma~\ref{lem:ultraproduct_of_pseudo_full_groups} and Lemma~\ref{lem:ultraproductactions}, we have that $C_{i}^{(2)}(\Sigma_\ul)$ embeds into the metric ultraproduct of the Hilbert spaces $C_{i}^{(2)}(\Sigma_n)$ and under this embedding the $i$-th Laplacian $\Delta_{i,\ul}$ corresponds to $[\Delta_{i,n}]_\ul$. 
Thus  $\tau_1(P(\Delta_{i, \ul}))= \tau_2(P([\Delta_{i,n}]_\ul))= \lim_\ul \frac 1 {\vert X_n\vert} Tr( P(\Delta_{i,n}))$ for every polynomial $P$, where $\tau_1$ is the trace on the commutant of $\LL(\cg_\ul)$ acting on $C_{i}^{(2)}(\Sigma_\ul)$ and 
 $\tau_2$ is the trace on the commutant of $[\LL(\cg_n)]_\ul$ acting on $[C_{i}^{(2)}(\Sigma_n)]_\ul$.
Therefore we can use Lemma~\ref{lem:lueck_approximation} to complete the proof. 
In the general case when $\cg_n$ is only assumed to be sofic, we proceed similarily.
The Laplacian is a matrix whose entries are integer combinations of elements in the full pseudogroup $[[\cg_n]]$, thus approximated by integer matrices. We then apply Lemma~\ref{lem:lueck_approximation_sofic} instead of Lemma~\ref{lem:lueck_approximation}.

\end{proof}

\subsection{\texorpdfstring{Limits of $\ell^2$ Betti numbers: the general case}{Limits of L2 Betti numbers: the general case}}

We will now try to understand what happens in Theorem~\ref{thm:limitbounded} if we do not assume that the sequence of complexes $\Sigma_n$ is $(\Phi_n,L)$-ULB. For doing so we will need to recall some further notation from \cite{gab_ihes}. 

Let $\cg$ be a \pmp groupoid and let $\Sigma_1\subset \Sigma_2$ be $\cg$-simplicial complexes. The inclusion of complexes induces an inclusion of the chain complexes $C_*^{(2)}(\Sigma_1)$ in $C_*^{(2)}(\Sigma_2)$ which induces a $\lvN$-module map $\overline H_i^{(2)}(\Sigma_1) \to \overline H_i^{(2)}(\Sigma_2)$. The von Neumann dimension over $\lvN$ of the closure of the image of this map is denoted $\nabla_i(\Sigma_1,\Sigma_2)$. We recall the following theorem, which follows from  \cite[Proposition 3.9, Th\'eor\`eme 3.13]{gab_ihes} for equivalence relations and \cite[Proposition 3.7]{Takimoto-cost-L2-2005}.

\begin{theo}[Gaboriau]
  Let $\cg$ be a \pmp groupoid over a standard probability space, let $\Sigma$ be a $\cg$-simplicial complex which is a countable union of ULB $\cg$-simplicial sub-complexes, $\Sigma=\cup_L\Sigma_L$. Then the quantity 
  \begin{equation}\label{eq:formula_for_betti_number} \lim_{L} \lim_{L' \geq L} \nabla_i(\Sigma_{L},\Sigma_{L'})\end{equation}
  is independant from the choice of the $(\Sigma_L)_L$, and is denoted $\beta^{(2)}_i(\Sigma,\cg)$, the $i$-th $\ell^2$-Betti number of $\Sigma$. It does not depend on $\Sigma$ as long as $\Sigma$ is $i+1$-connected, and coincides with the $i$-th Betti number of $\Gamma$ when $\cg$ comes from a \pmp action of $\Gamma$.
\end{theo}

We will prove the following theorem.

\begin{theo}\label{thm:limitgeneral}
  For every $n\in\Nmath$ let $(\cg_n,\Phi_n)$ be a finite (or more generally sofic) graphed groupoid over $X_n$, let $\Sigma_n$ be a $\cg_n$-simplicial complex and consider $(\Phi_n,L)$-ULB simplicial complexes $\Sigma_{n,L}$ such that $\Sigma_n=\cup_{L\in\mathbb N} \Sigma_{n,L}$. Let us denote by $\Sigma_{\ul,L}$ the $\cg_\ul$-simplicial complex which is the ultraproduct of the sequence $(\Sigma_{n,L})_n$ and set $\Sigma_\ul:=\cup_L\Sigma_{\ul,L}$. Then for every $i \geq 0$ and $L' \geq L$,
\[ \nabla_i(\Sigma_{\ul,L},\Sigma_{\ul,L'}) = \lim_\uu \nabla_i(\Sigma_{n,L},\Sigma_{n,L'}).\]
In particular,
\[ \beta^{(2)}_i(\Sigma_\ul,\cg_\ul) = \lim_{L} \lim_{L' \geq L} \lim_\uu \nabla_i(\Sigma_{n,L},\Sigma_{n,L'}).\]
\end{theo}

\begin{rema}
  Even if for every $n\in\mathbb N$ the simplicial complex $\Sigma_n$ is $k$-connected we can not deduce that $\Sigma_\ul$ is $k$-connected. However if we assume that for all natural number $N$ and for almost every $[x_n]_\ul\in X_\ul$ there exists $L_0$ such that for all $n$ in a subset of $\ul$ and $L\geq L_0$ we have that $\Sigma_{n,L}[x_n]\bigr|_N=\Sigma_{n,L_0}[x_n]\bigr|_N$, then Proposition~\ref{prop:contract} will tell us that $\Sigma_\ul$ is $k$-connected.
\end{rema}

\begin{coro}
 For every $n\in\Nmath$ let $(\cg_n,\Phi_n)$ be a finite (or sofic) graphed \pmp groupoid.  Let $\Sigma_{n,L}$ be the simplicial complex constructed in the Example~\ref{ex:universal}. Then \[\beta^{(2)}_i(\cg_\ul) = \lim_{L} \lim_{L' \geq L} \lim_\uu \nabla_i(\Sigma_{n,L},\Sigma_{n,L'}).\]
\end{coro}
\begin{proof}
 The corollary just follows from the fact that $\Sigma_{\ul,L}$ is again the simplicial complex constructed in the Example~\ref{ex:universal} and hence $\Sigma_\ul=\cup_L\Sigma_{\ul,L}$ is contractible. 
\end{proof}

Let us now prove Theorem~\ref{thm:limitgeneral}.

\begin{lemm}\label{lem:vNdimension} Let $\cM$ be a von Neumann algebra with a finite faithful normal trace $\tau$. Let $\cH$ be a Hilbert $\cM$-module with finite von Neumann dimension and $V,W$ (not necessarily closed) sub $\cM$-modules of $\cH$. Then
\[ \overline{V \cap W} = \overline{V} \cap \overline{W}.\]
\end{lemm}
\begin{rema} It is easy to see that the equality $\overline{V \cap W} = \overline{V} \cap \overline{W}$ is not true for arbitrary subspaces of $\ell^2$. Therefore, already for $\cM=\Cmath$, the lemma is not true for Hilbert modules $\cH$ of infinite dimension.
\end{rema}
\begin{proof}
Denote by $\cM'$ the commutant of $\cM$ in $B(\cH)$. Since the module $\cH$ has finite dimension also $\cM'$ carries a finite normal faithful trace, denoted $\tau'$. Write $V$ as the union of an increasing net of Hilbert modules $V_n$. Let $p_n \in M_N(\cM')$ be the orthogonal projection on $V_n$ and set $p:=\lim_n p_n$, the projection on $\overline{V}$. We can do the same for $W$ to obtain $W_n,q_n$ and $q$. The orthogonal projection on $\overline{V} \cap \overline{W}$ is $ p\wedge q$. Similarly, since $\overline{V \cap  W}$ is the closure of the increasing net of closed modules $V_n \cap  W_n$, the orthogonal projection on it is $\lim_n p_n \wedge q_n$. So we have to prove the equality $p \wedge q = \lim_n p_n \wedge q_n$.

By \cite[Proposition V.1.6]{takesakiVol1} $\tau'( p \wedge q) = \tau'(p+q-p\vee q)$ and same for $p_n$, $q_n$, so we have 
\[ \tau'( p \wedge q - p_n \wedge q_n) = \tau'(p-p_n) +\tau'(q-q_n) - \tau'(p \vee q - p_n \vee q_n).\] 
Both terms $\tau'(p-p_n)$ and $\tau'(q-q_n)$ go to zero, whereas $\tau'(p \vee q - p_n \vee q_n) \geq 0$, so we get $\limsup_n \tau'( p \wedge q - p_n \wedge q_n)\leq 0$. On the other hand $p_n \wedge q_n \leq p \wedge q$ and therefore we obtain that $p \wedge q =\lim_n p_n \wedge q_n$ as required.
\end{proof}

 It is important for us that we have formulas in terms of ``dimensions of kernels'' for applying again Lemma~\ref{lem:lueck_approximation}.

\begin{lemm}\label{lem:formula_for_nabla} Let $\Sigma \subset \Sigma'$ be two uniformly locally bounded  $\cg$-simplicial complexes, with associated boundary maps $\partial'_i \colon C^{(2)}_i(\Sigma') \to C^{(2)}_{i-1}(\Sigma')$ and $\partial_i \colon C^{(2)}_i(\Sigma) \to C^{(2)}_{i-1}(\Sigma)$, and denote by $p_i$ the orthogonal projection onto $C^{(2)}_i(\Sigma)\subset C^{(2)}_i(\Sigma')$. Then
\[ \nabla_i(\Sigma,\Sigma') = \mathrm{dim}(\ker \partial_i) + \mathrm{dim}(\ker \partial'_{i+1}) - \mathrm{dim}(\ker((1-p_i)\partial'_{i+1})),\]
where $\dim$ stands for the von Neumann dimension over $\mathcal L(\cg)$.
\end{lemm}
\begin{proof}
Let $V$ denote the Hilbert $\mathcal L(\cg)$-module $\ker((1-p_i)\partial'_{i+1})$ and observe that if $f\in V$, then $\partial'_{i+1}f=p_i\partial'_{i+1}f$ and hence \[\partial_i\partial'_{i+1}f=\partial_i p_i\partial'_{i+1}f=\partial'_ip_i\partial'_{i+1}f=\partial'_i\partial'_{i+1}f=0,\] since $\partial_i$ and $\partial_i'$ coincide on the subspace generated by the image of $p_i$. Let $\kappa:\ker \partial_i\rightarrow \overline H_i^{(2)}(\Sigma')$ be the quotient map by $\overline{\Ima \partial'_{i+1}} \cap \ker \partial_i$. Consider the sequence of Hilbert $\mathcal L(\RR)$-modules
\[ 0 \to \ker \partial'_{i+1} \xrightarrow{\iota} V \xrightarrow{\partial'_{i+1}} \ker \partial_i \xrightarrow{\kappa}  W \to 0,\]
where $\iota$ is just the inclusion. We claim that it is weakly exact. The only point which deserves a justification is that the closure of the image of $\partial'_{i+1}$, namely $\overline{\Ima \partial'_{i+1} \cap \ker \partial_i}$ agrees with the kernel of $\kappa$, that is $\overline{\Ima \partial'_{i+1}} \cap \ker \partial_i$ which we already established in the previous lemma.

 The lemma follows by the additivity of the von Neumann dimension.
\end{proof}

\begin{proof}[Proof of Theorem~\ref{thm:limitgeneral}]
  Let us denote by $\partial_{n,L}^i$ the boundary map from $C_{i}^{(2)}(\Sigma_{n,L})$ to $C_{i-1}^{(2)}(\Sigma_{n,L})$. Since $\partial_{n,L}^i$ restricted to every orbit $\Sigma_{n,L}[x]$ is a classical boundary operator it can be expressed as a matrix with integer coefficients (with coefficients in $[[\cg_n]]$ if $\cg_n$ is only assumed to be sofic). By Lemma~\ref{lem:ultraproduct_of_pseudo_full_groups} and Lemma~\ref{lem:ultraproductactions} we have that $\partial_{\ul,L}^i=[\partial_{n,L}^i]_\ul$ is the boundary map from $C_{i}^{(2)}(\Sigma_{\ul,L})$ to $C_{i-1}^{(2)}(\Sigma_{\ul,L})$. We also clearly have that $[(\partial_{n,L}^i)^*\partial_{n,L}^i]_\ul=(\partial_{\ul,L}^i)^*\partial_{\ul,L}^i$ and therefore we can conclude from Lemma~\ref{lem:lueck_approximation} (or Lemma~\ref{lem:lueck_approximation_sofic}) that \[ \dim\ker(\partial_{\ul,L}^i) = \lim_n  \dim\ker(\partial_{n,L}^i).\]
Using a similar argument to all the operators involved in the expression of $ \nabla_i(\Sigma,\Sigma')$ in Lemma~\ref{lem:formula_for_nabla} one obtains the desired formula.
\end{proof}

\begin{rema}[$\ell^2$-Betti numbers for random unimodular ULB rooted simplicial complex]
\label{rem: random unimodular simplicial complexes}
If $\Theta$ is random unimodular ULB rooted simplicial complex with distribution $\nu$, we consider the Laplace operator $\Delta_{d,\Theta}$ on the space $C^{(2)}_d(\Theta)$ of $L^2$-chains for each dimension $d$ and the orthogonal projection $p_{d,\Theta}: C^{(2)}_d(\Theta)\to \ker \Delta_{d,\Theta}$. 
We consider 
$\frac{1}{d+1} \sum_{\sigma\sim \rho_{\Theta}} \langle p_{d,\Theta}(\indic_{\sigma}) , \indic_{\sigma}\rangle=\frac{1}{d+1} \sum_{\sigma\sim \rho_{\Theta}} \Vert p_{d,\Theta}(\indic_{\sigma}) \Vert^2$
where the sum is extended to all $d$-cells $\sigma$ that are adjacent to the root $\rho_\Theta$ and integrate this quantity with respect to $\nu$
\begin{equation}\beta_d^{(2)}(\nu):= \int \frac{1}{d+1} \sum_{\sigma\sim \rho_{\Theta}} \langle p_{d,\Theta}(\indic_{\sigma}) , \indic_{\sigma}\rangle d\nu(\Theta).\end{equation}
Compare \cite{Schrodel-L2-Betti-simpl-cplx-2018}. See also \cite{Elek-10-Betti-test}.
If for instance $\nu$ is the Dirac measure at the Cayley complex of some finitely presented group, one recovers the first $\ell^2$-Betti number of the group as $\beta_1^{(2)}(\nu)$. Compare with the formulas from \cite{Gab05}.

If for instance $\cg$ is a \pmp groupoid on $(X,\mu)$ and $\Sigma$ is a ULB $\cg$-simplicial complex, then pushing forward the measure $\mu$ by the map $\pi:x\mapsto \Sigma_x$ defines a random unimodular ULB rooted simplicial complex $\pi_*\mu$ and the change of variable formula immediatly gives:
\[\beta^{(2)}_d(\Sigma, \cg)= \beta^{(2)}_d(\pi_*\mu).\]

Just like \pmp groupoids admit $\ell^2$-Betti numbers when they act on $\cg$-simplicial complexes, but they also admits their own  $\ell^2$-Betti numbers, there is also an absolute version of $\ell^2$-Betti numbers for random unimodular ULB rooted simplicial complexes. They are defined by taking the $L$-Rips complexes $R^L(\Theta)$ of the random complex $\Theta$, the closure $\nabla_d(R^{L}(\Theta),R^{L'}(\Theta))$ of the image  of the map (induced  by inclusion) between the kernels of the Laplace operators for the various values of $L$ and the absolute $\ell^2$-Betti numbers:
\begin{equation}\label{eq:formula_for_betti_number random complexes} 
\widehat{\beta_d^{(2)}(\nu)}:=\lim_{L} \lim_{L' \geq L} \nabla_d(R^{L}(\Theta),R^{L'}(\Theta)).
\end{equation}

If  the random complexe $\Theta$ is for instance $k$-connected almost surely, then 
$\widehat{\beta_d^{(2)}(\nu)}=\beta_d^{(2)}(\nu)$ for all $d\leq k$.
When $\theta$ is produced as above by a ULB $\cg$-simplicial complex, then 
\[\beta^{(2)}_d(\cg)= \widehat{\beta^{(2)}_d(\pi_*\mu)}.\]
In particular, if $(\cg,\Phi)$ and $(\ch,\Psi)$ are two finite size graphed \pmp groupoids that produce the same random unimodular (unlabelled) network, then  for every~$d$: \[\beta^{(2)}_d(\cg)=\beta^{(2)}_d(\ch).\]

\end{rema}

\newcommand{\etalchar}[1]{$^{#1}$}
\def\cprime{$'$} \def\cprime{$'$}

\end{document}